\documentclass[reqno, a4paper, 11pt, oneside]{article}
\usepackage[11pt]{extsizes}
\usepackage[a4paper,hmargin=2.5cm,vmargin=2.5cm, bindingoffset=0cm]{geometry}

\usepackage{mathtools}
\usepackage[utf8]{inputenc}
\usepackage[T2A]{fontenc}
\usepackage{amsmath,amssymb,amsthm}
\usepackage[english]{babel}
\usepackage{textcomp}
\usepackage[usenames]{color}
\usepackage{colortbl}
\usepackage{float}
\usepackage{relsize}
\usepackage{amsthm}

\usepackage{thmtools}
\usepackage{thm-restate}
\makeatletter
\newcommand*{\rom}[1]{\expandafter\@slowromancap\romannumeral #1@}
\makeatother

\usepackage[pagebackref,bookmarksopen]{hyperref}
\hypersetup{
	colorlinks=true,
	linkcolor=blue,
	citecolor=magenta,
	urlcolor=cyan,
}
\usepackage{cleveref}

\usepackage{graphicx}
\DeclareGraphicsExtensions{.pdf,.png,.jpg}

\usepackage{amsthm}

\theoremstyle{definition}

\newtheorem{theorem}{Theorem}[section]

\newtheorem*{remark}{Remark}
\newtheorem{claim}[theorem]{Claim}
\newtheorem{lemma}[theorem]{Lemma}
\newtheorem{cor}[theorem]{Corollary}
\newtheorem{sign}[theorem]{Notation}
\newtheorem{con}[theorem]{Conjecture}

\newtheorem{ex}[theorem]{Example}
\newtheorem*{exa}{Example \ref{B_2-facet_example}}
\newtheorem{definition}[theorem]{Definition}

\newtheorem{problem}{Problem}
\newtheorem*{theorem-non}{Theorem}

\definecolor{applegreen}{rgb}{0.0, 0.42, 0.24}
\definecolor{applegreen}{rgb}{0.55, 0.71, 0.0}

\newcommand{\tl}{\widetilde}

\newcommand{\Sset}{\mathcal F}

\newcommand{\Z}{\mathbb{Z}}
\newcommand{\R}{\mathbb{R}}

\renewcommand{\C}{\mathbb{C}}

\newcommand{\D}{\Delta}

\usepackage{abstract}
\addto{\captionsenglish}{}

\begin{document}
	
\begin{center}
{\Large \sc
Arnold's monotonicity problem}
\vspace{3ex}

{\sc Fedor Selyanin} 

\end{center}
	\begin{abstract}
        According to the Kouchnirenko formula, the Milnor number of a generic isolated singularity with given Newton polyhedron is equal to the alternating sum of certain volumes associated to the Newton polyhedron. In this paper we obtain a non-negative analogue (i.e. without negative summands) of the Kouchnirenko formula. The analogue relies on the non-negative formula for the monodromy operator from \cite{Stap17} and formulas for the Milnor number from \cite{Fur04}. As an application we give a criterion for the Arnold's monotonicity problem (1982-16) in arbitrary dimension, which leads to complete solution in dimension up to $4$ and partial solution in dimension $5$. The latter relies on the classification of thin triangulations (or vanishing local h-polynomial) in dimension $2$ and $3$ from \cite{dMGP+20} and \cite{GKZ94} and contains examples which differ dramatically from the ones which arise in dimension up to $3$ in \cite{BKW19} (see also \cite{L-AMS20}). Some of the $4$-dimensional examples were first described in \cite{ELT22} in the context of the local monodromy conjecture.


	\end{abstract}

\setcounter{tocdepth}{2}
\tableofcontents

\section{Introduction}
     A germ of analytic function $f: (\C^n,0) \to (\C^1,0)$ is called a \emph{singularity} if the gradient $(\frac{\partial f}{\partial x_1} , \dots , \frac{\partial f}{\partial x_n})$ vanishes at the origin. Singularity is called \emph{isolated} if there are no other singular points in a neighbourhood of the origin. \emph{Milnor number} $\mu (f)$ of an isolated singularity $f: (\C^n,0) \to (\C^1,0)$ is the number of Morse singular points in the neighbourhood of the origin of a generic small perturbation $\tl f$ of $f$ (e.g. see \cite{AGLV}). See \cite{Mil68} for the original definition of the Milnor number.

     The \emph{support} of the analytic function $f$ is the subset in $\mathbb Z_{\ge 0}^n$ consisting of the exponents of all monomials that have non-zero coefficients in the power series of $f(x)$, and the \emph{Newton polyhedron} $\Gamma_+ (f)$ of $f$ is the convex hull of the union of the positive octants with vertices in the points of the support of $f$. For a Newton polyhedron $\Gamma_+$ denote by $\Gamma$ the Newton diagram (the union of its compact faces) and by $\Gamma_-$ the complement $\R^n_{\ge 0} \setminus \Gamma_+$. The Newton polyhedron $\Gamma_+$ is called \emph{convenient} if it contains a point on each coordinate axis. The following theorem describes the Milnor number of a singularity in terms of its Newton polyhedron in non-degenerate case (for a comprehensible survey on the connections between algebraic geometry and convex geometry see \cite{EKK}).
	
    \begin{theorem-non}[Kouchnirenko \cite{Kou76}]\label{th_kush}
		If the Newton polyhedron $\Gamma_+$ is convenient, then generic singularity $f$ with Newton polyhedron $\Gamma_+$ is isolated and the Milnor number $\mu(f)$  is equal to the Newton number
		\begin{equation}\label{KF}
		\nu (\Gamma_-) = V - \sum\limits_{1\le i \le n} V_i + \sum\limits_{1\le i < j \le n} V_{i,j} - \dots + (-1)^n
		\end{equation}
		of the polytope $\Gamma_-$, where $V$ is the $n$-dimensional lattice volume
  \footnote{We use the lattice volume $Vol(\cdot)$ which means that the volume of the unit lattice cube in $\mathbb R^n$ equals $n!$.}
of the polytope $\Gamma_-$, $V_i$ is the $(n-1)$-dimensional lattice volume of the polytope 
		$\Gamma_- \cap \{x_i = 0\}$, $V_{i,j}$ is the $(n-2)$-dimensional lattice volume of the polytope $\Gamma_- \cap \{x_i = x_j = 0\}$, and so on.
	\end{theorem-non}

The Newton number is monotonous (i.e. $\Gamma_+ \subset \Gamma_+^\prime$ implies $\nu(\Gamma_-) \ge \nu(\Gamma_-^\prime)$) because of monotonicity of the Milnor number under small perturbations. But it is hard to deduce this monotonicity from the combinatorics of the Kouchnirenko formula because it contains negative summands. In his list of problems, V. I. Arnold posed the following one (\cite{Ar}, № 1982-16):

\begin{quote}
     ``Consider a Newton polyhedron $\Gamma_+$ in $\R^n$ and the number $\nu(\Gamma_-) = V -\sum\limits_{1\le i \le n} V_i + +\sum\limits_{1\le i < j \le n} V_{ij} - \dots$\footnote{We write the formulas using our notations to avoid confusion.}, where $V$ is the lattice volume under $\Gamma_+$, $V_i$ is the lattice volume under $\Gamma_+$ on the hyperplane $x_i = 0$, $V_{ij}$ is the lattice volume under $\Gamma_+$ on the hyperplane $x_i = x_j = 0$, and so on.
Then $\nu(\Gamma_-)$ grows (non strictly monotonically) as $\Gamma_-$ grows (whenever $\Gamma_+$ remains convex
and integer?). There is no elementary proof even for $n = 2$.''
\end{quote}

In this paper we give a non-negative analogue of the Kouchnirenko formula which is useful for solving problems on the monotonicity of the Newton number. For example, we do not know a combinatorial proof even for the following problem without using the non-negative analogue.

\begin{problem}\label{Naive_problem}
    Prove that $\nu(\Gamma_-) \ge 0$ for any convenient Newton polyhedron $\Gamma_+$.
\end{problem}

The non-negative analogue of the Kouchnirenko formula automatically solves the Arnold's problem (1982-12) in its original form. But related problems, such as the monodromy conjecture, motivate a refined version of this question: how to classify extensions of a given Newton polyhedron preserving the Newton number? So from now on we call the following problem the ``Arnold's monotonicity problem''. Denote by $\Gamma_{+G}$ the polyhedron $Conv(\Gamma_+, G)$.

\begin{problem}[Arnold's monotonicity problem]\label{Arnold_problem}
    Consider a convenient Newton polyhedron $\Gamma_+ \subset \R^n_{\ge 0}$ and a lattice point $G\in \Gamma_-$. Give a criterion for $$\nu(\Gamma_-) = \nu(\Gamma_{-G})$$
\end{problem}

The Arnold's monotonicity problem is solved in dimension up to $3$ in \cite{BKW19} (see also \cite{O89} and \cite{Fur04} if $\Gamma_{+G} \setminus \Gamma_+ \subset \R^n_{\ge 0}$ is a simplex for $n = 3$ and $n=4$ respectively) and has very simple answer. Namely, the Newton number is preserved if and only if the difference $\Gamma_{+G} \setminus \Gamma_+$ is a $B_1$-pyramid, i.e. a pyramid of height $1$ with base on a coordinate hyperplane. It is conjectured in \cite{BKW19} and claimed in \cite[Theorem 2.25 and Remark 2.22]{L-AMS20} that the same criterion works in arbitrary dimension (we call it the ``\emph{$B_1$-conjecture}''). But most of the examples in dimension $4$ do not satisfy this conjecture. The authors of \cite{L-AMS20} however have informed the author that they believe that their approach to this problem is working as well, and the gap will be fixed.

\begin{ex}(\cite[$B_2$-facet]{ELT22})\label{B_2-facet_example}
Let us give the simplest example of the Morse singularity in dimension $4$ which contradicts the $B_1$-conjecture from \cite{BKW19}. Note that there are similar examples of not Morse singularities. Consider the Newton polyhedron $\Gamma_+$ spanned by the points from the table below, except for the point $G$.

\begin{table}[H]
    \centering
    \begin{tabular}{c|c|c|c|c|c|c|c|c}
        C & B & P & Q & X & Y & Z & T & G  \\
         (1,0,1,0) & (1,0,0,1) & (0,1,1,0) & (0,1,0,1) & (0,0,5,0) & (0,0,0,5) & (3,0,0,0) & (0,3,0,0) & (0,0,2,2)
    \end{tabular}
    \label{B_2-table}
\end{table}

\begin{figure}[H]
		\begin{center}
		\includegraphics[scale=0.25]{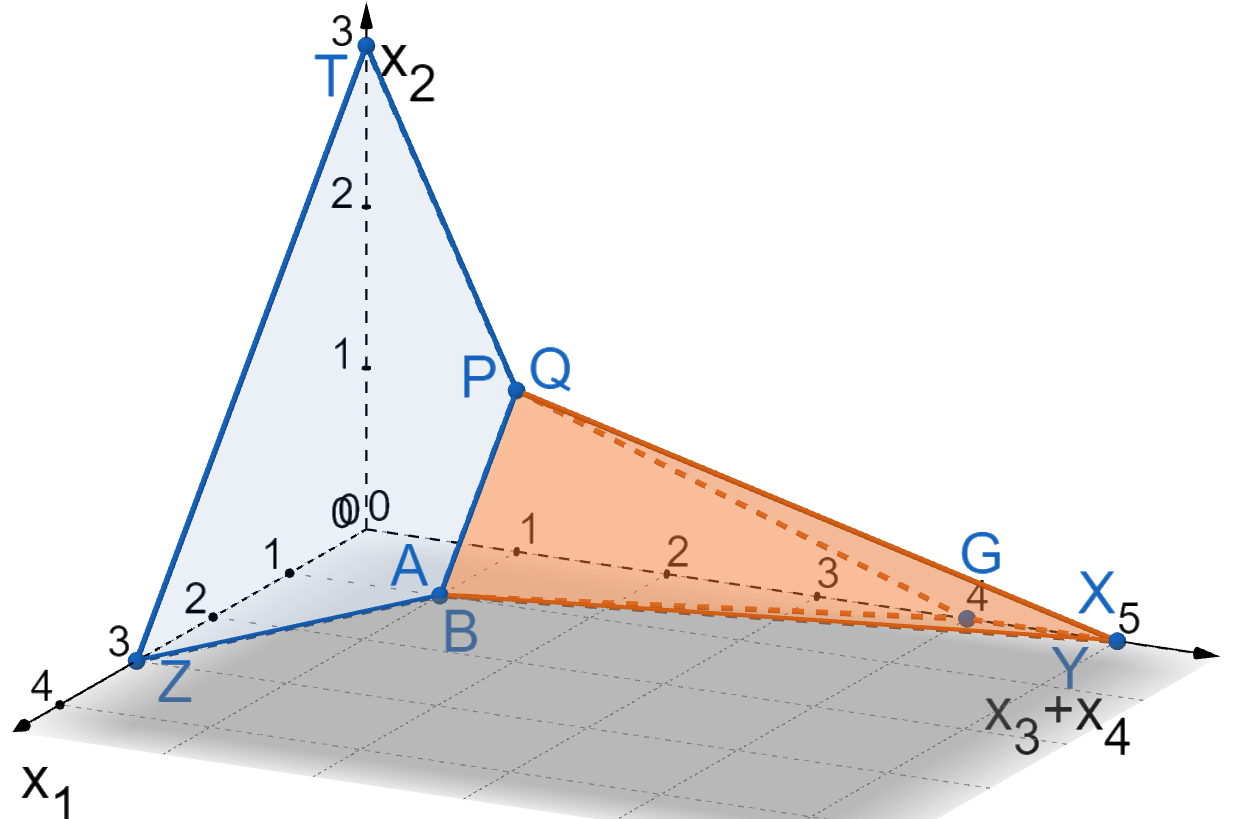}
		\includegraphics[scale=0.25]{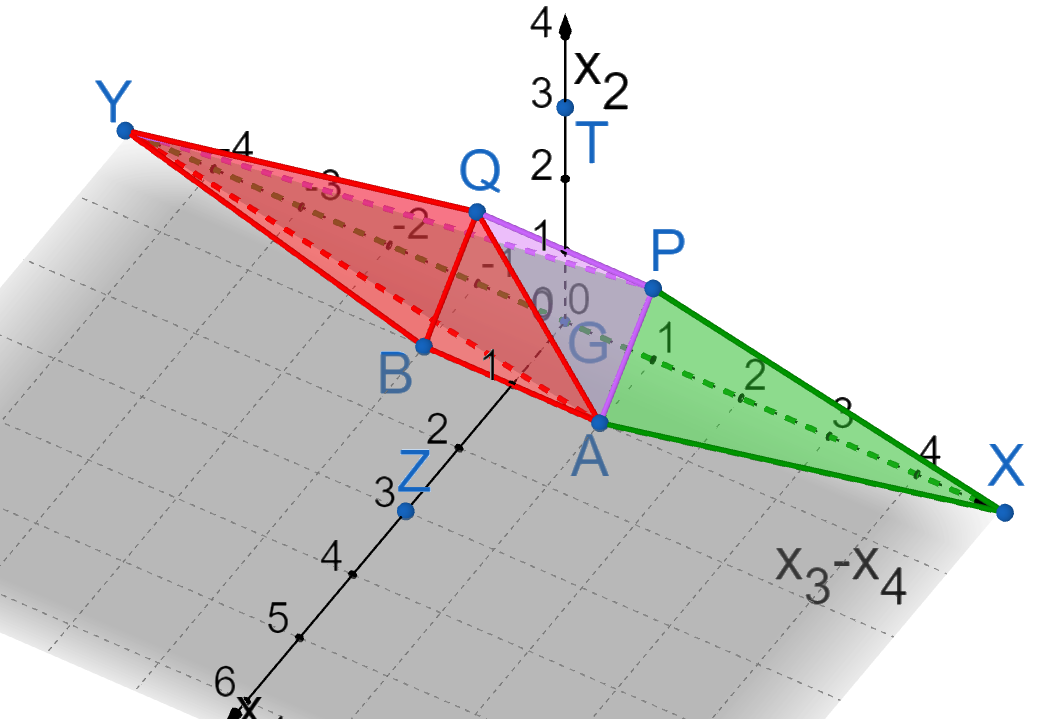}
		\end{center}
		\caption{\label{B_2-picture} Two projections of $\Gamma_+$ and $G$: left $(x_1,x_2,x_3+x_4)$, right $(x_1,x_2, x_3-x_4)$. 
  \\Left: The projection of the polytope $\Gamma_{+G} \setminus \Gamma_+$ is colored brown. 
  \\Right: The polytope is later denoted by $\Pi_\mathsmaller{G}(\Gamma_+)$. Its splitting into $3$ simplices is needed for the continuation of this example on p.\pageref{continued_example}.}
\end{figure}
The Newton polyhedron $\Gamma_+$ defines the Morse singularity if the coefficients are non-degenerate (for example, since  $(1/2,1/2,1/2,1/2)\in \Gamma$, see \cite{Y19}), therefore $$\nu(\Gamma_-) = \nu(\Gamma_{-G} )= 1$$

Note that the set $\Gamma_{+G}\setminus \Gamma_+$ contains no lattice points except for $G$, so there is surely no way to extend the $B_1$-conjecture to this case.
\end{ex}

\begin{remark}
First examples of such kind (called $B_2$-facets) were discovered in \cite{ELT22} in the context of the local monodromy conjecture in dimension $4$. The local monodromy conjecture (defined in \cite[Conjecture 3.3.2]{DL92}; see \cite{Ve} for a survey) is a mysterious open problem which relates arithmetic and topological/geometric properties of a polynomial $f$. The case of interest is when the zero set of $f$ has singular points. A good class of examples of polynomials to verify the local monodromy conjecture are non-degenerate polynomials with respect to their Newton polyhedra and with isolated singularity at the origin. The local monodromy conjecture is confirmed for simplicial Newton polyhedra in \cite{LPS22a} and for arbitrary Newton polyhedra in dimension up to $4$ in \cite{ELT22}. The local monodromy conjecture is still not verified for arbitrary Newton polyhedra and is closely related to the Arnold's monotonicity problem.
\end{remark}

The core results of the paper are the following ones (for more details see \S \ref{main_results_section}).
\begin{itemize}
    \item Theorem \ref{Non-negative_analogue_theorem}: Non-negative analogs of the Kouchnirenko formula.
    \item Theorem \ref{Arnold_monotonicity_criterion_theorem}: Criterion for the Arnold's monotonicity problem in arbitrary dimension.
    \item Theorem \ref{4_5_Arnold_theorem}: Complete solution of the Arnold's monotonicity problem in dimension $4$ and partial in dimension $5$.
\end{itemize}

\subsection{Arnold's monotonicity problem in arbitrary dimension}

In this section we give a criterion for the Arnold's monotonicity problem \ref{Arnold_problem} in arbitrary dimension $n$. Consider the Arnold's monotonicity problem with convenient a Newton polyhedron $\Gamma_+$ and a point $G$. Denote by $\pi_{\mathsmaller{OG}}$ the projection along the line $OG$. Denote by $\Pi_G(\Gamma)$ the polytope 
$$\Pi_G(\Gamma) = \pi_{\mathsmaller{OG}}(\Gamma_{+G} \setminus \Gamma_+) \subset \pi_{\mathsmaller{OG}}(\R^n_{\ge 0})$$
The polytope $\Pi_G(\Gamma)$ is not necessarily coconvex but it is a \emph{complete star-polytope}, i.e. it admits a triangulation (provided by the Definition \ref{bottom_triang_def}  of \emph{bottom triangulation}) all whose maximal simplices contain the origin  and the cone over it is the whole $\pi_{\mathsmaller{OG}}(\R^n_{\ge 0})$. Note that 
$$\pi_{\mathsmaller{OG}}(\R^n_{\ge 0}) = \R^{m}_{\ge 0} \oplus \R^{n-m-1},$$
where $m$ is the codimension of the minimal face of the positive octant $\R^n_{\ge 0}$ containing the point $G$. The Newton number can be defined for pairs of arbitrary cone and polytope inside the cone as the alternating sum of lattice volumes (see Definition \ref{Arbutrary_cone_Newton_number_def}). If the Newton number of a polytope inside a cone is equal to $0$ then the polytope is called \emph{negligible} in the cone (see Definition \ref{star-negligible_def}). \\

\textbf{Theorem \ref{Arnold_monotonicity_criterion_theorem}(1)} (Negligible criterion). A pair $(\Gamma_+,G) \subset \R^n_{\ge 0}$ satisfies the Arnold's monotonicity problem if and only if the complete star-polytope $\Pi_G(\Gamma) = \pi_{\mathsmaller{OG}}(\Gamma_{+G} \setminus \Gamma_+)$ is negligible in the cone $C_G = \pi_{\mathsmaller{OG}}(\R^n_{\ge 0})$.\\

Thus the Arnold's monotonicity problem in dimension $n$ is reduced to the classification of negligible complete star-polytopes in the cones $\R^m_{\ge 0} \oplus \R^{n-m-1}$. The main tool for the classification of negligible complete star-polytopes is called the \emph{Thin criterion} (see Theorem \ref{Arnold_monotonicity_criterion_theorem}(2)).

\subsection{Arnold's monotonicity problem in dimension 4} \label{4-Arnold_intro_sec}

One of the main results of this paper is the complete solution of the Arnold's monotonicity problem in dimension $4$. Here we introduce a naive and demonstrative way of describing the solutions which can hardly be generalized for higher dimension (for the regular description of the result using \textit{conical facet refinements} see Theorem \ref{4_5_Arnold_theorem}(1)).

Recall that the Arnold's monotonicity problem in dimension $n$ is reduced to the classification of negligible complete star-polytopes in the cones $\R^{m}_{\ge 0} \oplus \R^{n-m-1}$. Let us describe the classification of negligible complete star-polytopes in dimension $3$ for all pairs $(m,3-m)$ up to reordering of the coordinates.

\begin{enumerate}
    \item[\textbf{0,3:}] There are no negligible complete star-polytopes in $\R^3$.
    \item[\textbf{1,2:}] All negligible complete star-polytopes in $\R^1_{\ge 0} \oplus \R^2$ are $B_1$-pyramids i.e. pyramids of height $1$ with the base on $\R^2$.

 \begin{figure}[H]
		\begin{center}
		\includegraphics[scale=0.3]{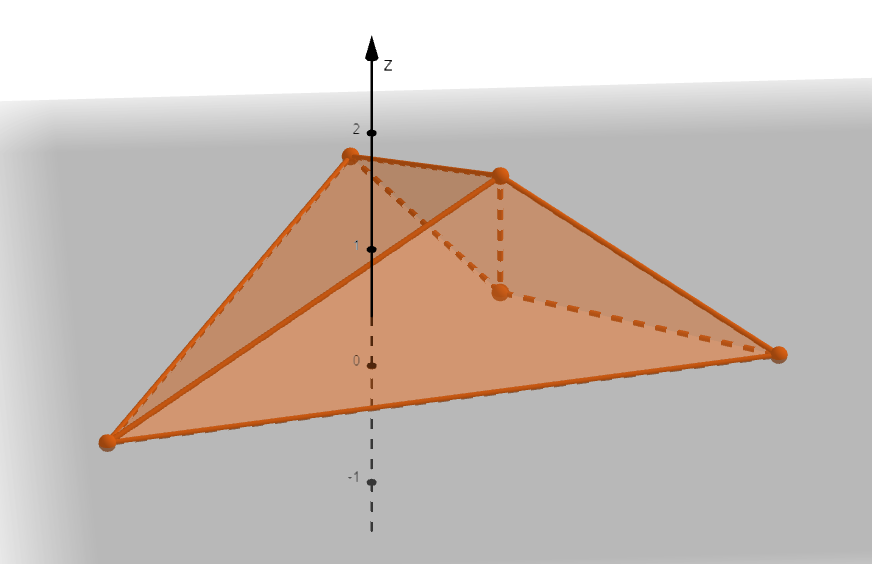}
		\end{center}
		\caption{Example of \textbf{$\mathbf{B_1}$-pyramid}}
\end{figure}
    
    \item[\textbf{2,1:}] All negligible complete star-polytopes in $\R^2_{\ge 0} \oplus \R^1$ are \emph{greenhouses}.
    
    Consider two lines $L_1 = (1,0,\star), L_2 = (0,1,\star)$ on the boundary of $\R^2_{\ge 0} \oplus \R^1$. \emph{Greenhouse} is a polytope bounded by a diagram which can be split to parts, where each part is a pyramid with an apex on $L_1$ (or $L_2$) and a base on the opposite $2$-dimensional face of $\R^2_{\ge 0} \oplus \R^1$.
    
 \begin{figure}[H]
		\begin{center}
		\includegraphics[scale=0.2]{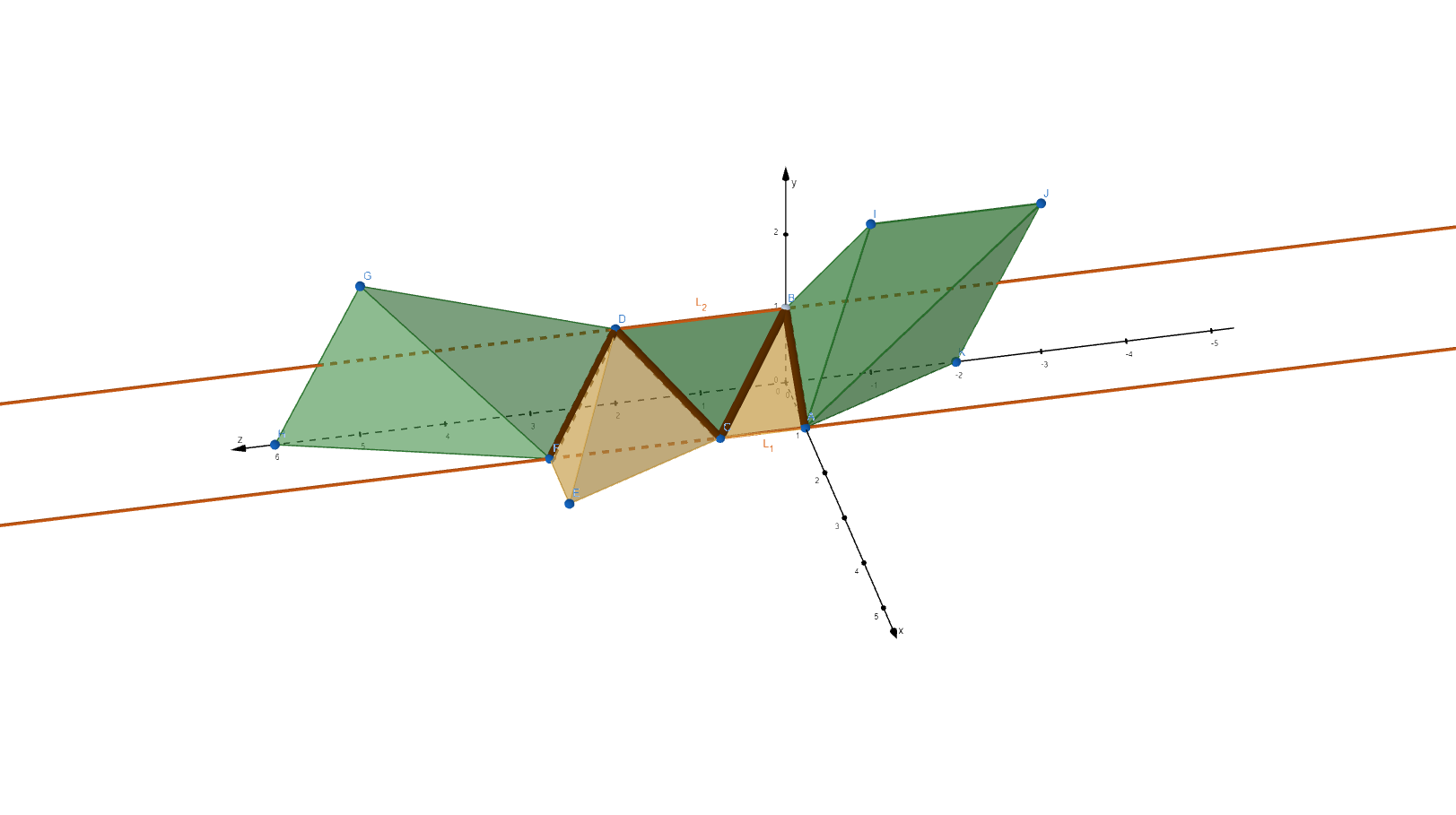}
		\end{center}
		\caption{\label{Greenhouse} Example of \textbf{greenhouse}}
\end{figure}
    
    \item[\textbf{3,0:}]  All negligible complete star-polytopes in $\R^3_{\ge 0}$ are at most two \emph{greenhouses} connected with a \emph{crook}.

    Consider the point $X = (0,0,1)$ and a diagram on the plane $x_3 = 0$ with endpoints $Y,Z$ distant from the corresponding coordinate rays by at most  $1$. Consider the pyramid $D$ with the apex $X$ and the base on the diagram. \emph{Crook} is the pyramid with the base $D$ and the apex in the origin. If the diagram's end $Y$ (or $Z$) is not on the coordinate ray, then there is also a greenhouse on the corresponding ray glued by the triangle $OXY$ (or $OXZ$) to the crook.
     \begin{figure}[H]
		\begin{center}
		\includegraphics[scale=0.4]{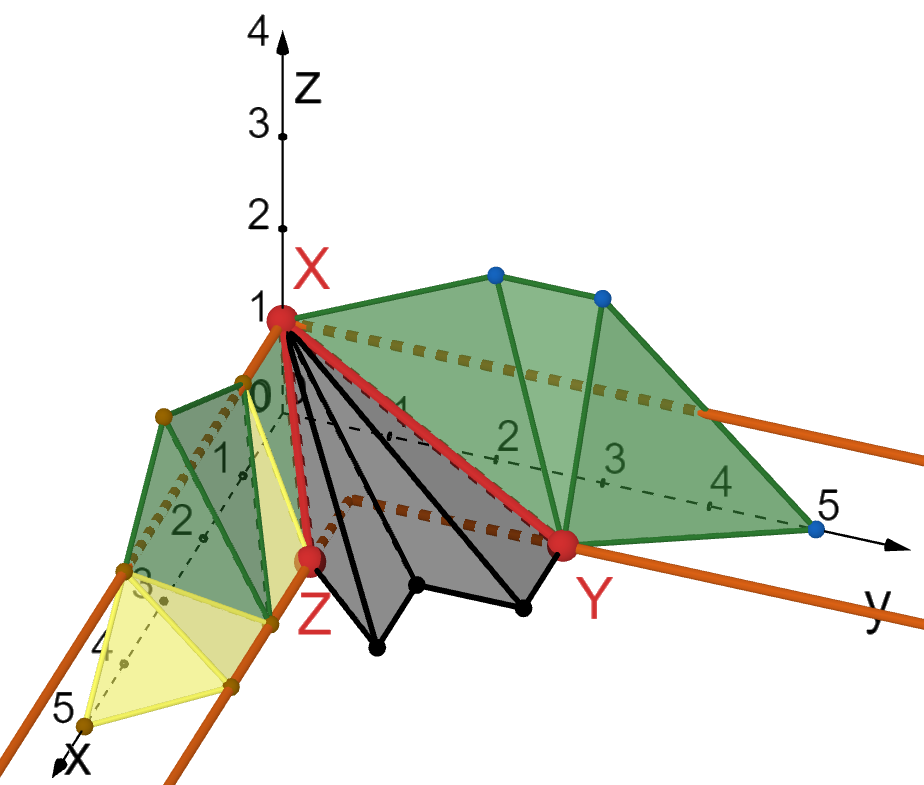}
		\end{center}
		\caption{\label{Crook} Example of a \textbf{crook} (black) connecting two \textbf{greenhouses}}
\end{figure}

\end{enumerate}

Note that the answers in the cases \textbf{3,0} and \textbf{2,1} do not satisfy the $B_1$-conjecture from \cite{BKW19}. The simplest Example \ref{B_2-facet_example} \textbf{contradicting} the $B_1$-conjecture is a \textbf{greenhouse} (the case  \textbf{2,1}). If the polytope $\Gamma_{+G} \setminus \Gamma_+$ in the case \textbf{2,1} contains only one facet of the polyhedron $\Gamma_+$ then the facet is the \textit{$B_2$-facet} from \cite{ELT22}. If the crook in the case \textbf{3,0} consists of only one segment (i.e. $Y=Z = (1,1,0)$) then the preimage of the triangle $OXY = OXZ$ under the projection $\pi_{\mathsmaller{OG}}$ is the \textit{$B_2$-border} from \cite{ELT22}.

Note that greenhouses from the case \textbf{2,1} are parts of the answer in the case \textbf{3,0}. A more detailed consideration of this phenomena in arbitrary dimension is given in \S \ref{Sprigs_subsecion}. It plays the crucial role in the solution of the Arnold's monotonicity problem in dimension $4$ and $5$.

\subsection{Acknowledgment}

I am grateful to Alexander Esterov for supervising me, for the suggestions to apply my considerations to the Arnold's monotonicity problem and to pay more attention to the Ehrhart theory. The initial considerations were devoted to generalizing the results of \cite{Vas17} to the case of Newton non-degenerate isolated singularities with arbitrary Newton polyhedron.

\section{Preliminaries}

\subsection{Polytopes and triangulations}

Let us mention the main definitions and notations concerning the polyhedral geometry used in this paper. For basic definitions see e.g. \cite{Ew}. Recall that lattice convex polytope in $\R^n$ is a convex polytope with vertices in $\Z^n \subset \R^n$.

\begin{definition}\label{simplicial_complex_def}
\emph{Pure lattice geometric simplicial complex} $T$ in $\R^n$ is a simplicial complex consisting of lattice simplices such that the dimension of all its maximal simplices equals $n$. Note that the simplicial complex forms a partially ordered set with respect to inclusion.
\end{definition}

\begin{definition}\label{link_def}
    The link $\text{lk}_T (\D)$ of a simplex $\D\in T$ is the pure lattice geometric simplicial complex spanned by the simplices $\D^\prime \in T: \D^\prime > \D$. 
\end{definition}

\begin{definition}\label{polytope_def}
\emph{Polytope} $P\subset \R^n$ is the union of all simplices of a pure lattice geometric simplicial complex $T$. We call pure lattice geometric simplicial complex $T$ a \emph{triangulation} of the polytope $P$ if the union of all simplices of $T$ is the polytope $P$. Polytope $P^\prime \subset P$ is called a subpolytope of the polytope $P$ if there is a triangulation of $P$ such that the union of some of its simplices is $P^\prime$.
\end{definition}

\begin{definition}\label{regular_triangulation_def}
    \emph{Regular (or coherent)} triangulation of a convex polytope $P\subset \R^n$ with a support $\mathbf P \subset \Z^n: Conv(\mathbf P) = P$ is a triangulation induced by a height function $\psi: \mathbf P \to \R^1$, i.e. it consists of projections of the bounded simplices of the convex hull of the rays $(p,t), p \in \mathbf P, t \ge \psi(p)$.
\end{definition}

\begin{sign}\label{projection_sign}
    Consider a set $X\subset \R^n$. Denote by $\pi_\mathsmaller{X}$ the affine projection $\R^n\to \R^n/\text{aff}(X)$ with the kernel $\text{aff}(X) \to 0$.
\end{sign}

\begin{definition}\label{localized_triangulation_def}
    Suppose that $T$ is a pure lattice geometric simplicial complex and $\D \in T$ is a simplex. The \emph{localized triangulation} $T/\D$ is the triangulation obtained as the projection $\pi_\D(\text{lk}_T(\D))$ of the link lk$_T(\D)$.
\end{definition}

\subsection{Local $h$- and $h^*$-polynomials}\label{Ehrhart_theory_section}

In this subsection we recall the local $h$- and $h^*$-polynomials coming from the Ehrhart theory. In the paper we use only the specialization of the locale $h$-polynomial for the cones $\R^m_{\ge 0} \oplus \R^{n-m}$ and the local $h^*$-polynomial of a simplex (it is also called the \emph{box polynomial}). These specializations can also be described in terms of the Newton number and the combinatorial Newton number from \cite{GKZ94} and \cite{GKZ90} which we recall in \S \ref{Newton_section}. But we believe that for the genegalizations of the results of this paper for arbitrary cones the language of the local $h$- and $h^*$-polynomials is more appropriate then the one coming from \cite{GKZ94} and \cite{GKZ90}.

A good elementary introduction to the Ehrhart theory is given in \cite{BR}, but it does not contain the definitions of the local $h$- and $h^*$-polynomials (though it contains the \emph{box polynomial} which is $h^*$-polynomial of the simplex). An extremely intelligible survey of the local $h$- and $h^*$-polynomials for lattice polytopes and lattice polyhedral subdivisions can be found in the introduction of \cite{BKN23}. A survey to the local $h$-polynomials of (not necessarily geometric) simplicial subdivisions of the simplex is given in \cite{Ath}.

\subsubsection{Local $h$- and $h^*$-polynomials in the Ehrhart theory}

The $h$-polynomial encodes the number of faces of different dimension of simplicial complex such that the Dehn-Sommerville equations between the numbers of faces of simplicial convex polytopes are equivalent to the palindromicity of the $h$-polynomial. The $h$-polynomial was significantly generalized by Stanley \cite{Stan87} to the case of so called Eulerian posets. For instance, faces of any convex polytope form an Eulerian poset and the coefficients of the corresponding $h$-polynomial are the dimensions of intersection cohomology groups of the projective toric variety corresponding to the polytope (see e.g. \cite{Ful} for a book on toric varieties). The $h^*$-polynomial of a lattice polytope $P$ is a good way to rewrite the Ehrhart polynomial (see \cite{Ehr62}), i.e. it encodes the number of lattice points in integer dilation of $P$. The connection between $h$- and $h^*$- polynomials was first observed by Betke and McMullen in \cite{BM85} and later generalized by Katz and Stapledon \cite{KS16}. For example, if a lattice triangulation $T$ of the polytope $P$ is unimodular (i.e. the lattice volume of any simplex equals $1$), then $h^* (P;t) = h(T;t)$. For arbitrary triangulation there appears the box-polynomial which is the local $h^*$-polynomial (or the $\ell^*$-polynomial) of a simplex.

The $\ell$- and $\ell^*$-polynomials are the local versions of the $h$- and $h^*$-polynomials. Same as the usual $h$- and $h^*$-polynomials they have non-negative coefficients. These local polynomials are used to decompose the $h$- and $h^*$-polynomials into non-negative contributions when considering polyhedral subdivisions (see \cite{KS16}). 

The local $h$-polynomials naturally appear when applying the decomposition theorem from \cite{BBD81} to the toric morphisms (see \cite{dCMM18}), also see \cite[Notes and References, Part 2, Chapter 11; p. 514]{GKZ94}. The Jordan normal form of the monodromy operator of Newton non-degenerate singularity is given in \cite{Stap17} in terms of the generalizations from \cite{Stap08} of the local $h$- and $h^*$-polynomials.

The local $h$-polynomial for simplicial subdivisions of the simplex was defined in \cite{Stan92}. The relative version of the local $h$-polynomial was introduced independently in \cite[Remark 3.7]{Ath12} for simplicial subdivisions of the simplex and in \cite[Section 3]{N12} for regular triangulations of polytopes. 

\begin{remark}
    We avoid the notion of the relative local $h$-polynomial $\ell(\Sigma,\Delta;t)$ since it is the local $h$-polynomial of the localized triangulation $\ell(\Sigma,\Delta;t) = \ell(\Sigma/\Delta;t)$.
\end{remark}

\subsubsection{Box polynomial or local $h^*$-polynomial $\ell^*(\D;t)$ of a simplex}

Suppose that $\D \subset \R^n$ is a lattice simplex with vertices $e_1,\dots, e_k$.

\begin{definition}(see \cite{BM85})\label{box_def}
    The \emph{box polynomial} (or the \emph{local $h^*$-polynomial}) of a simplex $\D$ is defined as $\ell^*(\D;t) = \sum_i \lambda_i t^i$, where $\lambda_i$ is the number of vectors $(a_1,\dots,a_k): 0 < a_j< 1$ such that $\sum_j a_j e_j \in \Z^n$ and $\sum_j a_j = i$.
\end{definition}

\subsubsection{Local $h$-polynomial $\ell(\Sigma;x)$ of strong formal subdivision}

In \cite{KS16} Katz and Stapledon generalized the notion of the local $h$-polynomial of simplicial subdivisions for the case of so called \emph{strong formal subdivisions}. In this paper we racall only a very particular case of this generalisation, namely, the local $h$-polynomial $\ell(\Sigma;x)$ of strong formal subdivision in case of a simplicial poset subdividing a Boolean poset. For an introduction to the posets (i.e. \textit{partially ordered sets}) see e.g. \cite{Stan}.

\begin{definition}\label{Boolean_simplicial_poset}
    \textit{Boolean} poset is the poset of faces of a simplex. A poset $B$ with minimal element $\hat 0$ is \textit{simplicial} if for every $x$ in $B$, the interval $[\hat 0, x]$ is a Boolean poset. The rank $\rho$ of an element of the poset is the dimension of the corresponding simplex.
\end{definition}

For instance, the poset of the simplices in any simplicial complex is a simplicial poset.

\begin{definition}\cite[Definition 3.16]{KS16}\label{strongly_sur_def}
    Let $\sigma: S \to B$ be an order-preserving, rank-increasing function between simplicial and Boolean posets. Then $\sigma$ is \textit{strongly surjective} if it is surjective and for all $y \in S$ and $x\in B$ such that $\sigma(y) \le x$, there exists $y \le y^\prime \in S$ such that $\rho(y^\prime) = \rho(x)$ and $\sigma(y^\prime) = x$. 
\end{definition}


\begin{definition}\label{strong_formal_def}
    \cite[Definition 3.17]{KS16}
Let $\sigma : S \to B$ be an order-preserving, rank-increasing function between simplicial posets with rank functions $\rho_S$ and $\rho_B$ respectively. Then $\sigma$ is a \textit{strong formal subdivision} if it is strongly surjective and for all $y \in S$ and $x \in B$ such that $\sigma (y) \le x$,
\begin{equation}
    \sum\limits_{\substack{y \le y^\prime \\ \sigma(y^\prime) = x} } (-1)^{\rho_B(x) - \rho_S (y^\prime)} = 1
\end{equation}
\end{definition}

For instance (see \cite[Lemma 3.25]{KS16}) any triangulation of a cone induce a strong formal subdivision.

\begin{definition}\label{local_h_poset_def}
    \cite[Lemma 4.12]{KS16}
Let $\sigma : S \to B$ be a strong formal subdivision between a simplicial poset $S$ and a Boolean poset $B$ with rank functions $\rho_S$ and $\rho_B$ respectively. Then the local $h$-polynomial can be defined as $$\ell_B(S;t) = \sum_{y\in S} (-1)^{\text{rk} (S) - \rho_S(\hat{0}_S,y)} t^{\text{rk} (S) - e(y)} (t - 1)^{e(y)},$$
where $e(y) = \rho_B (\sigma(y)) - \rho_S (y)$ is the \textit{excess} of $y$.
\end{definition} 

Note that the local $h$-polynomial can be negative in this setting, see \cite[Example 5.6]{KS16}.

\begin{claim} \label{linear_coef_claim} \cite[Example 4.10]{KS16}
    The linear coefficient of $\ell_B(S;t)$ for is equal to 
$$\begin{cases}
   \quad \beta - \text{rk} (S) & \text{rk} (B) = 0, \\
   \quad \beta - 1 & \text{rk} (B) = 1, \\
   \quad \beta & \text{rk} (B) > 1.
\end{cases},$$
    where 
    $$\beta = \# \{y\in S| \sigma(y) = \hat 1_B , \rho_S(\hat 0_S, y) = 1\}$$
\end{claim}

\subsubsection{Local $h$-polynomial $\ell(\Sigma;x)$ in $\R^{m}_{\ge 0} \oplus \R^{n-m}$}

We apply the definitions from the previous subsubsection to the case of triangulation $\Sigma$ of a cone $C = \R^{m}_{\ge 0} \oplus \R^{n-m}$ and obtain the local $h$-polynomial $\ell (\Sigma;t)$. Let $\Sigma$ be a triangulation of $C$. For $Q$ a cone in $\Sigma$, let $\sigma_C (Q)$ be the smallest face of the cone $C$ containing $Q$. Set

$$e_C(Q) := \dim \sigma_C(Q) - \dim Q$$

\begin{definition}\label{local_h_def} (cf. Definition \ref{local_h_poset_def})
    Let $\Sigma$ be a triangulation of the cone $C= \R^{m}_{\ge 0} \oplus \R^{n-m}$. Then the local $h$-polynomial $\ell(\Sigma;t)$ can be defined as 
    $$\ell (\Sigma;t) = \sum\limits_{Q \in \Sigma} (-1)^{codim(Q)}t^{n-e_C(Q)}(t-1)^{e_C(Q)}$$
\end{definition}

Note that in the case $n=0$ we have $\ell (\{0\};t) = 1$. This definition works only for triangulations of the cones $\R^m_{\ge 0}\oplus \R^{n-m}$ since their faces form a Boolean poset. In this paper we generalize in Definition \ref{local_hC_def} the local $h$-polynomial to the case of simplicial fans with contractible \emph{sticky} (see Definition \ref{sticky_polytope_def}) support in $\R^m_{\ge 0}\oplus \R^{n-m}$.

\begin{claim}(see e.g. \cite{BKN23})
    The local $h$-polynomial $\ell(\Sigma;t)$ has non-negative coefficients.
\end{claim}


\begin{claim}[cf. Claim \ref{linear_coef_claim}]\label{internal_ray_claim}
    Denote by $\beta$ the number of \textit{interior rays} (i.e. which are not contained in any proper face of $C$) of a triangulation $\Sigma$ of the cone $C = \R_{\ge 0}^m \oplus \R^{n-m}$. Then the coefficient of the local $h$-polynomial $\ell(\Sigma;t)$ of $t$ is equal to

    $$\begin{cases}
    \quad \beta - n & m = 0, \\
   \quad \beta - 1 & m = 1, \\
   \quad \beta & m > 1.
\end{cases}$$

\end{claim}



Note that in \cite[Example 4.10]{KS16} this claim is given in much more generality (for instance, for any polyhedral subdivisions). For regular (coherent) triangulations of $\R_{\ge 0}^n$ it was proved in \cite[Proposition 4.13]{GKZ94}. For arbitrary triangulation of $\R_{\ge 0}^n$ it is proven in \cite[Example 2.3 (f)]{Stan92}. The same argument using \cite[Proposition 1.1]{Mac71} works for arbitrary triangulation of $\R^m_{\ge 0} \oplus \R^{n-m}$, so Claim \ref{internal_ray_claim} can be proved using more primitive technique.

\subsection{Newton numbers and thin simplices}\label{Newton_section}

In this subsection we recall definitions from \cite{GKZ90}; the definitions are also given in \cite[\S 11]{GKZ94}. The Newton numbers arise when considering the rational function called \emph{regular $A$-determinant} (see \cite[Example 2.5]{GKZ94} for the case when the regular $A$-determinant is not polynomial) which is equal to the $A$-discriminant in the case if the projective toric variety corresponding to $A$ is smooth (\cite[Theorem 1.3]{GKZ94}).

\begin{definition}\label{Arbutrary_cone_Newton_number_def}
The \emph{Newton number} $\nu_C$ of a polytope $P$ in a convex lattice cone $C$ is the sum 
$$\nu_C(P) =  \sum_{F\subset C} (-1)^{codim (F)} Vol(F\cap P)$$
over all faces $F \subset C$.
\end{definition}

\begin{definition}\label{V-def}
    Polytope (or cone) in a cone $C$ is called \emph{$V_C$-}polytope (or $V_C$-cone) if its affine span coincides with the affine span of a face of the cone $C$. If it is obvious which cone $C$ is assumed then we write just $V$-polytope.
\end{definition}

In \cite{GKZ94} $V$-simplices are called \emph{massive}.

\begin{definition}\label{Combinatorial_Newton_number_def}
The \emph{combinatorial Newton number} $\ell(\Sigma;1)$ of a triangulation $\Sigma$ in a convex cone $C$ is the sum 
$$CN(\Sigma) =  \sum_{\substack{Q\in \Sigma\\ F \text{ is } V_C\text{- cone}}} (-1)^{codim (Q)}$$
over the cones $Q \subset \Sigma$.
\end{definition}

\begin{remark}\label{difference_remark}
    In this paper we consider only the cases $C = \R^{m}_{\ge 0} \oplus \R^{n-m}$. In these cases the combinatorial Newton number equals the evaluation of the local $h$-polynomial (see Definition \ref{local_h_def}) at $t = 1$. For general cones they are different, for instance, because the local $h$-polynomial has non-negative coefficients, whereas the combinatorial Newton number can be negative (for example, consider a triangulation of the $4$-dimensional cone over the $3$-dimensional cube or the $3$-dimensional triangle prism).
\end{remark}

Denote by $\D^{+1} \subset \R^{n+1}$ the pyramid over $(\D,1)$ with the vertex $O$.

\begin{definition}\label{thin_simplex_def}
    Lattice simplex $\D \subset \R^n$ is called \emph{thin} if the Newton number of the simplex $\D^{+1}$ in its cone equals $0$.
\end{definition}

This definition in case of simplices (but not arbitrary polytopes) coincides with the definition of thin polytopes from \cite{BKN23}. Thin polytopes are completely classified there in dimension up to $3$ in \cite{BKN23}. In our paper we need the following more restrictive definition of $0$-thin simplices instead of the usual thin simplices.

\begin{definition}\label{0_thin_simplex_def}
    Lattice simplex $\D \subset \R^n$ one of whose vertices is the origin is called \emph{0-thin} if the Newton number of the simplex $\D$ in its cone $\Sigma(\D)$ equals $0$.
\end{definition}

\begin{sign}\label{cap_sign}
Consider a lattice simplex $\D$ with vertices $O,w_1,\dots,w_n$ and denote
$$Box_{\D}^\circ = \{w\in \Z^n: w = [ \sum\limits_{i=1}^n \lambda_i w_i], 0<\lambda_i<1\}, \quad Cap (\D) = |Box_{\D}^\circ|$$
We also set $Cap (O) = 1$.
\end{sign}

Denote by $\D^{-O}$ the facet of $\D$ which does not contain the origin.

\begin{remark}
    \begin{enumerate}
        \item Simplex $\D$ is thin is equivalent to $Cap(\D^{+1}) = 0$ and to $\ell^*(\D;t) = 0$.
        \item Simplex $\D$ is 0-thin is equivalent to $Cap(\D) = 0$ and to $\ell^*(\D;t) = \ell^*(\D^{-O};t) = 0$.
    \end{enumerate}
\end{remark}

\subsection{Thin triangulations or vanishing local $h$-polynomial} \label{thin_sec}

\begin{definition}\label{thin_triangulation_def}
    Triangulation $\Sigma$ of a cone $C = \R^m_{\ge 0} \oplus \R^{n-m}$ is called \emph{thin} if $\ell(\Sigma;1) = 0$.
\end{definition}

From Claim \ref{internal_ray_claim} we obtain the following corollary.

\begin{cor}\label{internal_ray_cor}
    Thin triangulations of the cones $\R^m_{\ge 0} \oplus \R^{n-m}$ for $m > 1$ contain no interior rays. For $m=0$ and $n>0$ there are no thin triangulations. If $m = 1$ then triangulation is thin if and only if it contains exactly one interior ray.  
\end{cor}

Thin triangulations are also called \emph{vanishing local $h$-polynomials} in \cite{dMGP+20}. In this subsection we recall the classification of thin triangulations in small dimensional.

\subsubsection{Conical facet refinements}

There is a class of refinements of triangulations called \emph{conical facet refinement} that preserves the local $h$-polynomial. So the classification of simplicial fans with vanishing combinatorial Newton numbers is always considered up to a sequence of conical facet refinements. In this subsection we briefly recall the required definitions and claims from \cite{dMGP+20} and \cite{GKZ94}.

Consider a simplicial fan $\Sigma$ in a cone $C \supset Supp(\Sigma)$. Recall that fan $\Sigma^\prime$ \emph{refines} fan $\Sigma$ if every cone of $\Sigma$ is the union of some cones of $\Sigma^\prime$.  \emph{Facet} of a fan is its top-dimensional cone.

\begin{definition}\label{facet_ref_def}
   \emph{Facet refinement} $\Sigma^\prime$ of simplicial fan $\Sigma$ along facet $F$ is a simplicial refinement such that for any other facet $K$ of $\Sigma$ the restriction $\Sigma^\prime_K$ of $\Sigma^\prime$ on $K$ is trivial. If $Supp (\Sigma) \subsetneq C$ then we also assume that the restriction of $\Sigma^\prime$ on the closure of $C \setminus Supp(\Sigma)$ is trivial.
\end{definition}

\begin{claim}[\cite{dMGP+20}]
    Let $\Sigma^\prime$ be a facet refinement of $\Sigma$ along $F$. Then we have
    $$\ell(\Sigma^\prime;1) = \ell(\Sigma;1) + \ell(\Sigma^\prime_F;1)$$
\end{claim}

\begin{definition}\label{conical_def}
    Facet refinement $\Sigma^\prime$ of $\Sigma$ along $F$ is called \emph{conical facet refinement} if $\Sigma^\prime_F$ is the cone over $\Sigma^\prime_H$ for some codimension $1$ face $H\subset F$. The edge of $F$ outside the face $H$ is denoted by $W(F)$.
\end{definition}

\begin{figure}[H]
		\begin{center}
		\includegraphics[scale = 5]{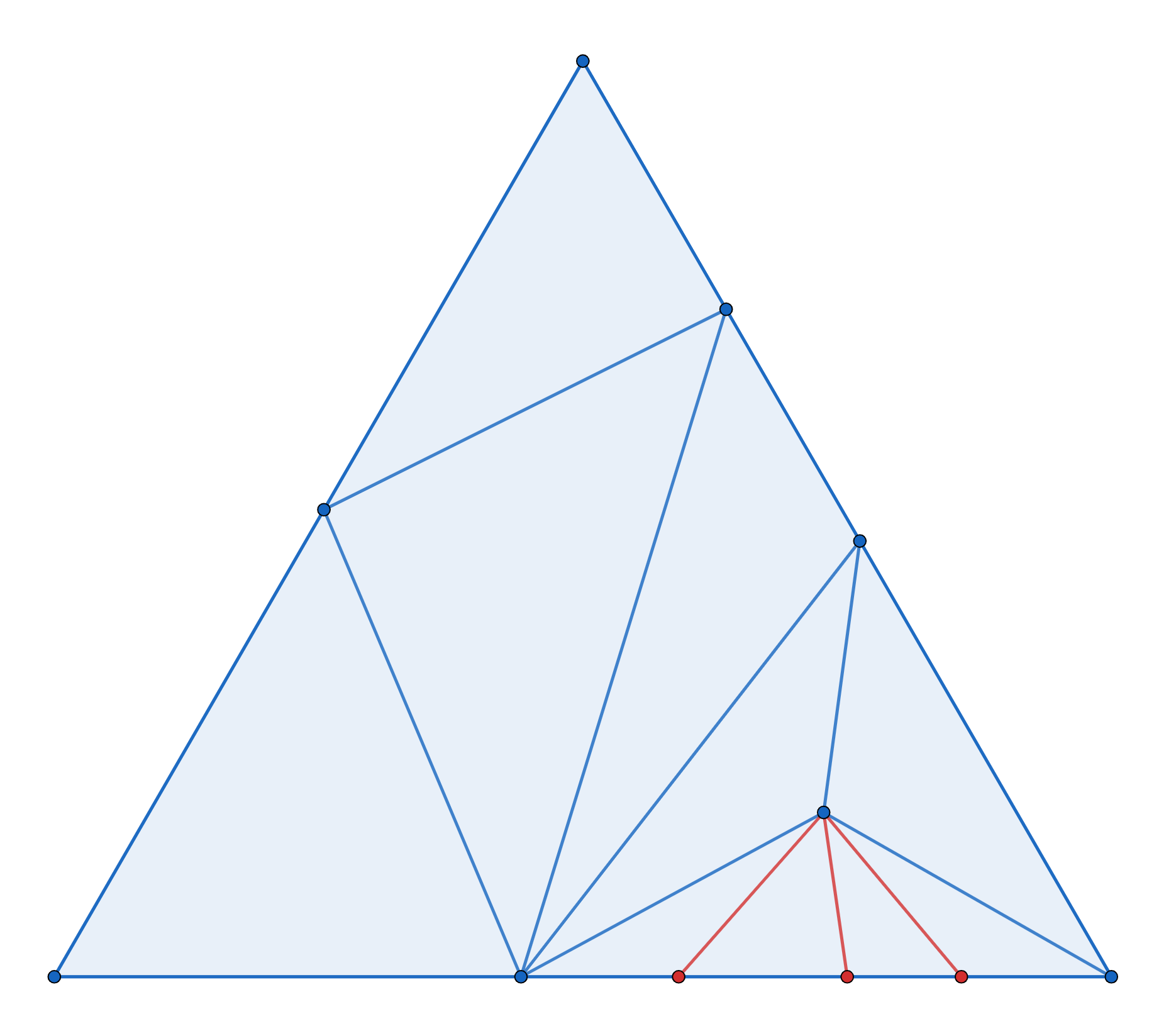}
		\end{center}
    \caption{\label{Conical_facet_refinement_picure} Projectivisation of a \textbf{conical facet refinement}}
\end{figure}

\begin{cor}[\cite{dMGP+20}]\label{conical_corollary}
    If $\Sigma^\prime$ is a conical facet refinement of $\Sigma$, then $\ell(\Sigma^\prime;1) = \ell(\Sigma;1)$. In particular, if $\Sigma^\prime$ is obtained from a fan with vanishing combinatorial Newton number by a sequence of conical facet refinements then $\ell(\Sigma;1) = 0$.
\end{cor}

\subsubsection{Classification of thin triangulations in small dimension} \label{Classification_of_vanish_subsec}

The problem of classifying thin triangulations was posed in \cite{GKZ90} for regular (coherent) triangulations and in \cite{Stan92} for triangulations of the cones $\R^n_{\ge 0}$. This classification for the cones $C = \R^{m}_{\ge 0} \oplus \R^{n-m}$ is important for the local monodromy conjecture (see \cite{LPS22a}) and for the Arnold's monotonicity problem. The classification is still not obtained. In Definition \ref{local_hC_def} we generalize the local $h$-polynomial to the case of so called \emph{strongly contractible} simplicial fans and (conjecturing its non-negativity) ask the same question concerning them. 

This problem is solved for regular triangulations for $m$ equals $0$ or $1$ and $(m,n)$ equals $(2,0), (3,0)$ and $(2,1)$ in \cite[Proposition 4.14]{GKZ94} (we cite the paper \cite{GKZ94} instead of \cite{GKZ90} because there are some inaccuracies in \cite{GKZ90}). But their arguments work for arbitrary triangulation since their classification is implied by Corollary \ref{internal_ray_cor}. In \cite{dMGP+20} the classification is obtained in the case $(4,0)$ and renewed in the case $(3,0)$ using the conical facet refinement. There is no complete classification for any other pair $(m,n)$. For other results concerning thin triangulations see \cite{LPS22a} and \cite{LPS23}.

Let us recall the low-dimensional classifications of thin triangulations and forecast the classification in the cases $(2,2)$ and $(3,1)$.

\begin{enumerate}
    \item[$\mathbf{0 , n:}$]

    No thin triangulations.
    
    \item[$\mathbf{1 , n-1:}$]

    Triangulation $\Sigma$ is thin if and only if there is unique ray $\tau \in \Sigma$ which is not contained in $\R^{n-1}$ (see \cite{GKZ94}).

    \item[$\mathbf{2 ,0:}$]

    Triangulation $\Sigma$ is thin if and only if it is the trivial triangulation of the positive quadrant (see \cite{GKZ94}).
    
    \item[$\mathbf{3,0:}$] Triangulation $\Sigma$ is thin if and only if it is obtained from the trivial or the triforce triangulation by a sequence of conical facet refinements (see \cite{dMGP+20} and \cite{GKZ94}).
\begin{figure}[H]
		\begin{center}
		\includegraphics[scale = 3]{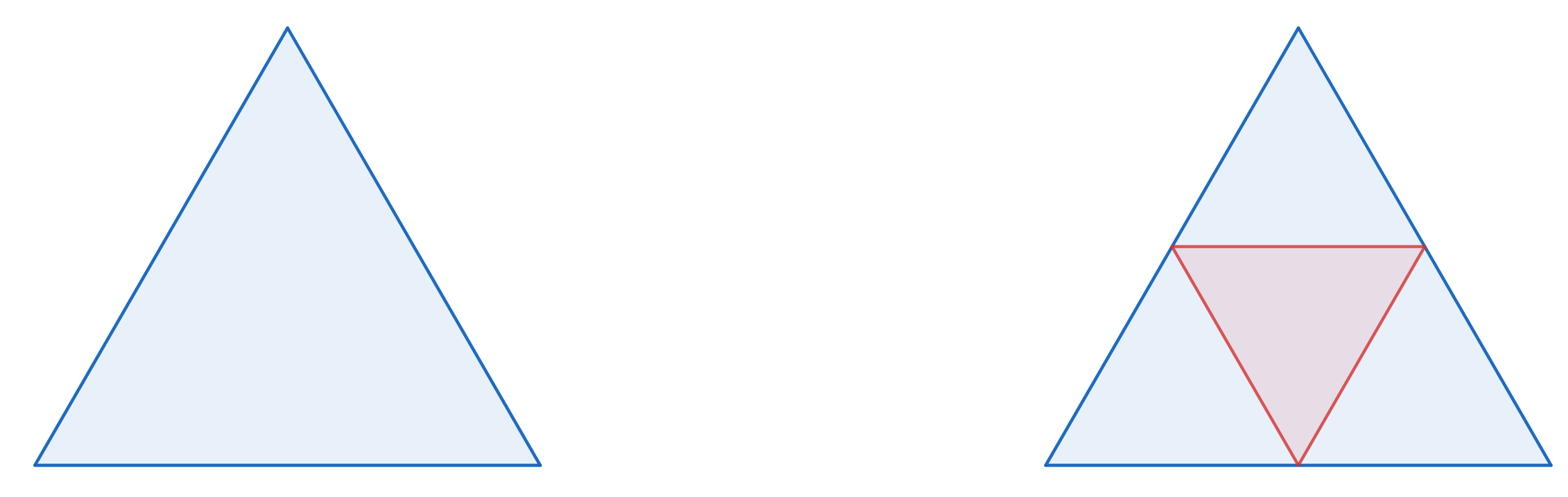}
		\end{center}
    \caption{\label{Trivias_and_Trifold} Projectivisation of the trivial and the triforce triangulations}
\end{figure}

    \item[$\mathbf{2,1:}$]
    Triangulation $\Sigma$ is thin if and only if it is obtained from a \emph{trivial triangulation} by a sequence of conical facet refinements, where a \emph{trivial triangulation} is a triangulation of the cone $\R^2_{\ge 0} \oplus \R^1$ containing only two facets (see \cite{GKZ94}).

\begin{figure}[H]
		\begin{center}
		\includegraphics[scale = 2]{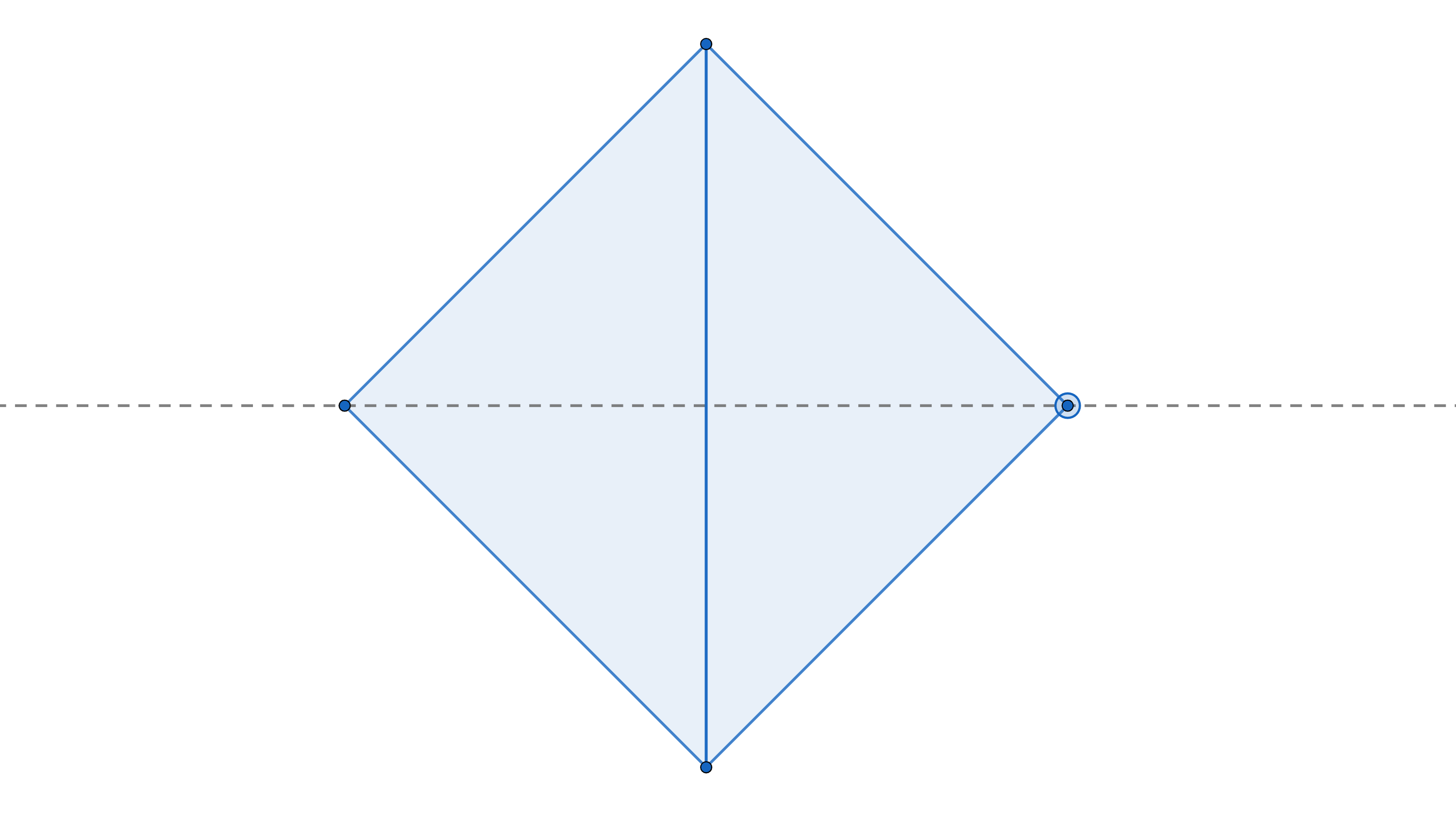}
		\end{center}
    \caption{\label{1+2-1} Projectivisation of a trivial triangulation}
\end{figure}

    \item[$\mathbf{4,0:}$] Triangulation $\Sigma$ is thin if and only if it is obtained from the trivial one by a sequence of conical facet refinements (see \cite{dMGP+20}).

    \item[$\mathbf{2,2 \text{ \& } 3,1}$]
    We believe that the same technique as in \cite{dMGP+20} can confirm the following conjecture. If the conjecture holds then we can extend our solution to the \textbf{complete solution of the Arnold's monotonicity problem in dimension $\mathbf{5}$} (see Theorem \ref{4_5_Arnold_theorem}).
    \begin{con}[Thin triangulations conjecture] \label{thin_triang_con}
        All thin triangulations of the cones $\R^2_{\ge 0} \oplus \R^2$ and $\R^3_{\ge 0} \oplus \R^1$ can be obtained from a trivial triangulation (see Definition \ref{Trivial_fan_def}) by a sequence of conical facet refinements.
    \end{con}

\end{enumerate}

\section{Main results} \label{main_results_section}

In this section we formulate the main definitions and results of the paper without proof.

\begin{definition}\label{star-polytope_def}
    A \emph{star-polytope} in a cone $C$ is a polytope admitting a star-triangulation, i.e. a triangulation all whose maximal simplices contain the origin. It is called \emph{complete} if the cone over it is the whole cone $C$.
\end{definition}

\begin{definition}\label{star-negligible_def}
    Polytope $P$ is called \emph{negligible} in the cone $C = \R^m_{\ge 0} \oplus \R^{k-m}$ if $\nu_C(P) = 0$.
\end{definition}

The following criterion reduces the Arnold's monotonicity problem in dimension $n$ to classifying the Negligible $(n-1)$-dimensional star-polytopes in the cones $\R^m_{\ge 0 } \oplus \R^{n-m-1}$ and provides a tool for classifying them.

\begin{theorem}[Arnold's monotonicity criterion]\label{Arnold_monotonicity_criterion_theorem} \ 

\begin{enumerate}
    \item (Negligible criterion) 
    A pair $(\Gamma_+,G) \subset \R^n_{\ge 0}$ satisfies the Arnold's monotonicity problem if and only if the complete star-polytope $\Pi_G(\Gamma) = \pi_{\mathsmaller{OG}}(\Gamma_{+G} \setminus \Gamma_+)$ is negligible in the cone $C_G = \pi_{\mathsmaller{OG}}(\R^n_{\ge 0})$.
    
    \item (Thin criterion) 
     Complete star-polytope $P\subset C$ containing a neighbourhood of the cone $C = \R^m_{\ge 0} \oplus \R^{k-m}$ with star-triangulation $T$ is negligible if and only if for each simplex $O\le \D\in T$ either $\D$ is $0$-thin or the localized triangulation $T/\D$ is thin.
\end{enumerate}

\end{theorem}

\begin{remark}
    The cone $C_G = \pi_{\mathsmaller{OG}}(\R^n_{\ge 0})$ equals the cone $\R^m_{\ge 0} \oplus \R^{n-m-1}$, where $n-m$ is the dimension of the minimal face of $\R^n_{\ge 0}$ containing the point $G$.
\end{remark}

The following Negligible example lemma \ref{Negligible_examples_lemma} provides a wide class of examples of the Negligible complete star-polytopes, but not all of them (see Example \ref{Unpredicted_example} in dimension $5$).

\begin{definition}\label{Trivial_fan_def}
    A triangulation $\Sigma$ of the cone $C = \R^m_{\ge 0}\oplus \R^{n-m}$ is called \emph{a trivial triangulation} if there is a cone $C_m \cong \R^m_{\ge 0} \subset C$ and a fan $\Sigma_{n-m}$ with support $\R^{n-m} \subset C$ such that $\Sigma = \Sigma_{C_m} + \Sigma_{n-m}$, where $\Sigma_{C_m}$ is the fan of faces of the cone $C_m$.
\end{definition}

Consider a complete star-polytope $P$ with star-triangulation $T$. For a ray $W\in \Sigma(T/O)$ denote by $W^p\ne O$ the vertex of $W \cap P$.

\begin{lemma}[Negligible example lemma]\label{Negligible_examples_lemma}
    Suppose that the triangulation $T/O$ is obtained from a trivial triangulation of the cone $C=\R^m_{\ge 0}\oplus \R^{n-m}$ by a sequence of conical facet refinements along some facets $F_i, i = 1,\dots ,k$ with the corresponding codimension $1$ faces $H_i \subset F_i$ and rays $W(F_i)\subset F_i$ outside $H_i$. Denote $W_i= (W(F_i))^p$. Suppose that each $W_i$ has height $1$ with respect to the corresponding face $H_i$. Then the complete star-polytope $P$ is negligible in the cone $C$.
\end{lemma}

\begin{remark}
    If the triangulation $T/O$ is already a trivial triangulation then the star-polytope $P$ is negligible if and only if it contains a vertex of height $1$ with respect to some facet of the cone $C$.
\end{remark}

\begin{figure}[H]
		\begin{center}
		\includegraphics[scale=5]{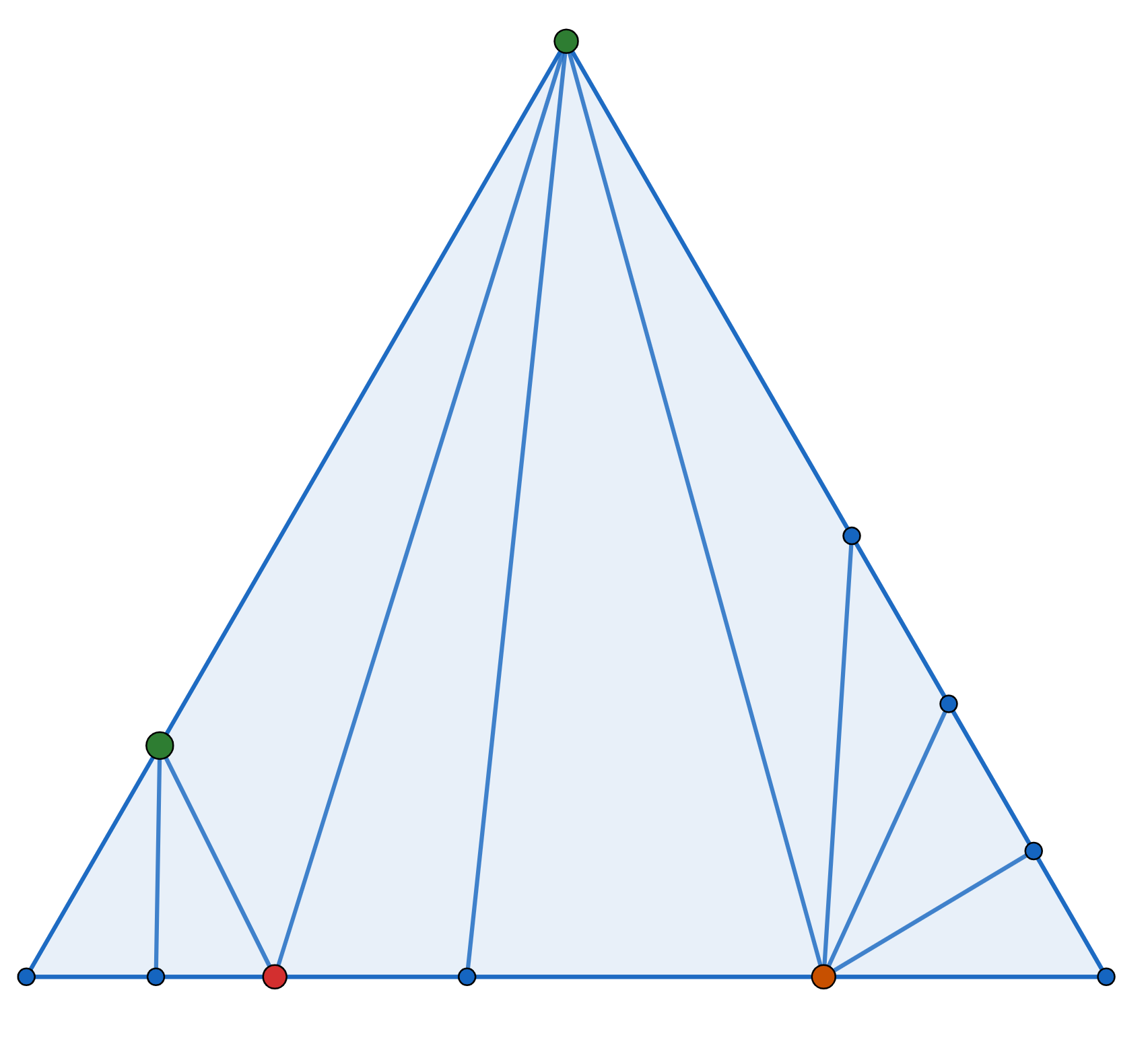}
		\end{center}
    \caption{\label{Examples_lemma_picture} Projectivisation of a negligible star-polytope. Colored vertices have height $1$}
\end{figure}

Together with the Negligible criterion from Theorem \ref{Arnold_monotonicity_criterion_theorem} the following theorem provides complete solution to the Arnold's monotonicity problem in dimension $4$ and partial solution in dimension $5$ (which is complete up to Thin triangulations conjecture \ref{thin_triang_con}).

\begin{theorem}[Arnold's problem in dimension $4$ and $5$]\label{4_5_Arnold_theorem} \
    \begin{enumerate}
        \item Negligible example lemma \ref{Negligible_examples_lemma} provides all  examples of Negligible complete star-polytopes in the cones $\R^m_{\ge 0} \oplus \R^{3-m}$.
        \item Negligible example lemma \ref{Negligible_examples_lemma} provides all  examples of Negligible complete star-polytopes in the cones $\R^m_{\ge 0} \oplus \R^{4-m}$ for $m \neq 2,3$.
    \end{enumerate}
\end{theorem}

\begin{remark}
The cases $m = 0,1$ are trivial for any $n$. Confirming the Thin triangulations conjecture \ref{thin_triang_con} proves that the Negligible example lemma \ref{Negligible_examples_lemma} provides complete solution to the Arnold's monotonicity problem in dimension $5$. The Arnold's monotonicity problem in dimension $6$ has solutions which can not be obtained by the Negligible example lemma (see Example \ref{Unpredicted_example}).
\end{remark}

The main tool in the proof of the Arnold's monotonicity criterion \ref{Arnold_monotonicity_criterion_theorem} is the non-negative analogue of the Kouchnirenko formula Theorem \ref{Non-negative_analogue_theorem}. It consists of three parts: the \textbf{Star-formula}, the \textbf{Links-formula} and the \textbf{Merged formula}. The non-negative analogue is also one of the main results and is valuable in itself. It has very close prototypes in \cite{Stap17} and \cite{LPS22a} for the Star-formula and in \cite{Fur04} for the Links-formula, we discuss the relations in more details in \S \ref{Non_negative_sec}.

The main definition of the Links-formula is \emph{sticky polytope} (see Definition \ref{sticky_polytope_def}). Sticky polytope is a generalization of the polytopes obtained as the difference of two convenient Newton polytopes. The formula is given for triangulations of sticky polytopes. Another key players is $V_{\min}$-simplex, i.e. $V$-simplex with no proper $V$-faces. The Newton number of the triangulated sticky polytope is decomposed into the sum of the Newton numbers of the links of the $V_{\min}$-simplices, which are also sticky polytopes. The Star-formula relies on the $\ell_C$-polynomial (see Definition \ref{local_hC_def}) which is a generalization of the local $h$-polynomials of triangulations to so called \emph{strongly contractible} fans. 

The proof of the Negligible criterion from Theorem \ref{Arnold_monotonicity_criterion_theorem} is tricky. It uses the Links-formula and the \emph{top and bottom} triangulations from Definitions \ref{top_triang_def} and \ref{bottom_triang_def}. The Thin criterion from Theorem \ref{Arnold_monotonicity_criterion_theorem} immediatly follows from the Star-formula.

The proof of the Theorem \ref{4_5_Arnold_theorem} is based on the $B_1$-simplex Lemma \ref{B_1-simplex_lem} which follows from the Sprig lemma \ref{sprig_lemma}. The Sprig lemma explains why the negligible polytopes for the cones $\R^{m^\prime}_{\ge 0} \oplus \R^{n-m^\prime}$ appear as parts of the negligible polytopes in the cone $\R^{m}_{\ge 0} \oplus \R^{n-m}$ along $(m-m^\prime)$-dimensional faces (as it happened in \S \ref{4-Arnold_intro_sec} when the crook connected two greenhouses). The Sprig lemma is the place where we essentially need the generalized version called $\ell_C$-polynomial of the local $h$-polynomial instead of the usual one. Confirming the $\ell_C$-positivity conjecture \ref{non_neg_local_sticky_con} would be useful for verifying the conditions of the Sprig lemma \ref{sprig_lemma}.

\newpage

\section{Further structure of the paper}

\begin{figure}[H]
		\includegraphics[scale=.36]{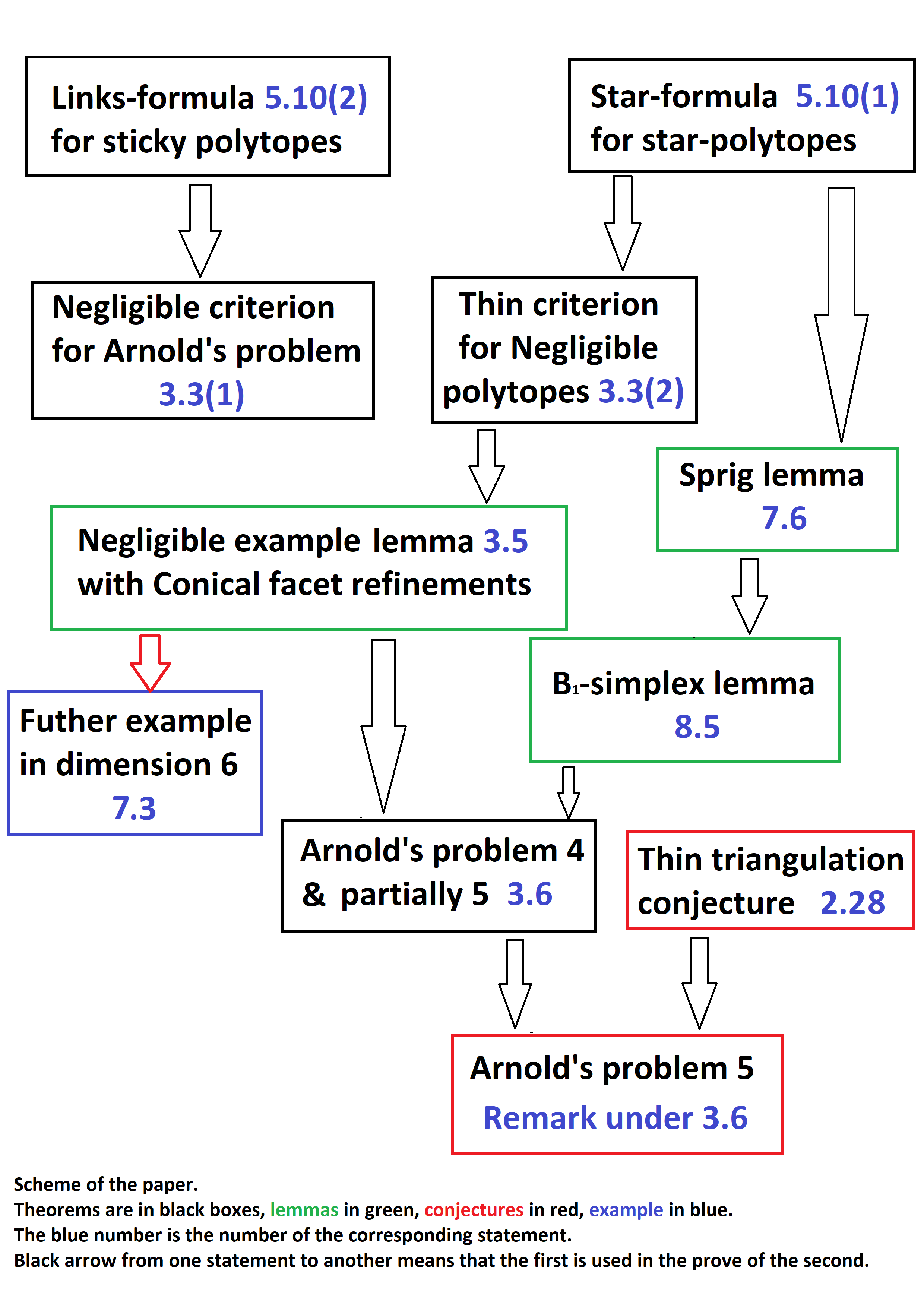}
\end{figure}

\newpage

\section{Non-negative analogue of the Kouchnirenko formula}\label{Non_negative_sec}

In this section formulate we give the non-negative analogue of the Kouchnirenko formula Theorem \ref{Non-negative_analogue_theorem} without proof. It consists of three parts: two disjoint parts, called the \textbf{Star-formula} for the Newton number of sticky star-polytopes and the \textbf{Links-formula} for the Newton number of sticky polytopes. The last part called the \textbf{Merged formula} is obtained by substituting the Star-formula into the Links formula. We prove it all in the next section \ref{proof_of_non_neg_section}. Definitions required for the Links-formula are given in \S \ref{sticky_subsec}, for the Star-formula in \S \ref{local_hC_subsec}. The non-negative analogue is formulated in \S \ref{Non_negative_analogue_subsection}

In the case of coconvex star-polytopes in $\R^n_{\ge 0}$ the Star-formula is the calculation of the degree of the non-negative formula for the monodromy operator (see \cite{Stap17}; see also \cite{LPS22a} for a simplified version). The original (with negative degrees) Varchenko formula for the monodromy operator in given in \cite{Var76}. For a basic introduction to the Picard–Lefschetz theory including the monodromy operator see \cite{AGLV}. For more detailed reading concerning the theory of geometric monodromies see \cite{T}. We give a short combinatorial proof of the star-formula in \S \ref{star_formula_subsec}.

The Links-formula for the Newton number has a close prototype in \cite{Fur04}. But for our applications we need to define and prove it for a different class of polytopes, called \textit{sticky polytopes} (Definition \ref{sticky_polytope_def}), although the proof doesn't change much. Sticky polytopes form a more natural class of polytopes for the input in the Links formula, for instance, in contrast to \cite{Fur04}, the difference of a sticky polytope and a sticky subpolytope is also a sticky polytope. We also note that in \cite{Fur04} in the proof of Theorem 2.3 a stronger version of the following statement (without any comment or proof) is used.
 
\begin{con}[Separation conjecture]\label{Separation_conj}
    Suppose that $P\subset \R^n_{\ge 0}$ is a (non necessarily convex) polytope containing the simplex $\R^n_{\ge 0} \cap \{\sum_i x_i \le 1\}$. Then the set $P \cap \{\sum_i x_i \ge 1\}$ is also a polytope (i.e. admits a lattice triangulation).
\end{con}

We believe that this conjecture is not true but do not have counterexamples. In the paper \cite{Fur04} the non-negativity of the Links decomposition formula relies on this statement. In our paper the non-negativity conditions are implied by the non-negativity of the Star-formula \ref{Non-negative_analogue_theorem}(1) for the Newton number of star-polytopes, so we do not need this statement. 

However the inaccuracy from \cite{Fur04} can also be fixed in the following way. The notion of regular triangulations can be extended to some non-convex polytopes, and this class of polytopes is preserved by the operations arising in the Links-formula. Namely, consider a regular triangulation induced by some height function of some convex polytope of dimension $n$. Define the height of $n$-simplex as the height of the intersection of the lifted simplex with the vertical line passing through the origin of $\R^n$. Consider the geometric simplicial complex spanned by all $n$-simplices whose height is contained in some segment. We say that these simplices form a regular triangulation of their union (their union is a polytope but not necessarily convex). In the case of these extended regular triangulations the Separation conjecture \ref{Separation_conj} can be easily confirmed by choosing an appropriate height function, so the Links decomposition formula can provide a non-negative analogue of the Kouchnirenko formula on its own.

\subsection{Sticky polytopes and sticky links decomposition} \label{sticky_subsec}

Consider cone $C = \R^m_{\ge 0} \oplus \R^{n-m}$ and a polytope $P \subset C$. Denote by $\partial_C (P)$ the boundary of $P$ with respect to the topology of $C$. For example, if $C = \R^n_{\ge 0}$ and $P = \{\sum_i x_i \le 1\} \cap \R^n_{\ge 0}$ is the unit simplex then $\partial_{\R^n_{\ge 0}}(P) = \{\sum_i x_i = 1\} \cap \R^n_{\ge 0}$.

\begin{definition}\label{sticky_polytope_def}
Polytope $P\subset C$ is called \emph {sticky} if its boundary $\partial_C(P)$ contains no $V_C$-faces. The same definition works for cones contained in $C$.
\end{definition}

For example, the complement to any convenient Newton polyhedron in any cone and the difference between any pair of embedded convenient Newton polyhedra are sticky polytopes. Here is an example of not sticky simplex in $\R^3_{\ge 0}$ with negative Newton number.

\begin{figure}[H]
		\begin{center}
		\includegraphics[scale=0.5]{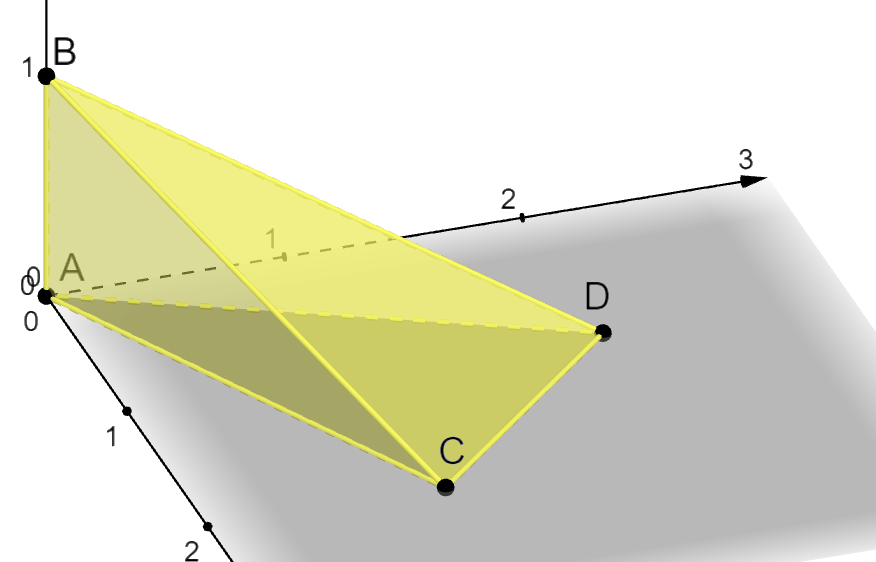}
		\end{center}
		\caption{\label{Not_sticky} Not \textbf{sticky} simplex in $\R^3_{\ge 0}$ with \textbf{negative Newton number} ($\nu = -1$).}
\end{figure}

\begin{definition}\label{V-min_def}
    Polytope is called $V_{\min}$ if it is a $V$-polytope and contains no proper $V$-faces.
\end{definition}

\begin{remark}\label{unique_min_lemma}
    Every $V$-simplex contains unique $V_{\min}$-face.
\end{remark}

\begin{sign}\label{P_T_sign}
    Consider a triangulation $T$ of a polytope $P$ and a simplex $\D \in T$. Denote by $P_T(\D)$ the union of the simplices in the link $\text{lk}_T(\D)$.
\end{sign}

\begin{lemma}[Links decomposition lemma]\label{Decomposition_lemma}
    Consider a sticky polytope $P$ in a cone $C$ with triangulation $T$. Then the triangulation can be coarsened to sticky links of the $V_{\min}$-simplices
    $$P = P_T(\D) := \bigsqcup \limits_{\D \in T \text{ is } V_{\min}} P_T(\D)$$
\end{lemma} 

\begin{figure}[H]
		\begin{center}
		\includegraphics[scale=1]{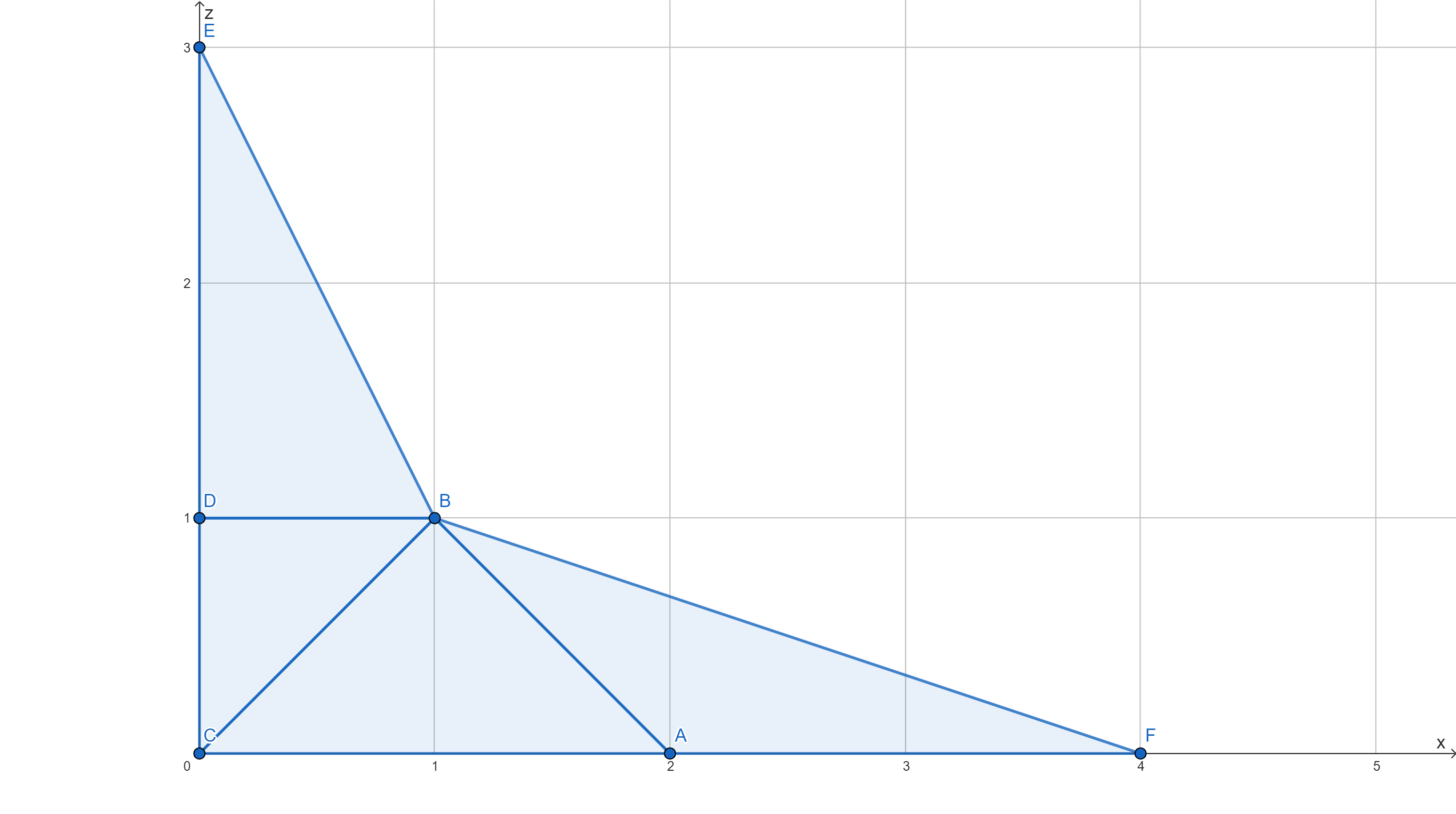} \quad
  \includegraphics[scale=1]{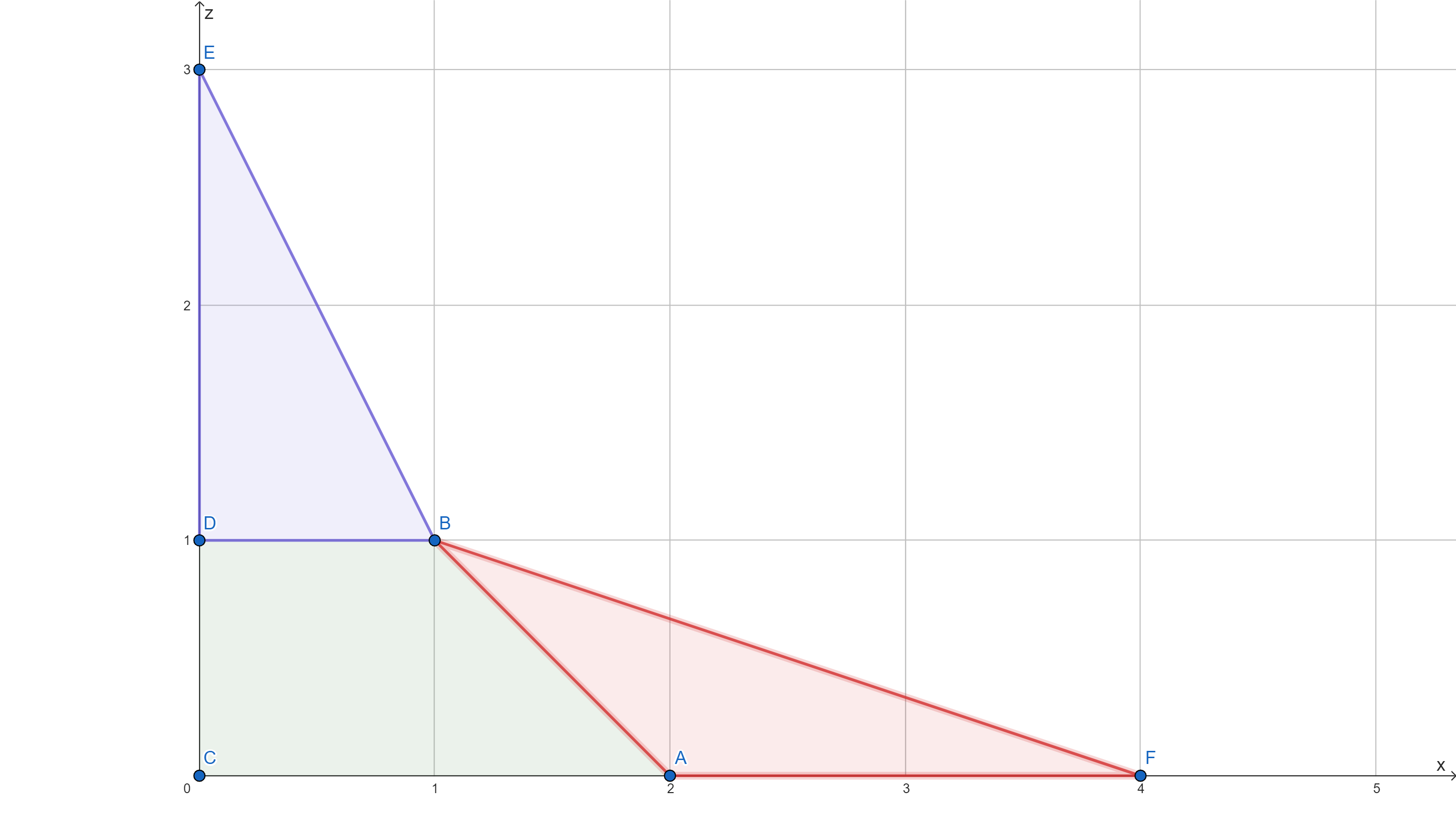}
		\end{center}
    \caption{\label{Decomposition_2} \textbf{Links decomposition} in dimension $2$}
\end{figure}

\begin{figure}[H]
		\begin{center}
\includegraphics[scale=.5]{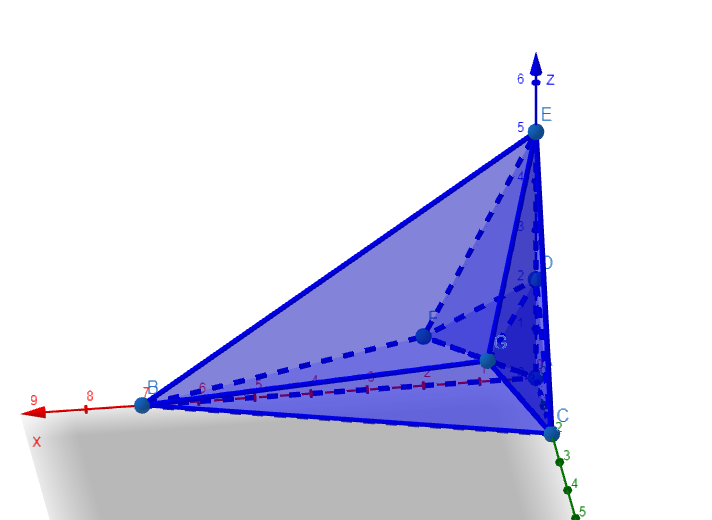} 

  \includegraphics[scale=.3]{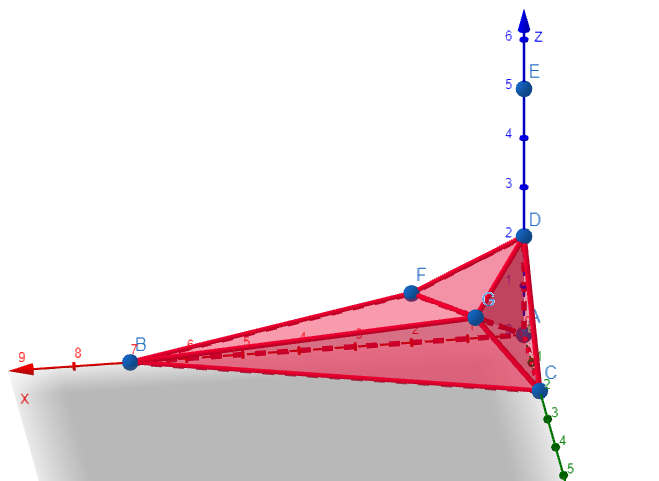}
  \includegraphics[scale=.3]{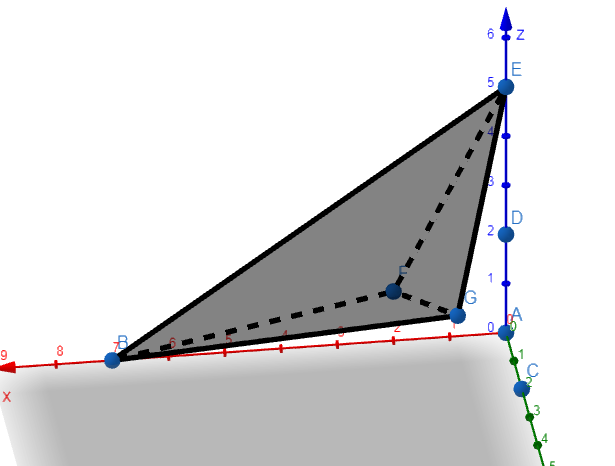}
  \includegraphics[scale=.3]{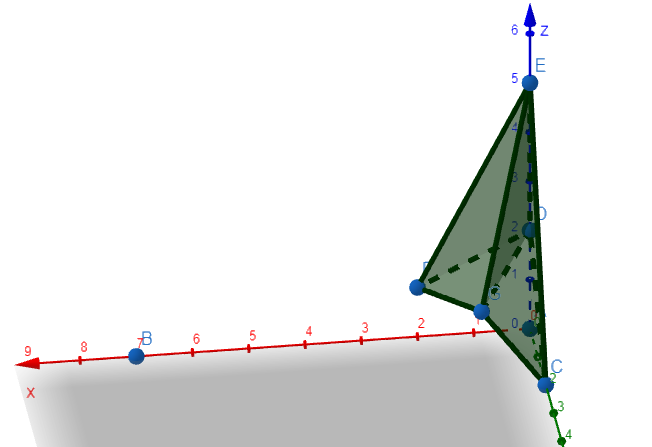}
		\end{center}
    \caption{\label{Decomposition_3} \textbf{Links decomposition} in dimension $3$}
\end{figure}

\begin{proof}[Proof of Lemma \ref{Decomposition_lemma}]
    Polytopes $P_T(\D)$ have non-intersecting interiors because each $V$-simplex (and hence each $n$-simplex) contains unique $V_{\min}$-simplex. The union of the polytopes $P_T(\D)$ is the whole polytope $P$ because each $n$-simplex is a $V$-simplex and contains a $V_{\min}$-simplex. So the only thing left to prove is that each polytope $P_T(\D)$ is sticky.
    
    If the polytope $P_T(\D)$ is not sticky then there is a $V$-simplex in the boundary $\partial_C (P_T(\D))$ which does not contain the simplex $\D$. Then there is a $V_{\min}$-simplex $\D^\prime$ in $\partial_C (P_T(\D))$ which does not contain the simplex $\D$. Note that any $n$-simplex in lk$_T(\D)$ containing $\D^\prime$ also contains $\D$. But a simplex can contain only one $V_{\min}$-simplex. So the polytope $P_T(\D)$ is sticky.
\end{proof}

\subsection{The $\ell_C$-polynomial for the cones $C = \R^m_{\ge 0} \oplus \R^{n-m}$}\label{local_hC_subsec}

In this subsection we extend the definition of the local $h$-polynomial for so called \emph{strongly contractible} sticky fans. In this paper we can use the Remark under Definition \ref{local_hC_def} as the definition of the $\ell_C$-polynomial (we do not use any of its properties in this paper). The following definitions are needed to explain the notion of the $\ell_C$-polynomial and its connection with the $\ell$-polynomial.

For a fan $\Sigma$ with sticky support in $C$ denote by $\Sigma_{\text{int}} \subset \Sigma$ the set of all cones which are not contained in the boundary $\partial_C (Supp (\Sigma))$. 

\begin{definition}\label{strongly_contractible_def}
 Fan $\Sigma$ with sticky support in a cone $C$ is called \emph{strongly contractible} if for any $Q\in \Sigma_{\text{int}}$ and any face $F$ of the cone $C$ containing $Q$ the set $\bigcup\limits_{\substack{Q^\prime \ge Q\\ \sigma(Q^\prime) = F}} \text{int}(Q^\prime)$
    is non-empty, open (in the topology of $F$) and contractible.
\end{definition}

For example, any fan with convex sticky support is strongly contractible, but not every fan with sticky support is strongly contractible.

\begin{claim}
        Suppose that $\Sigma$ is a strongly contractible sticky fan in the cone $C = \R^m_{\ge 0} \oplus \R^{n-m}$. Denote by $B$ the Boolean poset of the faces of the cone $C$. Denote by $\sigma_{\text{int}}: \Sigma_{\text{int}} \to B$ the function mapping simplex to the minimal face of $C$ containing it. Then $\sigma_{\text{int}}$ is a strong formal subdivision.
\end{claim}

The following definition is obtained if we input the strong formal subdivision corresponding to the set $\Sigma_{\text{int}}$ in the Definition \ref{local_h_poset_def}.

\begin{definition}\label{local_hC_def} (cf. Definition \ref{local_h_def})
    Let $\Sigma$ be a strongly contractible simplicial fan in the cone $C= \R^{m}_{\ge 0} \oplus \R^{n-m}$. Then the local $h$-polynomial $\ell(\Sigma;t)$ can be defined as 
    $$\ell_C (\Sigma;t) = \sum\limits_{Q \in \Sigma_{\text{int}} } (-1)^{codim(Q)}t^{n-e_C(Q)}(t-1)^{e_C(Q)}$$
    If $T$ is a star-triangulation, then by $\ell_C(T;t)$ we denote the polynomial $\ell_C(T/O;t)$.
\end{definition}

\begin{remark}
    We will only use the evaluation of $\ell_C$ at $t = 1$ which can be written in the following way:
    $$\ell_C(\Sigma;1) =  \sum_{\substack{Q\in \Sigma\\ F \text{ is } V_C\text{- cone}}} (-1)^{codim (Q)}$$
\end{remark}

Note that the only sticky cone in the cone $\R^n_{\ge 0}$ is the whole $\R^n_{\ge 0}$, but it is not true for not sharp cones. From Claim \ref{linear_coef_claim} it follows that the linear coefficient of $\ell_C(\Sigma;t)$ is the number of the interior rays of $\Sigma_{\text{int}}$. We believe that most of the properties of the local $h$-polynomials (see e.g. \cite{BKN23}) can be extended to the $\ell_C$-polynomials (even in case of polyhedral fans), in particular, the following one.

\begin{con}[$\ell_C$-positivity conjecture] \label{non_neg_local_sticky_con}
    The $\ell_C$-polynomial has non-negative coefficients for any strongly contractible simplicial fan $\Sigma$ in the cone $C$.
\end{con}

\begin{remark}
The Negligible criterion (Theorem \ref{Arnold_monotonicity_criterion_theorem} (1)) also works if instead of the projection along $OG$ we consider the projection along any line passing through $G$ and containing a ray inside the boundary of the cone with the apex $G$ and the base $\Gamma \setminus \Gamma_G$. If the ray is on the boundary of the cone then the assumption of the Thin criterion (Theorem \ref{Arnold_monotonicity_criterion_theorem} (2)) that $\Pi_G(\Gamma)$ is complete is not applied. We believe that the Thin criterion works for arbitrary star-polytope but in order to prove it we need to confirm the $\ell_C$-positivity Conjecture.
\end{remark}

The $\ell_C$-positivity Conjecture is also needed for better understanding of the connection between Negligible polytopes in different cones (see \S \ref{Sprigs_subsecion} and the Sprig lemma \ref{sprig_lemma}).

\subsection{Formulation of the non-negative analogue of the Kouchnirenko formula}\label{Non_negative_analogue_subsection}

Consider a $V$-simplex $\D$. Denote by $m(\D)$ that face of $\D$ which is a $V_{\min}$-simplex.

\begin{theorem}[Non-negative analogue of the Kouchnirenko formula]\label{Non-negative_analogue_theorem}

    Consider a triangulation $T$ of a sticky polytope $P$ in a cone $C = \R^m_{\ge 0} \oplus \R^{n-m}$. Then the following non-negative formulas for the Newton number of the polytope $P$ are held.
    
    \begin{enumerate}
        \item (Star-formula) If $T$ is a star-triangulation and $P$ is a complete star-polytope, then we have
        $$\nu_C(P) = \sum\limits_{0\le \D \in T} Cap(\D) \cdot \ell_{C/\D}(T/\D;1)$$
        \item (Links-formula) For any triangulation we have
        $$\nu_C(P) = \sum\limits_{\D \in T \text{ is $V_{\min}$}} Vol(\D) \cdot \nu_{C/\D} (P_T(\D)/\D)$$
        \item (Merged formula) For any triangulation we have
        $$\nu_C(P) = \sum\limits_{\D \in T \text{ is $V$}} Vol(m(\D))\cdot Cap(\D/m(\D)) \cdot \ell(T/\D;1)$$
    \end{enumerate}
\end{theorem}

\begin{remark}
    The star-formula also holds for sticky star-polytopes, but we cannot confirm its non-negativity in this case. Its non-negativity follows from the $\ell_C$-positivity Conjecture \ref{non_neg_local_sticky_con}.
\end{remark}

\begin{cor}\label{thin_corollary}The Thin criterion from Theorem \ref{Arnold_monotonicity_criterion_theorem} (2) follows immediately from the Star-formula.
\end{cor}

Note that if $T$ is a star-triangulation then the Links-formula is trivial, i.e. consists of only one summand $\nu_C(P)$. The Merged formula (3) is obtained by substituting the Star-formula (1) into the Links-formula (2). We prove the Star-formula in \S \ref{star_formula_subsec} and the Links-formula in \S \ref{Links_formula_subsec}.

\section{Proof of the non-negative analogue of the Kouchnirenko formula} \label{proof_of_non_neg_section}

\subsection{Proof of the Star-formula  Theorem \ref{Non-negative_analogue_theorem} (1)}\label{star_formula_subsec}


Note that if $0 \le \D \in T$ then $ \sum\limits_{\substack{0 \le \tl \D \le \D\\ \tl\D\in T}} Cap(\tl \D) = Vol(\D)$. This implies that 
$$\nu_C(P) = \sum \limits_{\text{$\D\in T$ is $V_C$}} (-1)^{n-\dim \D} Vol(\D) = \sum \limits_{\text{$\D\in T$ is $V_C$}} (-1)^{n-\dim \D} \sum\limits_{\substack{O \le \tl \D \le \D\\ \tl\D\in T}} Cap(\tl \D) = $$

$$= \sum\limits_{O \le \tl\D\in T} Cap (\tl \D) \sum\limits_{\substack{\tl \D \le \D \in T\\ \text{$\D$ is $V_C$}}} (-1)^{n-\dim \D} = \sum\limits_{O\le \tl \D \in T} Cap(\tl \D) \cdot \ell_{C/\D}(T/\tl \D;1)$$

\subsection{Proof of the Links-formula Theorem \ref{Non-negative_analogue_theorem} (2)} \label{Links_formula_subsec}

The proof of the Links-formula consists of three steps. First step is the already proved Decomposition lemma \ref{Decomposition_lemma} which states that the links of $V_{\min}$-simplices of any triangulated sticky polytope decompose the polytope into sticky polytopes. Second step is the Additive lemma \ref{Additive_lemma} which states that the Newton number is additive with respect to sticky polytopes decomposition. The third step is the Links projection lemma \ref{link_projection_lemma} which says that the Newton number $\nu_C(P_T(\D))$ of the link of a $V_{\min}$-simplex $\D$ equals the product $Vol(\D) \cdot \nu_{\mathsmaller{C/D}} (P_T(\D)) / \D$ of its volume and the Newton number of the link projected along the simplex. Together these three steps obviously provide the Links-formula. Let us prove them separately.

\subsubsection{Additive lemma for the Newton number of sticky polytopes}

Recall that the Newton number of polytopes behaves the same way as the Euler characteristic, namely (see \cite{O89}), if $P = P_1 \cup P_2$, then we have
\begin{equation}\label{O-add}
    \nu_C(P) = \nu_C(P_1) + \nu_C(P_2) - \nu_C(P_1 \cap P_2)
\end{equation}
In this subsection we prove that the Newton number is additive with respect to decompositions on sticky polytopes.

\begin{definition}\label{sticky_and_star_decomposition_def}
    \emph{Sticky polytopes decomposition} of polytope $P$ is a set of sticky polytopes with no intersecting interiors such that their union is the polytope $P$.
\end{definition}

\begin{lemma}[Additive lemma]\label{Additive_lemma}
The Newton number is additive for sticky polytopes decompositions.
Namely, suppose that $\mathcal S$ is a sticky polytopes decomposition of a sticky polytope $P$ in the cone $C$. Then we have

$$\nu_C(P) = \sum\limits_{P_i\in \mathcal S} \nu_C(P_i)$$
\end{lemma}

\begin{proof}

Formula (\ref{O-add}) on p.\pageref{O-add} implies the inclusion-exclusion formula for the Newton number:
$$
\nu_C(P) = \sum_{\mathcal S^\prime \subset \mathcal S}(-1)^{\# \mathcal S^\prime - 1} \cdot \nu_C(\bigcap\limits_{P_i\in \mathcal S^\prime} P_i)
$$
Note that the polytopes in the sticky polytopes decomposition $\mathcal S$ intersect only in its boundaries and, by the definition of sticky polytopes, their boundaries cannot afford a non-zero Newton number. This implies that if $\# \mathcal S^\prime > 1$ then $\nu_C(\bigcap\limits_{P_i\in \mathcal S^\prime} P_i) = 0$. So the inclusion-exclusion formula implies the additivity for sticky polytopes decompositions.

\end{proof}

\subsubsection{Link projection lemma}

In this subsection we give a formula for the Newton number of the links of $V_{\min}$-simplices. Recall that by $P_T(\D)$ we denote the union of the simplices in lk$_T(\D)$ (see Notation \ref{P_T_sign}).



\begin{lemma}[Link projection lemma]\label{link_projection_lemma}
    Consider a triangulation $T$ of a sticky polytope $P$ in the cone $C$ and a $V_{\min}$-simplex $\D\in T$. Then we have

    $$\nu_C(P_T(\D)) = Vol(\D) \cdot \nu_{C/\D} (P_T(\D)/\D)$$
    
\end{lemma}

\begin{proof}

Note that there is one-to-one correspondence of $V$-simplices in lk$_T(\D)$ and $V$-simplices in lk$_T(\D) / \D$ and the lattice volume of the simplex before the factorisation is equal to $Vol(\D)$ times 
the lattice volume of the factorized simplex. This consideration implies Link projection lemma.
\end{proof}

\subsubsection{Example of Links projection formula dimension 2}\label{sticky_pol_in_dim2_section}

\begin{figure}[H]
		\begin{center}
		\includegraphics[scale=1.5]{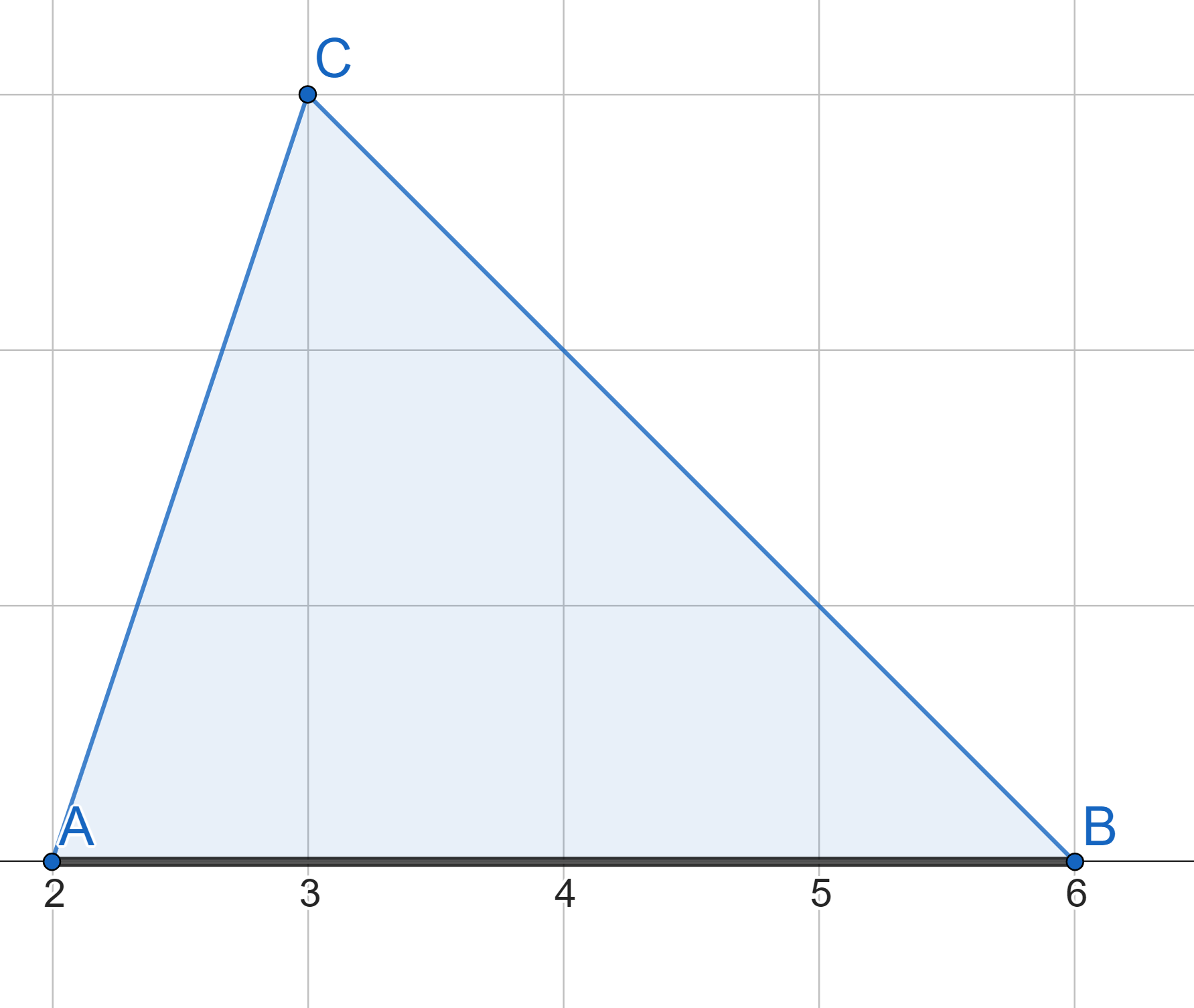}
		\end{center}
    \caption{\label{Sticky2d} $V$-triangle $ABC$ and its $V_{\min}$-segment $AB$}
\end{figure}
The triangulation consisting of one triangle provides the following calculation 
$$\nu(ABC) = Vol(AB) \cdot \nu(\pi_{\mathsmaller{AB}} (ABC)) =  4 \cdot(3-1) = 8$$

\subsubsection{Example of Links projection formula dimension 3}
\begin{figure}[H]
		\begin{center}
		\includegraphics[scale=0.5]{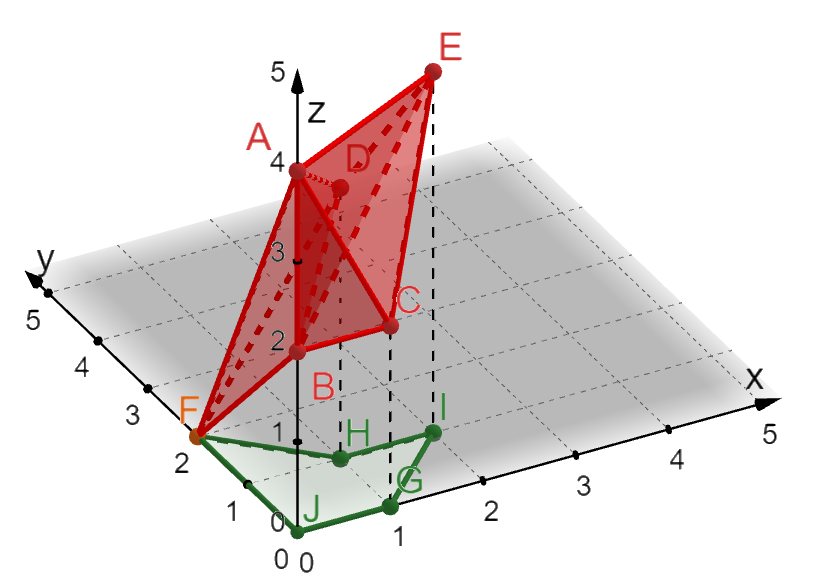}
		\end{center}
    \caption{\label{3dlink} Link $P_3$ (red) of the $V_{\min}$-segment $AB$ and its projection $JGIHF$ (green)}
\end{figure}
The Link projection formula implies the following calculation 
$$\nu(P_3) = Vol(AB) \cdot \nu(JGIHF) = 2 \cdot (4-1-2+1) = 4$$

\subsubsection{Example: sticky triangulation or good pavement}

Let us give an explicit formula for the Newton number of sticky triangulation. This case is very close to the case of \emph{good pavement} from \cite{ELT22}, so we also call it a good pavement.

Suppose that $\D \subset C\subset \R^n$ is a sticky $n$-simplex. Then $\D$ contains the unique $V_{\min}$-face $m(\D)$ and (since $\D$ is sticky) all faces of $\D$ containing $m(\D)$ are $V$-faces and all the others are not. Then the projection $\pi_{m(\D)}(\D)$ of the simplex $\D$ along aff$(m(\D))$ is a simplex of dimension $n - \dim (m(\D))$ with vertices in the origin and on the coordinate rays. Denote the heights of the points on the coordinate rays by $h_1(\D),\dots ,h_{n- \dim (m(\D))} (\D)$. Then we have 
$$\nu_C(\D) = Vol(m(\D)) \cdot\prod\limits_{i=1}^{n-\dim (m(\D))} (h_i (\D) - 1)$$
and this provides the following non-negative formula for the simplicial sticky polytopes decompositions.

\begin{cor}\label{Good_pavement_cor}
    If a triangulation $T$ of a sticky polytope $P$ in a cone $C$ is a good pavement, then the Newton number of $P$ is given by the following non-negative formula $$\nu_C(P) = \sum\limits_{\D\in T:\ \dim \D = n} Vol(m(\D)) \cdot\prod\limits_{i=1}^{n-\dim (m(\D))} (h_i (\D) - 1)$$
\end{cor}
Note that if a star-polytopes admits a good pavement all whose simplices are $B_1$ then it is negligible. Here is a good pavement giving another proof that the Example \ref{B_2-facet_example} provides a solution of the Arnold's monotonicity problem.

\begin{exa} [continued]\label{continued_example}

The good pavement of the polytope $\Gamma_{+G} \setminus \Gamma_+$ is the pyramid with the apex $G$ over the triangulation of the polytope $\Pi_G(\Gamma)$ on the right Figure \ref{B_2-picture} on p.\pageref{B_2-picture}. The good pavement consists of the following $3$ sticky simplices 
$$\Gamma_{+G} \setminus \Gamma_+ =  \quad GXYAP \quad \cup \quad GAYPQ \quad \cup \quad GAYQB$$
$$\nu(\Gamma_{+G} \setminus \Gamma_+) = Vol(GXY) \cdot (1-1) \cdot (1-1) + Vol(GPQY) \cdot (1-1) + Vol(GABY) \cdot (1-1)  = 0$$

All the simplices in the triangulation are $B_1$-simplices so the Newton number is preserved by Corollary \ref{Good_pavement_cor}.
\end{exa}

Note that there may not exist a good pavement even for the difference of convenient Newton polyhedra in $\R^3_{\ge 0}$, but we do not know examples of negligible complete star-polytopes  which do not admit a good pavement.

\section{Arnold's monotonicity problem in arbitrary dimension}

In this section we prove results concerning the Arnold's monotonicity problem in arbitrary dimension. First of all, we complete the proof of Theorem \ref{Arnold_monotonicity_criterion_theorem} by proving the Negligible criterion Theorem \ref{Arnold_monotonicity_criterion_theorem} (1) (the Thin criterion \ref{Arnold_monotonicity_criterion_theorem} (1) is already proved, see Corollary \ref{thin_corollary}). The Negligible criterion reduces the Arnold's monotonicity problem to the classification of complete negligible star-polytopes, while the Thin criterion is a criterion for the complete negligible star-polytopes. Then we demonstrate how Theorem \ref{Arnold_monotonicity_criterion_theorem} works by giving a short proof of the Arnold's monotonicity problem in dimension $3$. Then we prove the Negligible example lemma \ref{Negligible_examples_lemma} providing wide class of examples of complete negligible star polytopes and give an example in dimension $5$ which cannot be obtained by this lemma. Finally, we formulate and prove the Sprig lemma \ref{sprig_lemma} which generalizes for arbitrary dimension the observation from Figure \ref{Crook} when negligible polytopes in the cone $\R^2_{\ge 0} \oplus \R^1$ arise as parts of the negligible polytopes in the cone $\R^3_{\ge 0}$.

\subsection{Proof of the Negligible criterion Theorem \ref{Arnold_monotonicity_criterion_theorem} (1)} \label{Negligible_criterion_subsection}

In this subsection we give a criterion for the Arnold's monotonicity problem.
Consider a convenient Newton polyhedron $\Gamma_+ \subset \R^n_{\ge 0}$ and a lattice point $G\in \Gamma_-$. Denote by $C_G$ the cone $\pi_{\mathsmaller{OG}}(\R^n_{\ge 0})$. Denote by $\Pi_G(\Gamma) \subset C_G$ the polytope $\pi_{\mathsmaller{OG}}(\Gamma_{+G} \setminus \Gamma_+)$.

\begin{definition}\label{top_triang_def} \emph{Top triangulation $T^t$} is a triangulation of $(\Gamma \setminus \Gamma_G)$. For a top triangulation $T^t$ denote by $T^t_\mathsmaller {G}$ the triangulation of $ \Gamma_{+G} \setminus \Gamma_+$ which is the pyramid over $T^t$ with the vertex $G$.
\end{definition}

\begin{definition}\label{bottom_triang_def}
\emph{Bottom triangulation $T^b$} is a triangulation of $\Gamma_G \setminus \Gamma$ such that all its $(n-1)$-dimensional triangles contain the point $G$.
\end{definition}

Both projections $\pi_\mathsmaller{OG} (T^t)$ and $\pi_\mathsmaller{OG} (T^b)$ are triangulations of  $\Pi_G(\Gamma)$ and the triangulation $\pi_\mathsmaller{OG} (T^b)$ is a star-triangulation, so $\Pi_G(\Gamma)$ is a complete star-polytope.

\begin{figure}[H]
		\begin{center}
		\includegraphics[scale=.25]{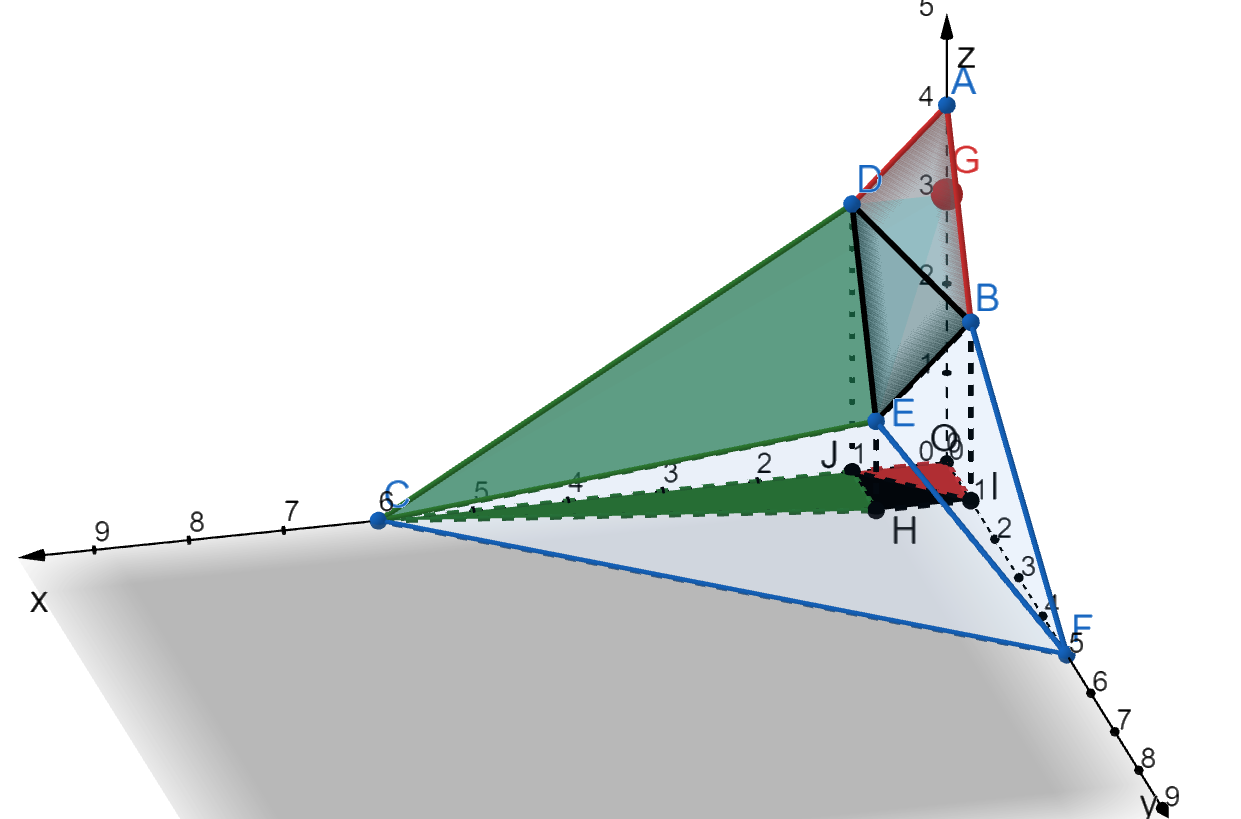}
		\includegraphics[scale=.25]{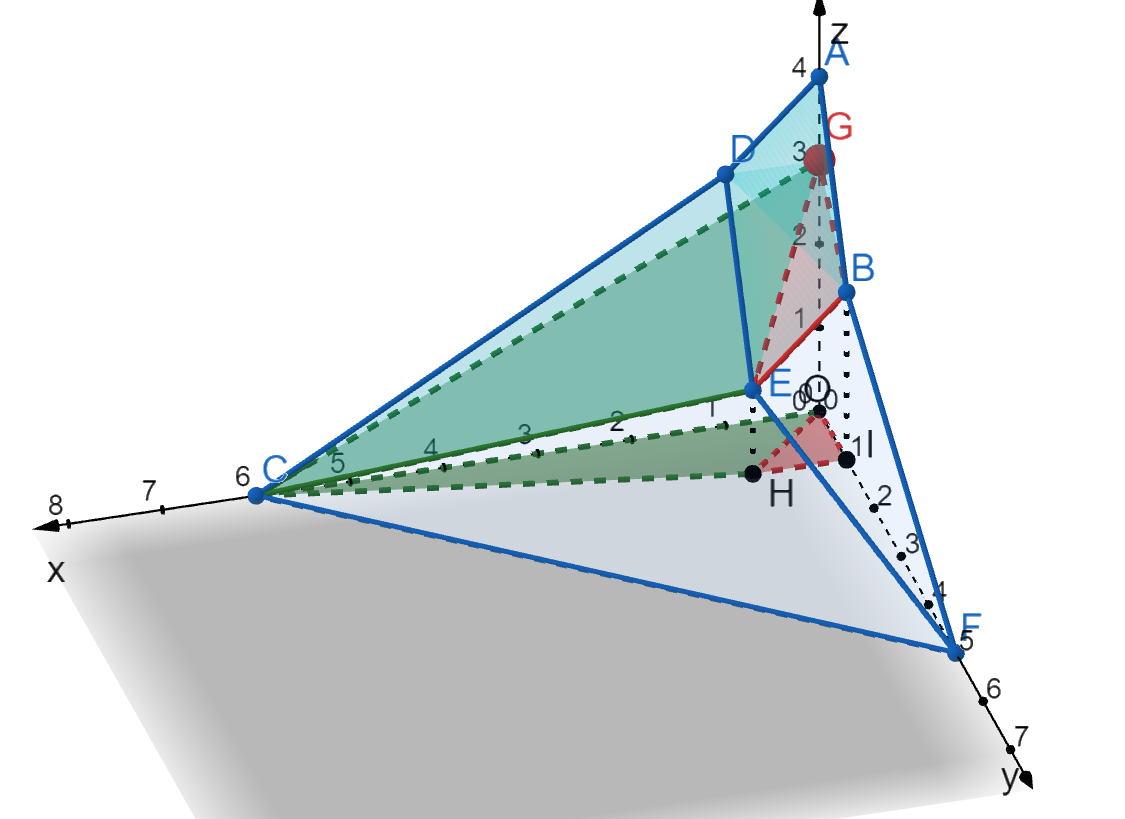}
		\end{center}
    \caption{\label{lower_upper_fig} \textbf{Top triangulation} $T^t$ (left) and its projection $\pi_\mathsmaller{OG} (T^t)$ and \textbf{Bottom triangulation} $T^b$ (right) and its projection $\pi_{\mathsmaller{OG}}(T^b)$} 
\end{figure}


\begin{proof}[Proof of Theorem \ref{Arnold_monotonicity_criterion_theorem} (1)]

From the Additive Lemma \ref{Additive_lemma} we have $$\nu(\Gamma_-) = \nu(\Gamma_{-G}) + \nu(\Gamma_{+G} \setminus \Gamma_+),$$
so it is sufficient to prove that the inequalities $\nu(\Gamma_{+G} \setminus \Gamma_+) > 0$ and $\nu_{C_G} (\Pi_G(\Gamma)) > 0$ are equivalent.

Consider a top-triangulation $T^t$ of $\Gamma \setminus \Gamma_G$. Denote by $\Sset$ the set of $V_{\min}$-simplices in the triangulation $T^t_G$ of $\Gamma_{+G} \setminus \Gamma_+$. Recall that the union of the simplices in the link $\text{lk}_{T^t_G} (\D), \D \in \Sset$ is denoted by $P_{T^t_G} (\D)$. For a simplex $\D^\prime \in T^t_G$ denote by $\tl \D^\prime = \pi_{\mathsmaller{OG}}(\D^\prime \cap \Gamma) \subset C_G$ the projection of the corresponding simplex of $T^t$. Note that a simplex $\D\in T^t_G$ is $V_{\min}$ if and only if the simplex $\tl \D \subset C_G$ is $V_{\min}$ and the projections $\pi_{\D} (P_{T^t_G}(\D)) = \pi_{\tl \D} (P_{\pi_\mathsmaller{OG} (T^t)} (\tl \D))$  are the same, hence the Links projection lemma \ref{link_projection_lemma} implies
 $$\nu(P_{T^t_G} (\D)) = \frac{Vol(\D)}{Vol(\tl \D)} \cdot \nu_{\mathsmaller{C_G}} (P_{\pi_{\mathsmaller{OG}} (T^t)} (\tl \D))$$

\begin{figure}[H]
		\begin{center}
		\includegraphics[scale=.3]{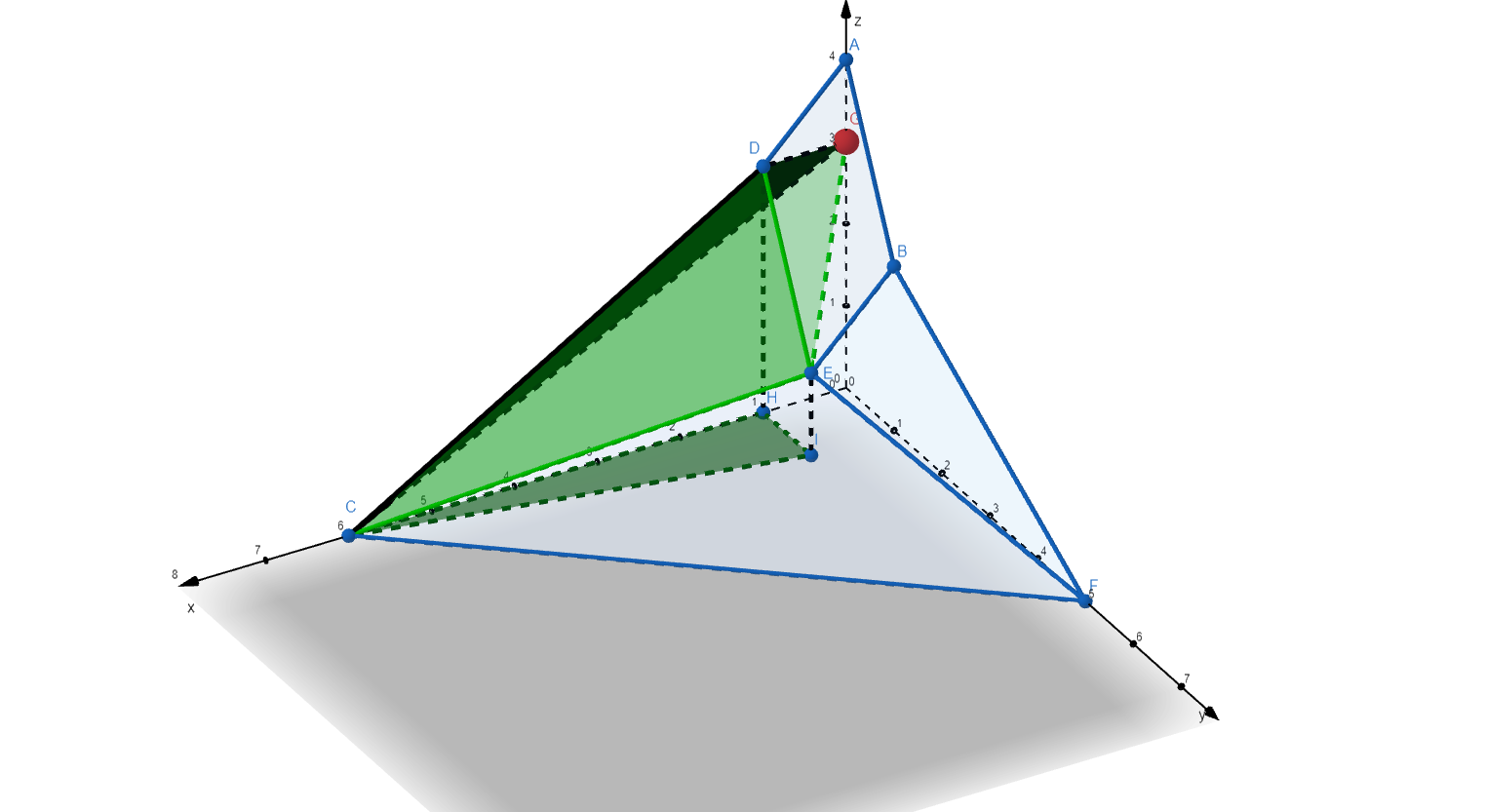}
		\end{center}
    \caption{\label{A_Dupdown_fig} $V$-simplex $\mathbf{\D}$ (black triangle), polytope $\mathbf{P_{T^t_G} (\D)}$ (light-green tetrahedron) and polytope $\mathbf{P_{\pi_\mathsmaller{OG} (T^t)} (\tl \D)}$ (dark-green triangle)}
\end{figure}

Note that the polytopes $P_{\pi_\mathsmaller{OG} (T^t)} (\tl \D), \D \in \Sset$ form a sticky polytopes decomposition of the polytope $\Pi_G(\Gamma)$ since they are sticky and have non-intersecting interiors. So due to the Additive lemma \ref{Additive_lemma} the condition $\nu_{\mathsmaller{C_G}} (\Pi_G(\Gamma)) > 0$ is equivalent to the existence of $V_{\min}$-simplex $\D \in T^t_G$ such that $\nu_{\mathsmaller{C_G}} (P_{\pi_\mathsmaller{OG} (T^t)} (\tl \D)) > 0$ which is equivalent to $\nu(P_{T^t_G} (\D)) > 0$. The latter condition is equivalent to the positiveness $\nu(\Gamma_{+G} \setminus \Gamma_+) > 0$ by the Additive Lemma \ref{Additive_lemma}. 
\end{proof}

\subsection{Example: proof of the Arnold's monotonicity problem in dimension 3} \label{3_Ar_pr_subsec}

In this subsection we give a short proof of the Arnold's monotonicity problem in dimension $3$ using the technique developed in this paper.

\begin{theorem-non}[\cite{BKW19}]\label{3_Arnold_monotonicity_th} 
        A pair $(\Gamma_+, G) \subset \R^3_{\ge 0}$ satisfies the Arnold's monotonicity problem if and only if the polytope $\Gamma_{+G} \setminus \Gamma_+$ is a $B_1$-pyramid.
\end{theorem-non}

Due to the Negligible criterion Theorem \ref{Arnold_monotonicity_criterion_theorem}(1), in order to solve this problem we need to classify negligible complete star-polytopes in dimension $2$.
    
\subsubsection{Negligible complete star-polytopes in dimension 2}\label{Negligible_classification_dim2_sec}

Note that there no negligible star-polytopes in $C = \R^n$. If $C = \R^1_{\ge 0}$ then the only negligible star-polytope is the segment $[0,1]$. Let us classify all negligible star-polytopes in dimension $2$.

\begin{enumerate}
    \item $\mathbf{C = \R^1_{\ge 0} \oplus \R^1}$: all negligible complete star-polytopes are $B_1$-pyramids.
    
\begin{figure}[H]
		\begin{center}
		\includegraphics[scale=1.5]{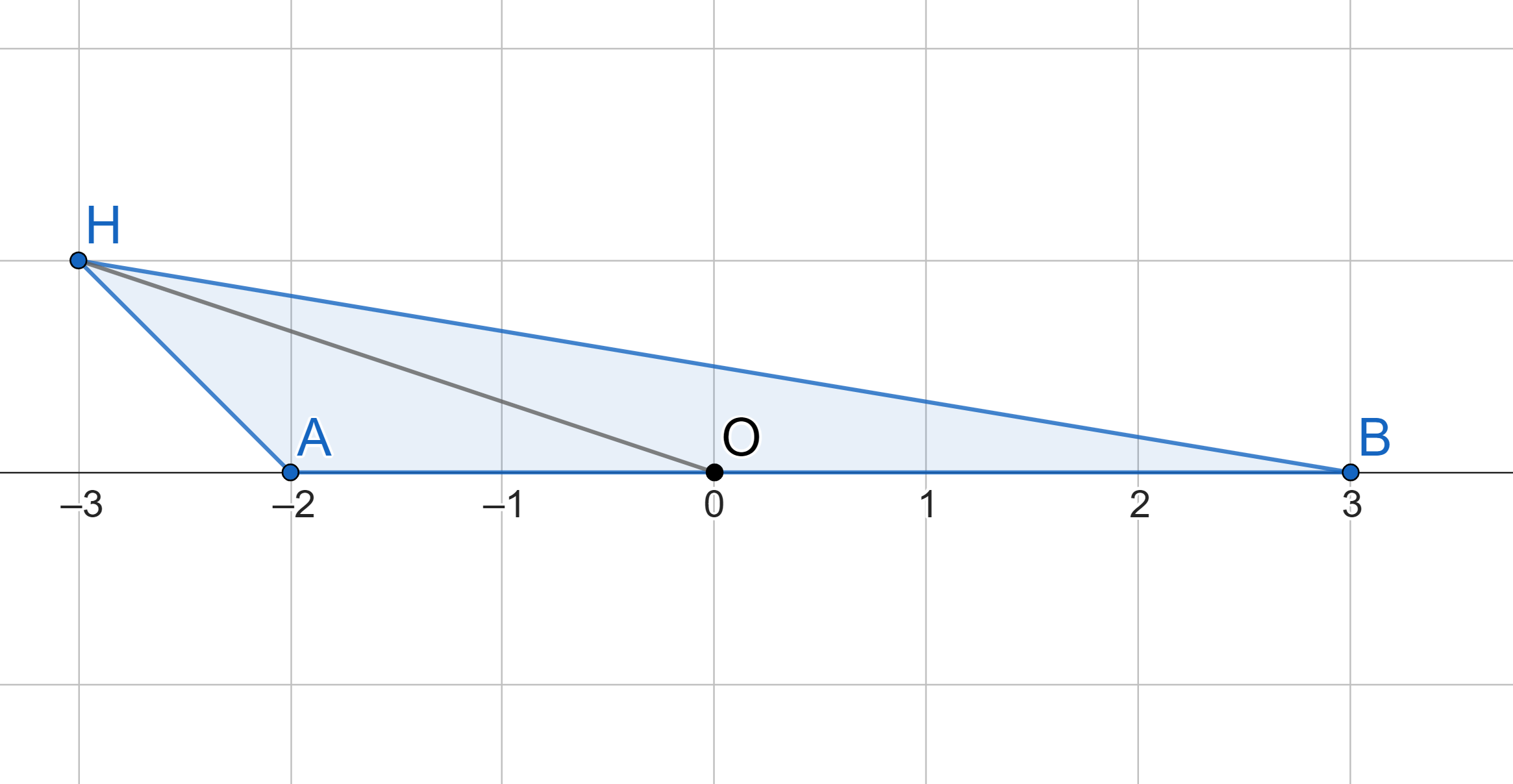}
		\end{center}
    \caption{\label{1+1} \textbf{Negligible} star-polytopes in the cone $C= \R^1_{\ge 0} \oplus \R^1$}
\end{figure}

    If the polytope $P \subset C$ with star-triangulation $T$ is not a pyramid, then the origin $O$ is not $0$-thin and the triangulation $T/O$ is not thin by Corollary \ref{internal_ray_cor}. So due to the Thin criterion $P$ is a pyramid.

    The Newton number of a pyramid with height $h$ equals $(h - 1) \cdot Vol(AB)$, so $h = 1$ and $P$ is $B_1$-pyramid. This fact can also be proved using the Thin criterion. Namely, if $h > 1$, then by the Pick's formula the triangle $ABH$ contains a lattice point $Q \notin AB, Q \ne H$. If $Q \in OC$, then $OH$ is not $0$-thin and $T/OH$ is not thin. Otherwise $Q$ is in the $Box^\circ$ of one of two triangles $OAH$ or $OBH$ and the corresponding triangle is not $0$-thin.
    
    \item $\mathbf{C = \R^2_{\ge 0}}$: all negligible complete star-polytopes are $B_1$-pyramids.

\begin{figure}[H]
		\begin{center}
		\includegraphics[scale=1.5]{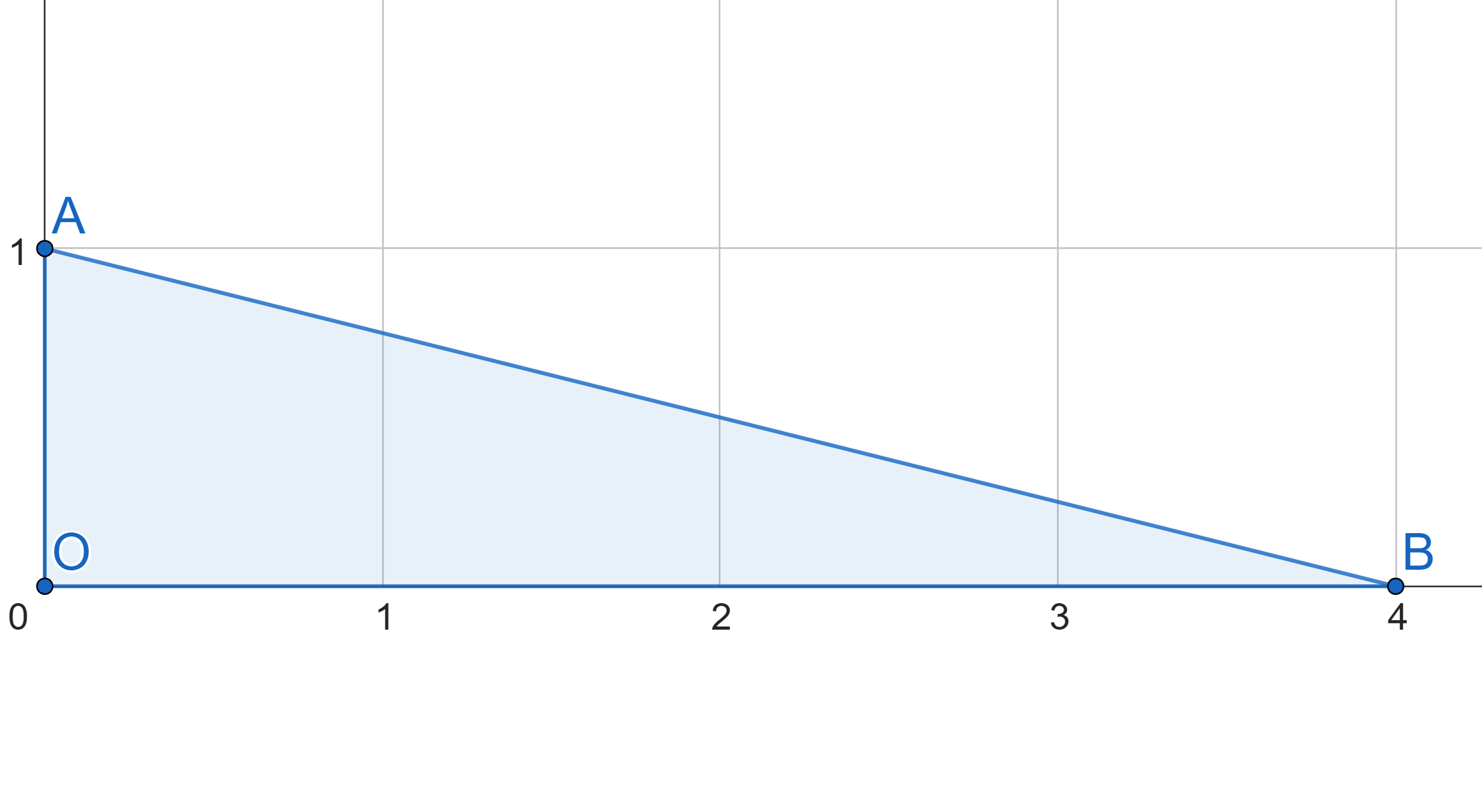}
		\end{center}
    \caption{\label{2+0} \textbf{Negligible} star-polytopes in the cone $C= \R^2_{\ge 0}$}
\end{figure}
    
    It is sufficient to prove that if the negligible star-polytope $P$ contains a lattice point $Q$ outside the coordinate cross, then $\nu(P) > 0$. Suppose there is such point $Q$, then consider a (not necessarily star-) triangulation $T$ of the star-polytope $P$ containing the segment $[O,Q]$. Then $\ell(T/O;1) > 0$ by Corollary \ref{internal_ray_cor}. Thus $\nu(P)\ge \ell(T/O;1) > 0$ by the Merged formula \ref{Non-negative_analogue_theorem} (3).
    
\end{enumerate}

\subsubsection{Proof: Arnold's monotonicity problem in dimension 3}

Due to the Negligible criterion Theorem \ref{Arnold_monotonicity_criterion_theorem} (1), the Arnold's monotonicity problem in dimension $3$ is implied by the classification of negligible star-polytopes in dimension 2 from \S \ref{Negligible_classification_dim2_sec}. Namely, there are 3 cases depending on the position of the point $G$ with respect to the positive octant $\R^3_{\ge 0}$. Considering these cases prove that the only solutions of the Arnold's monotonicity problem in dimension $3$ are $B_1$-pyramids.

\begin{enumerate}
    \item \textit{Minimal face of $\R^3_{\ge 0}$ containing $G$ is of dimension $3$.}

    In this case $C_G = \R^2$ and there are no negligible complete star-polytopes.
    
    \item \textit{Minimal face of $\R^3_{\ge 0}$ containing $G$ is of dimension $2$.}

    In this case $C_G = \R^1\oplus \R^1_{\ge 0}$ and all negligible complete star-polytopes are $B_1$-pyramids.

    \item \textit{Minimal face of $\R^3_{\ge 0}$ containing $G$ is of dimension $1$.}

    In this case $C_G = \R^2_{\ge 0}$ and all negligible complete star-polytopes are  $B_1$-pyramids.
\end{enumerate}

\subsection{Proof of the Negligible example lemma \ref{Negligible_examples_lemma}}

     
     The triangulation $T/O$ is thin because it is obtained from a trivial triangulation by a sequence of conical facet refinements. Let us prove that for any simplex $\D\in T: \D > 0$ either $\D$ is $B_1$-pyramid or $P_T(\D)$ is a pyramid (i.e. contains only one vertex outside some hyperplane). Note that in the first case the simplex $\D$ is $0$-thin and in the second the localized triangulation $T/\D$ is thin, so, by the Thin criterion \ref{Arnold_monotonicity_criterion_theorem}(2), this statement implies the desired Lemma \ref{Negligible_examples_lemma}.

    Consider a simplex $\D\in T: \D > 0$. Suppose there exists $i$ such that $W_i \in \D$ and $\D \subset F_i$. Then $\D$ is $B_1$-pyramid with vertex $W_i$ and base on $H_i$. 
    
    If there is no such $i$ then consider the first moment $j\ge 1$ in the sequence of conical facet refinements when the simplex $\D$ arises. Denote the triangulation at this moment by $T_j$. Note that $P_{T_j}(\D)$ is a pyramid with vertex $W_j$ and that the triangulations $T/D$ and $T_j/D$ coincide since a neighbourhood of any interior point of the simplex $\D$ remains the same during each further step in the sequence of conical facet refinements. So $T/D$ also contains only one ray outside a hyperplane, thus $P_T(\D)$ is a pyramid.
    
        

\begin{remark}
    Negligible star-polytope obtained by the Negligible example lemma \ref{Negligible_examples_lemma} is $B_1$-polytope if and only if it can be obtained by at most $1$ conical facet refinement.
\end{remark}

Note that all negligible star-polytopes in dimension up to $2$ are obtained by Lemma \ref{Negligible_examples_lemma}. In the next section we prove that it also holds in dimension $3$. We believe that it is also true in dimension $4$ (it is implied by the Thin triangulations Conjecture \ref{thin_triang_con}). The simplest example in dimension $5$ we know which is not predicted by Lemma \ref{Negligible_examples_lemma} is as follows.

\subsection{Example outside the Negligible example lemma \ref{Negligible_examples_lemma} in dimension 5}

\begin{ex}[Further example] \label{Unpredicted_example} \ 

\begin{figure}[H]
		\begin{center}
		\includegraphics[scale=1]{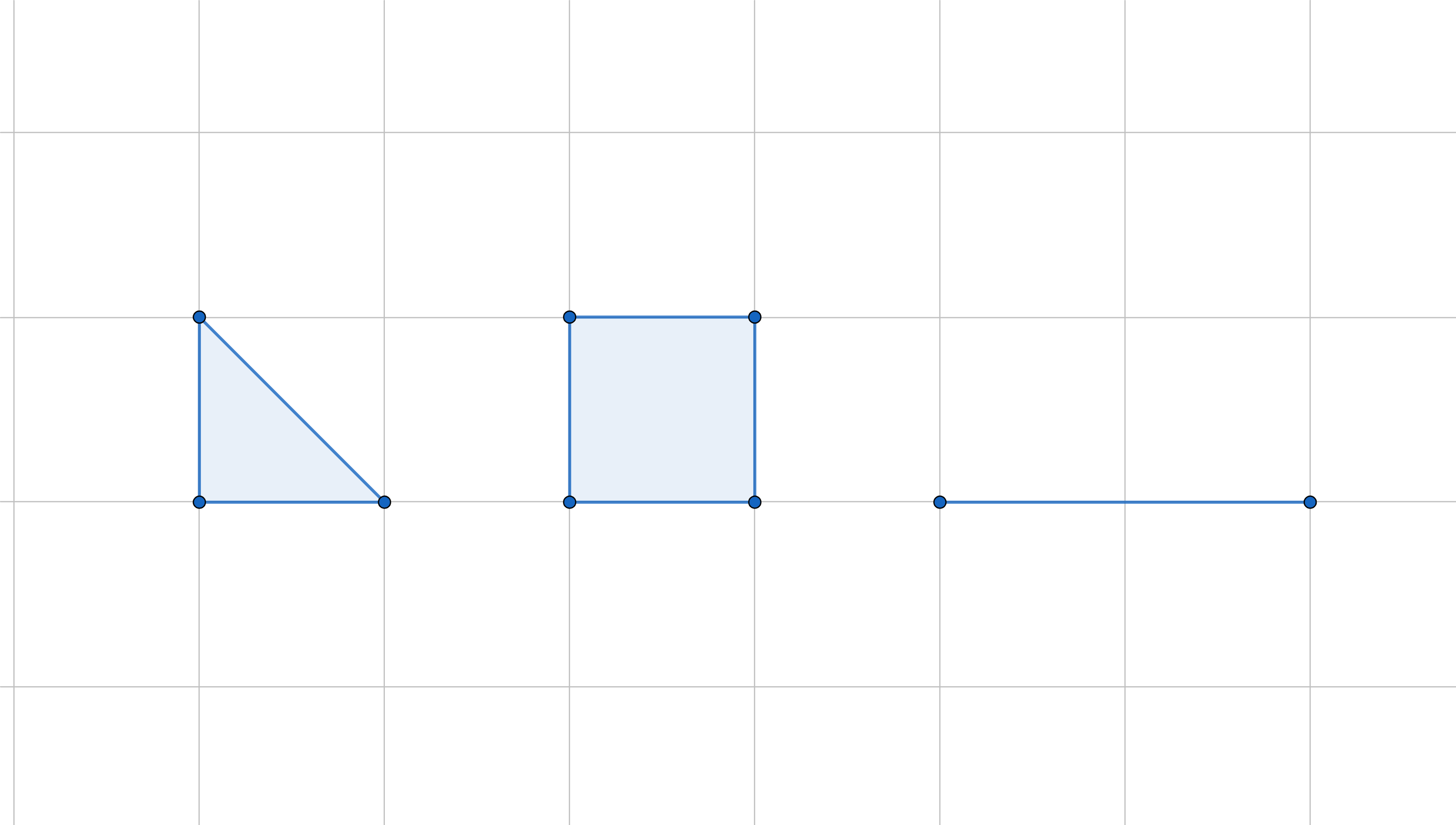}
		\end{center}
    \caption{\label{P_5} The triangle $\vartriangleright$, the square $\Box$ and the segment $[-1,1]$}
\end{figure}

    Consider the cone $C = \R^4_{\ge 0} \oplus \R^1$
    and the polytope 
    $$P_5 = (\vartriangleright_{1,2} \oplus\ \Box_{3,4}) + [-1,1]_5,$$
    where polygon $\vartriangleright_{1,2}$ is the unit triangle $$\vartriangleright_{1,2} = Conv((0,0), (1,0),(0,1)) \subset \R^2_{\ge 0}$$ in the coordinate plane of $\R^5$ spanned by the first two coordinates, polygon $\Box_{3,4}$ is the unit square $$\Box_{3,4} = Conv((0,0),(0,1),(1,0),(1,1))\subset \R^2_{\ge 0}$$ in the coordinate plane of $\R^5$ spanned by the second two coordinates, the segment $[-1,1]_5$ is on the last coordinate ray of $\R^5$, $X\oplus Y$ is the \emph{free sum}, i.e. the convex hull of $X$ and $Y$ lying in the orthogonal spaces, $X+Y$ is the Minkowski sum of $X$ and $Y$ lying in the orthogonal spaces.

    Let us prove that this example does not satisfy Lemma \ref{Negligible_examples_lemma}. Note that if a triangulation $\Sigma$ of $\R^m_{\ge 0} \oplus \R^{n-m}, n-m>0$ is obtained from a trivial triangulation of $C_m$ (see Definition \ref{Trivial_fan_def}) then the triangulation $\Sigma$ also contains the simplicial cone $C_m$ since every proper face of $C_m$ is contained in more then one facets and the cone $C_m$ is not a $V$-cone. So there is a cone $C_m$ in $\Sigma$ such that $\pi_{\mathsmaller{\R^{n-m}}} (C_m) = \R^m_{\ge 0}$. But there is no star-triangulation $T$ of $P_5$ with simplex $\D_4 \in T$ such that the cone of $\pi_{\mathsmaller{\R^1}} (\D_4) \subset \R^4_{\ge 0}$ is the whole $\R^4_{\ge 0}$.

    The sums of the lattice volumes of the intersections of the free sum $\vartriangleright_{1,2} \oplus\ \Box_{3,4} \subset \R^4_{\ge 0}$ with the coordinate subspaces of the corresponding dimensions are as follows $$2 \cdot (1,4\cdot 1, 5 \cdot 1 + 2, 2\cdot 2 + 2\cdot 1, 2) = 2 \cdot (1,4,7,6,2)$$ Then the same vector for $P_5$ is as follows $$2 \cdot (1\cdot 1, 4\cdot 2, 7\cdot 3, 6\cdot 4, 2\cdot 5) = 2 \cdot (1,8,21,24,10)$$ and the alternating sum $2\cdot(1-8+21-24+10) = 0$, so $P_5$ is indeed a negligible star-polytope in the cone $C$.
    
\end{ex}

\begin{remark}
    If we replace the segment $[-1,1]$ with $[0,1]$ and the cone $\R^4_{\ge 0} \oplus \R^1$ with $\R^5_{\ge 0}$ then the polytope will be also negligible. But we can not verify whether this example can be obtained by Lemma \ref{Negligible_examples_lemma} or not.
\end{remark}

\subsection{Negligible sprigs of Negligible polytopes}\label{Sprigs_subsecion}

Recall that in the classification from \S \ref{4-Arnold_intro_sec} the negligible star-polytopes in the cone $\R^3_{\ge 0}$ contain the negligible star-polytopes (called greenhouses) in the cone $\R^2_{\ge 0} \oplus \R^1$. In this subsection we observe the same phenomena for the cones $\R^m_{\ge 0} \oplus \R^{n-m} \subsetneq \R^k_{\ge 0}\oplus \R^{n-k}, k<m$, i.e.a connection between negligible polytopes in different cones.

\subsubsection{Sprig Lemma}

Let us start with an observation which follows directly form the Star-formula Theorem \ref{Non-negative_analogue_theorem} (1).

\begin{claim}[ReConing principle]\label{reConing_claim}
Consider pair of cones $C = \R^m_{\ge 0} \oplus \R^{n-m}, C^\prime = \R^{m^\prime}_{\ge 0} \oplus \R^{n-m^\prime}$, a complete star-polytope $P \subset C$ with star-triangulation $T$ and a sticky star-polytope $P^\prime \subset P, P^\prime \subset C^\prime$ whose star-triangulation $T^\prime$ consists of some simplices of the triangulation $T$. Suppose that for any $\D \in T^\prime$ we have $\ell_C(T/\D;1) \ge \ell_{C^\prime} (T^\prime/ \D;1)$. Then $\nu_C(P) \ge \nu_{C^\prime} (P^\prime)$.
\end{claim}

Note that in case the $\ell_C$-positivity conjecture \ref{non_neg_local_sticky_con} is confirmed, then the ReConing claim also holds if $P$ is strongly contractible. It would also help to confirm the conditions of the ReConing claim. Applying the ReConing principle \ref{reConing_claim} to an arbitrary sharp cone containing the simplex $\D$ in its interior we obtain the following Corollary.

\begin{cor}\label{Vol=1_cor}
    Suppose that $T$ is a star-triangulation of some negligible complete star-polytope in a cone $C = \R^m_{\ge 0} \oplus \R^{n-m}$ and $\D\in T$ is a simplex.
    If for every its face $O<0 \tl \D \le \D \in T$ we have $\ell_{C/\tl \D}(T/\tl \D;1) > 0$, then $Vol(\D) = 1$.
\end{cor}

Denote by $C^\D, \D \subset C$ the cone $C + \text{aff} (\D)$. Note that if $C$ is of the form $\R^{m}_{\ge 0} \oplus \R^{n-m}$ then $C^\D$ is of the same form (but $m$ may change). The following lemma is a useful specification of the ReConing principle.

\begin{lemma}[Sprig lemma]\label{sprig_lemma}
Consider a complete star-polytope $P\subset C = \R^m_{\ge 0} \oplus \R^{n-m}$ with star-triangulation $T$. Suppose that $\D\in T$ is a $V_C$-simplex and there is only one $n$-simplex $\tl \D \in T$ containing $\D$. Denote by $R$ the minimal simplex such that $O\le R\le \tl \D$ and $Conv(\D,R) = \tl \D$.
Suppose that for any simplex $\D^\prime \ge R$ we have $\ell_{C/\D^\prime} (T/\D^\prime;1) > 0$, then $\nu_C(P) \le \nu_{C^\D}(\tl \D)$. For instance, if $P$ is negligible, then $\tl \D$ is $B_1$-simplex in $C^\D$.
\end{lemma}

\begin{proof}
    Note that  $\ell_{C^\D / F}(\tl \D/ F;1) = \begin{cases}
      0, F \ngeq R\\
      1, F \ge R
    \end{cases}$, so this lemma is a corollary of the ReConing principle \ref{reConing_claim}.
\end{proof}

\begin{remark}
    If the simplex $R$ is interior (i.e. not contained in any proper face of $C$), then the conditions $\ell_{C/\D^\prime} (T/\D^\prime;1) > 0$ (for any simplex $R \le \D^\prime \le \tl \D$) of the Sprig lemma \ref{sprig_lemma} are automatically applied.
\end{remark}

\begin{remark}
    Applying a sequence of conical facet refinements along the simplex $\tl \D$ (with respect to the cone $C^\D$) preserves the Sprig lemma \ref{sprig_lemma}.
\end{remark}

The latter Remark explains the naming of ``Sprig lemma'' and why some negligible polytopes in $\R^m\oplus \R^{n-m}$ arise as parts of negligible polytopes in $\R^k\oplus \R^{n-k}$ for $k>m$ (see Figure \ref{sprig_figure} below).

\begin{figure}[H]
		\begin{center}
		\includegraphics[scale=.3]{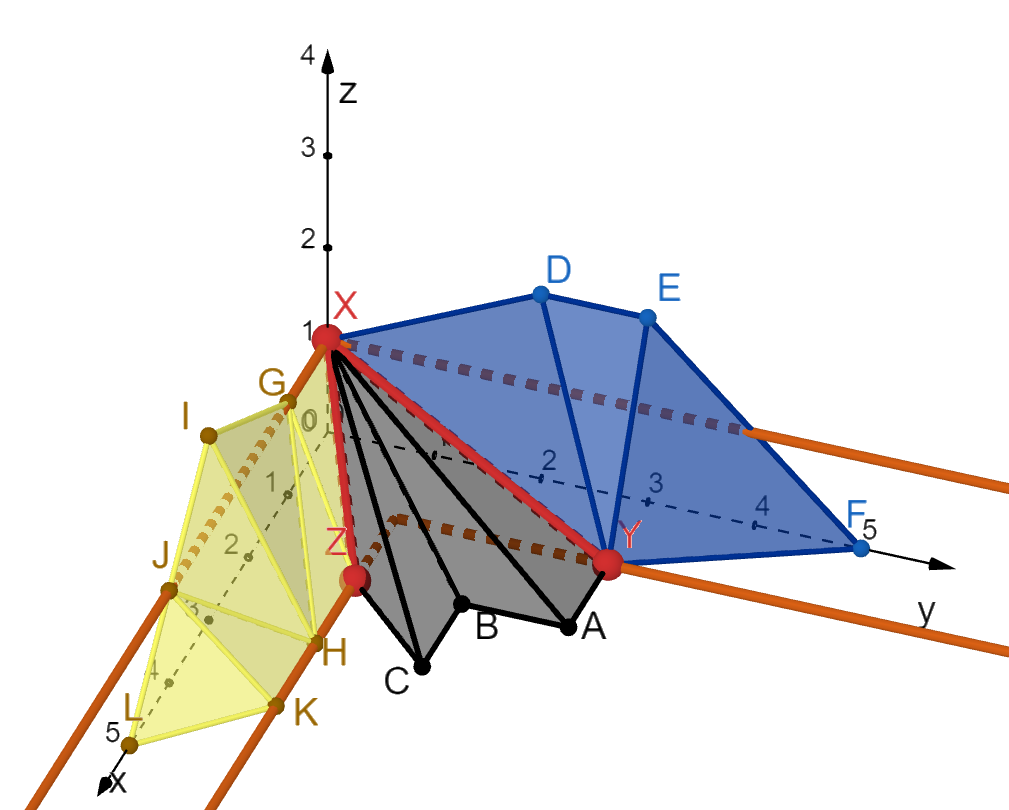}
		\includegraphics[scale=5]{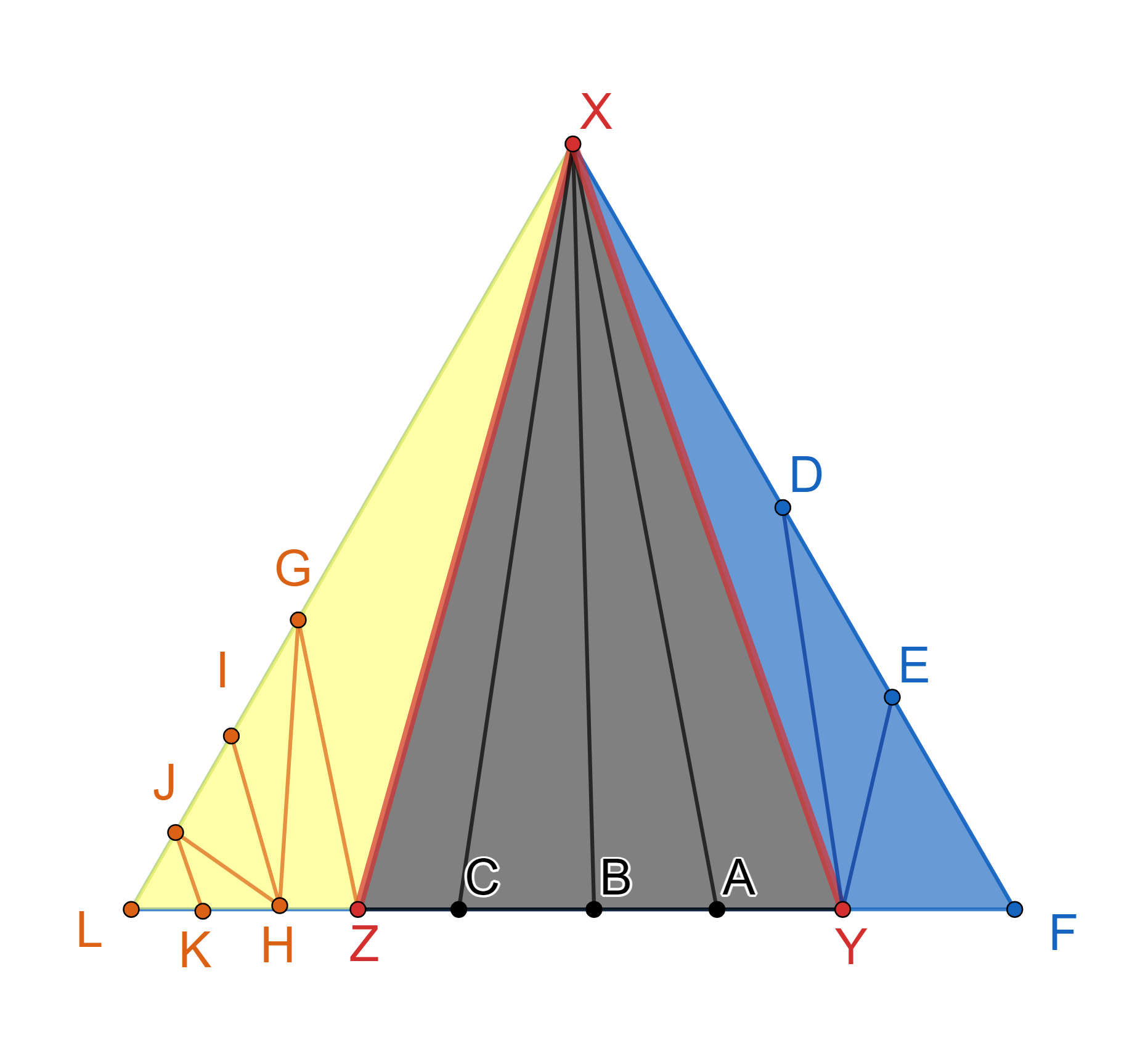}
		\end{center}
    \caption{\label{sprig_figure} Yellow and blue sprigs for $\tl \D = XYFO; \D = OL$ and $\tl \D = XZLO; \D = OF$. The sprigs are greenhouses.}
\end{figure}

\section{Arnold's monotonicity problem in dimension 4 and 5}\label{4_5_Ar-sec}

In this section we prove Theorem \ref{4_5_Arnold_theorem}. Together with the Negligible Criterion Theorem \ref{Arnold_monotonicity_criterion_theorem} (1) it provides complete solution of the Arnold's monotonicity problem in dimension $4$ and partial solution of the Arnold's monotonicity problem in dimension $5$, namely, for the case when the additional point $G$ is on a coordinate ray. It also provides complete solution in dimension $5$ if the Thin triangulation conjecture \ref{thin_triang_con} is applied.

Theorem \ref{4_5_Arnold_theorem} (and Remark below it concerning complete solution in dimension $5$) is a corollary of the following two Lemmas.

\begin{figure}[H]
		\begin{center}
		\includegraphics[scale=5]{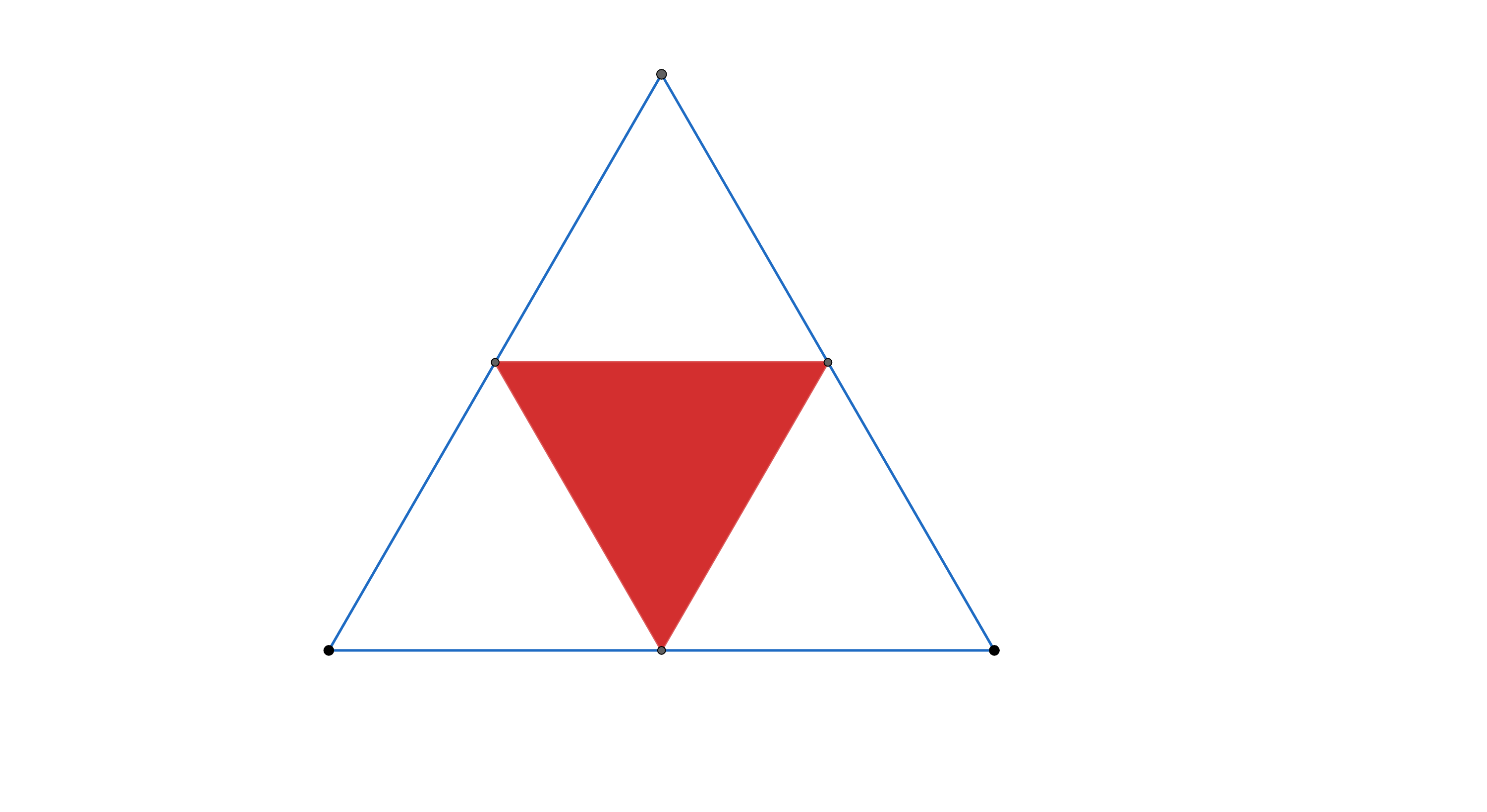}
		\end{center}
    \caption{\label{sentral_triforc} Projectivisations of the triforce and the \textbf{``central'' cone} (red)}
\end{figure}

\begin{lemma}[No triforce lemma]\label{sur_no_triforce_lem}
    Suppose that $P \subset \R^3_{\ge 0}$ is a negligible complete star-polytope with star-triangulation $T$. Then the triangulation $T/O$ does not contain a ``central cone'', i.e. a cone with three rays in the interior of different $2$-dimensional faces of $\R^3_{\ge 0}$ (see Figure \ref{sentral_triforc}).
\end{lemma}

\begin{lemma}\label{4_5_Arnold_lemma}
    Consider a negligible complete star-polytope $P$ with star-triangulation $T$ in a cone $\R^m_{\ge 0} \oplus \R^{n-m}$ for $n \le 4$. If the triangulation $T/O$ is obtained from a trivial triangulation by a sequence of conical facet refinements then the star-polytope $P$ can be obtained by the Negligible example lemma \ref{Negligible_examples_lemma}.
\end{lemma}

\begin{remark}
    Note that there are cases when given sequence of conical facet refinements (with fixed data $\{F_i, H_i\}$) does not apply Lemma \ref{Negligible_examples_lemma}. But the polytope $P$ can still be obtained by Lemma \ref{Negligible_examples_lemma} with a bit modified data $\{F_i, H_i\}$.
\end{remark}

The latter Lemma \ref{4_5_Arnold_lemma} is already proved in dimension up to $2$. We prove it in dimension $3$ and $4$ in \S \ref{complete_4_ar_sec} and \S \ref{5-Ar-sec} respectively. First we need to obtain some properties of Negligible star-polytopes in dimension $3$ and $4$. No triforce lemma \ref{sur_no_triforce_lem} is a corollary of Lemma \ref{no_triforce_lem} (see Corollary \ref{no_triforce_corollary} below).

\subsection{Properties of negligible star-polytopes in dimension 3 and 4}

In this subsection we prove some lemmas concerning negligible complete star-polytopes in dimension $3$ and $4$. The proofs significantly rely on the classification from \S \ref{Classification_of_vanish_subsec} and we believe that similar statements are not true in higher dimension. Using these lemmas we prove that Lemma \ref{Negligible_examples_lemma} provides all solutions to the Arnold's monotonicity problem in dimension $4$ and (under the assumption that the additional point is on a coordinate ray) in dimension $5$.

\begin{lemma}\label{no_triforce_lem}
    Consider any simplex $\D$ containing the origin whose cone is the ``central cone'' of a triforce. Then $Vol(\D) \ge 2$.
\end{lemma}

\begin{proof}
    $$Vol(\D) =  det
\begin{bmatrix}
   a & b & 0\\
   0 & c & d\\
   e & 0 & f
\end{bmatrix} = acf + bde \ge 2$$

\end{proof}

\begin{cor}\label{no_triforce_corollary}
    No triforce lemma \ref{sur_no_triforce_lem} is a corollary of Lemma \ref{no_triforce_lem}  and Corollary \ref{Vol=1_cor}.
\end{cor}

\begin{lemma}\label{B_1-simplex_lem}[$B_1$-simplex lemma]
    Consider some negligible complete star-polytope in dimension $n \le 4$ with star-triangulation $T$ in a cone $C = \R^m_{\ge 0} \oplus \R^{n-m}$. Then any $n$-simplex $\tl \D \in T$ containing a $V_C$-facet is $B_1$-simplex.
\end{lemma}

\begin{remark}
    We cannot prove this lemma for $n = 4$ and $\tl \D \in T$ containing no $V_C$-facets.
\end{remark}

\begin{proof}

    If the conditions of Sprig lemma \ref{sprig_lemma} (for $\D$ equals the intersection of all the $V_C$-facets of $\tl \D$) are applied then Lemma \ref{B_1-simplex_lem} follows from Sprig lemma \ref{sprig_lemma}. Comparing with the classification from \S \ref{Classification_of_vanish_subsec} we obtain that the only case when the conditions of Sprig lemma \ref{sprig_lemma} are not applied is the following one: 
    
    \begin{itemize}
        \item The cone $C = \R^4_{\ge 0}$
        \item The simplex $\tl\D$ contains exactly one $V_C$-facet $\D$
        \item The vertex $W: W \in \tl \D, W \notin \D$ is on a coordinate ray of $C$
        \item The triangulation $T/OW$ is obtained from the triforce with ``central'' tetrahedron $\tl \D / OW$ by a sequence of conical facet refinements
    \end{itemize}
    
    Suppose that $\tl \D$ is not $B_1$-simplex. By Lemma \ref{no_triforce_lem} the lattice volume $Vol(\D) \ge 2$, so there is a simplex $O<\D^\prime\le\D$ which is not $0$-thin. Then for the simplex  $Conv(\D^\prime, W)$ is also not $0$-thin since the point $W$ is on a coordinate ray and has height at least $2$. Note that (by the classification from \S \ref{Classification_of_vanish_subsec}) for any $O<\D^\prime\le\D$ the triangulation $T/ Conv(\D^\prime, W)$ is not thin. So by the Thin criterion Theorem \ref{Arnold_monotonicity_criterion_theorem} (2) our polytope is not negligible. Thus, $\tl \D$ is a $B_1$-simplex.
\end{proof}

\subsection{Theorem \ref{4_5_Arnold_theorem}(1): solution of the Arnold's problem in dimension $4$}\label{complete_4_ar_sec}

Complete solution of the Arnold's monotonicity problem in dimension $4$ (Theorem \ref{4_5_Arnold_theorem}(2)) is implied by No triforce lemma \ref{sur_no_triforce_lem}, Lemma \ref{4_5_Arnold_lemma} in dimension $3$ and the classification of thin triangulations in dimension $3$ from \S \ref{Classification_of_vanish_subsec}. No triforce lemma \ref{sur_no_triforce_lem} is proved in Corollary \ref{no_triforce_corollary}.

Let us prove Lemma \ref{4_5_Arnold_lemma} for all the cones $C = \R^m_{\ge 0} \oplus \R^{3-m}$ separately.

\begin{enumerate}
    \item \textbf{Lemma \ref{4_5_Arnold_lemma} holds for the cone $C = \R^3_{\ge 0}$}




$B_1$-simplex lemma \ref{B_1-simplex_lem} implies that every simplex $F_i, i > 1$ in the sequence of conical facet refinements is $B_1$-simplex with the corresponding base $H_i$ (see Figure \ref{Proving_triangles}). Note that for $i = 1$ the vertex $W_1$ is on a coordinate ray and $\ell(T/OW_1;1) > 0$, so  $OW_1$ is $0$-thin and the height of $W_1$ equals $1$. This implies that the star-polytope $P$ is obtained by Lemma \ref{Negligible_examples_lemma}.

\begin{figure}[H]
		\begin{center}
		\includegraphics[scale=5]{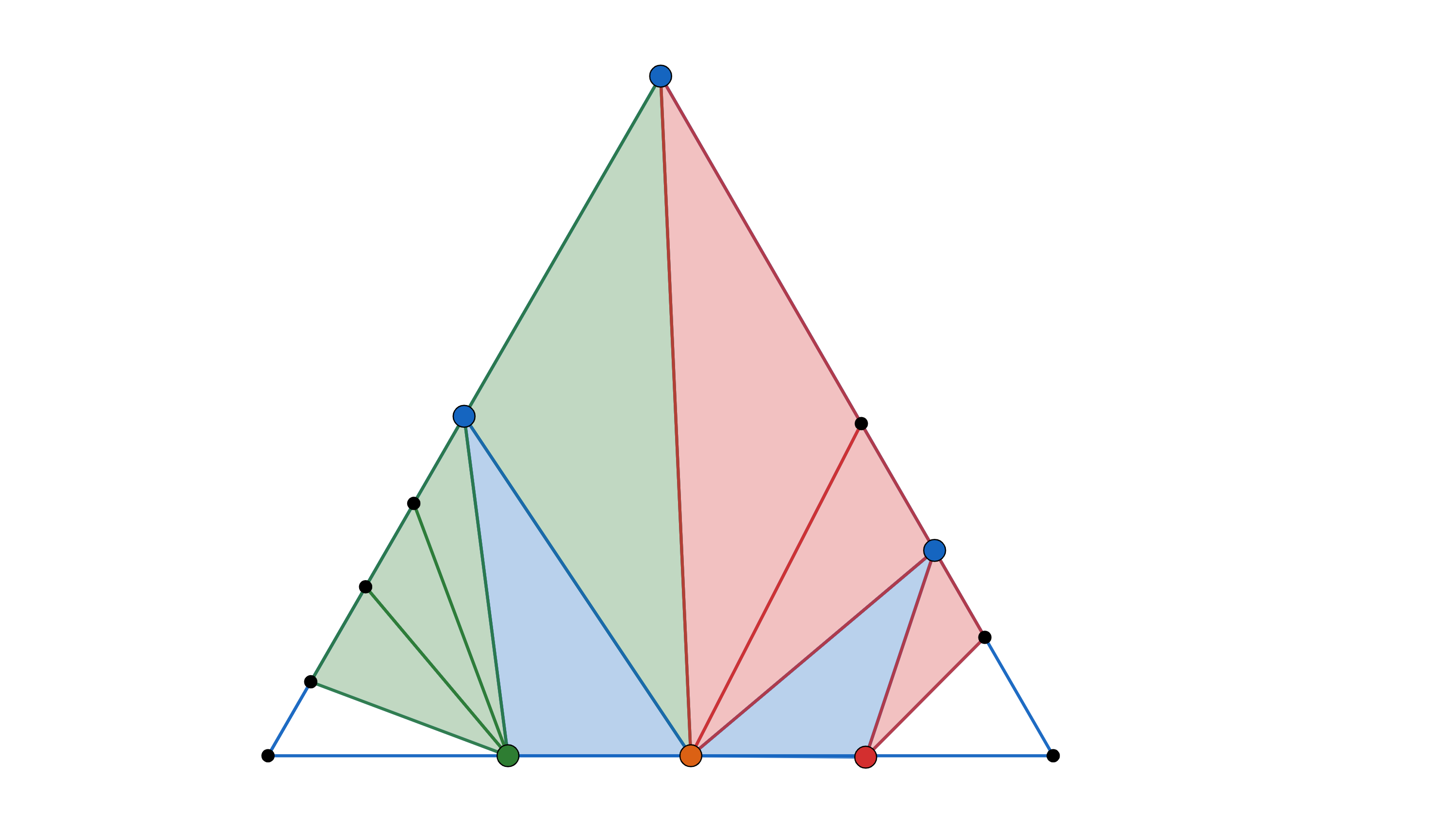}
		\end{center}
    \caption{\label{Proving_triangles} Lemma \ref{B_1-simplex_lem} is applied to blue, green and red triangles}
\end{figure}

\item \textbf{Lemma \ref{4_5_Arnold_lemma} holds for the cone $C = \R^2_{\ge 0} \oplus \R^1$}


\begin{figure}[H]
		\begin{center}
		\includegraphics[scale=2]{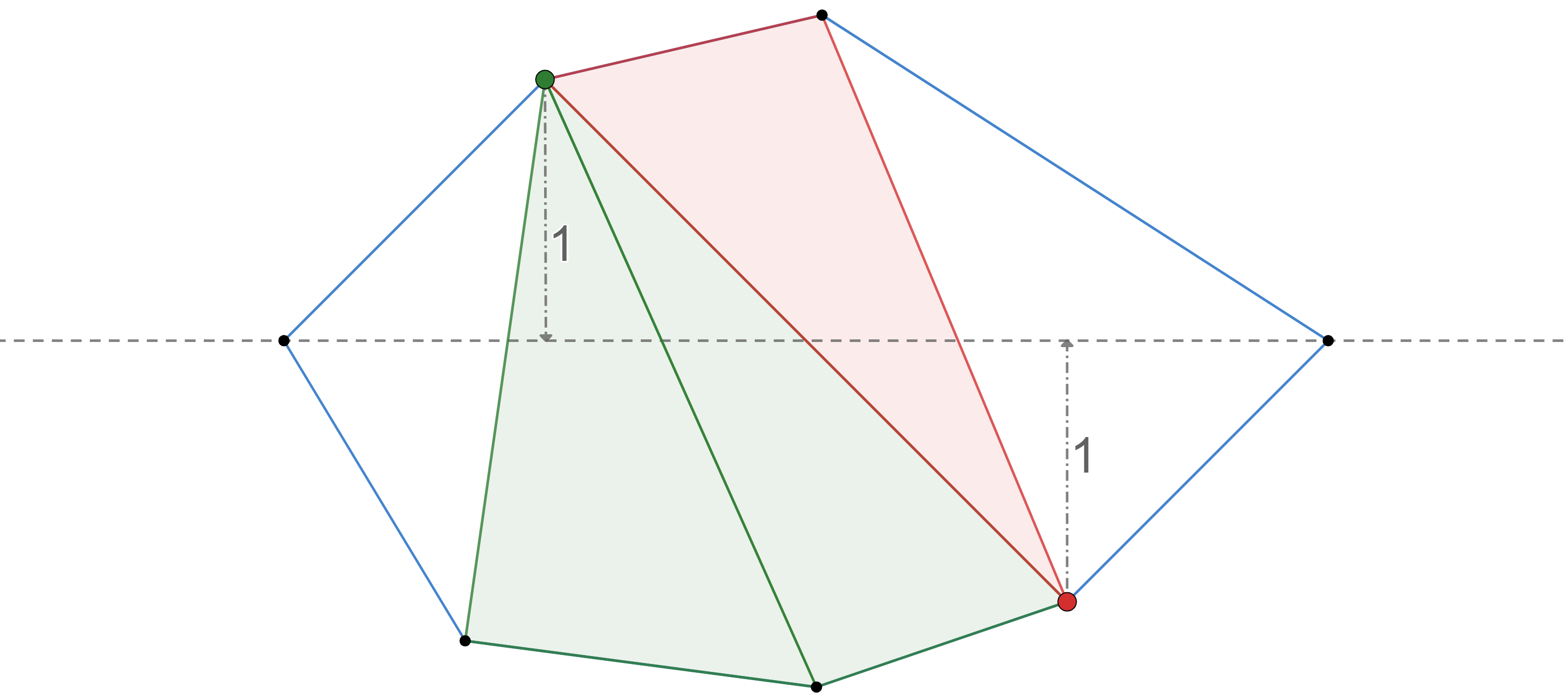}
		\end{center}
    \caption{\label{proj_triv_2+1} Projectivisation of a sequence of conical facet refinements of a trivial triangulation of the cone $C = \R^2_{\ge 0} \oplus \R^1$. Red vertices are of height $1$.}
\end{figure}

From $B_1$-simplex Lemma \ref{B_1-simplex_lem} we obtain that every vertex which lies in the interior of one of two facets of $\R^2_{\ge 0} \oplus \R^1$ and is connected with at least two points in the interior of the other facet has height $1$ with respect to the other facet. This implies that the star-polytope $P$ is obtained by Lemma \ref{Negligible_examples_lemma}.

\item \textbf{Lemma \ref{4_5_Arnold_lemma} holds for the cone $C = \R^1_{\ge 0} \oplus \R^2$}


Note that there is only one vertex $W$ of the negligible star-polytope $P$ outside $\R^2$. So we have $$\nu_C(P) = (h-1) \cdot Vol(P\cap \R^2),$$
where $h$ is the height of $W$ with respect to $\R^2$, so $h = 1$ and the polytope $P$ satisfies the Remark under Lemma \ref{Negligible_examples_lemma}.
\end{enumerate}

\subsection{Theorem \ref{4_5_Arnold_theorem}(2): partial solution to the Arnold's problem in dimension $5$} \label{5-Ar-sec}

Solution of the Arnold's monotonicity problem in dimension $5$ in the case if the additional point is on an coordinate subspace (Theorem \ref{4_5_Arnold_theorem}(2)) is implied by  Lemma \ref{4_5_Arnold_lemma} in dimension $4$ and the classification of thin triangulations in dimension $3$ from \S \ref{Classification_of_vanish_subsec}. Recall that if the Thin triangulation conjecture \ref{thin_triang_con} will be confirmed it will imply the complete solution in dimension $5$.

\begin{proof}[Proof of Lemma \ref{4_5_Arnold_lemma} in dimension $4$]

Let us prove that the triangulation $T/O$ of $C = \R^m_{\ge 0}
\oplus \R^{4-m}$ is obtained (maybe by a different) sequence of conical facet refinements along facets $F_i, i = 1,\dots, k$ such that each corresponding vertex $W_i$ has height $1$ with respect to the corresponding codimension face $H_i$ (as in Negligible example lemma \ref{Negligible_examples_lemma}).

\begin{enumerate}
    
    \item Assume there is a step $i$ in the sequence of conical facet refinements such that the corresponding facet $F_i$ contains exactly one codimension one  $V_C$-face. 
    
    Then after any sequence of conical facet refinements there will remain a top-dimensional simplex with exactly one $V_C$-facet and apex $W_i$. So by $B_1$-simplex Lemma \ref{B_1-simplex_lem} the height $\rho(W_i,H_i)$ is equal to $1$.
    
    \item Assume there is a step $i$ in the sequence of conical facet refinements such that the corresponding facet $F_i$ contains exactly two codimension one $V_C$-faces.
    
    Suppose that $\rho(W_i,H_i) > 1$ then by $B_1$-simplex Lemma \ref{B_1-simplex_lem} we have $\rho(W_i^\prime,H_i^\prime) = 1$, where $H_i^\prime$ is the second codimension one $V_C$-face of $F_i$ and $W_i^\prime$ is the corresponding vertex. There can not be a further conical facet refinement with number $j \ge i$ with apex $W_j = W_i$ such that the splitting of $F_i$ contains a facet with exactly one codimension one $V_C$-face because of the previous item. So all the further conical facet refinements with indices $j\ge i$ such that $W_j = W_i$ and $H_j\subset H_i$ can be considered as conical facet refinements with $W_j = W_i^\prime$ and $H_j \subset H_i^\prime$.
    
    \item Assume there is a step $i$ in the sequence of conical facet refinements such that the corresponding facet $F_i$ contains exactly three codimension one $V_C$-faces.
    
    Suppose that $\rho(W_i,H_i) > 1$ then here can not be a further conical facet refinement with number $j \ge i$ with apex $W_j = W_i$ such that the splitting of $F_i$ contains a facet with exactly one codimension one $V_C$-face because of the first item.
    
    Denote by $\D$ the triangle face of $P\cap H_i$ which does not contain the origin. Denote by $S<\D$ its segment whose cone is not a $V_C$-face. The only case when there is no $4$-simplex with exactly one codimension one $V_C$-face and base on $H_i$ after the conical facet refinement is when the segment $S$ is separated from the rest triangle $H_i^p$ by a line (as in the left Figure \ref{3v_case}).

   \begin{figure}[H]
		\begin{center}
		\includegraphics[scale=5]{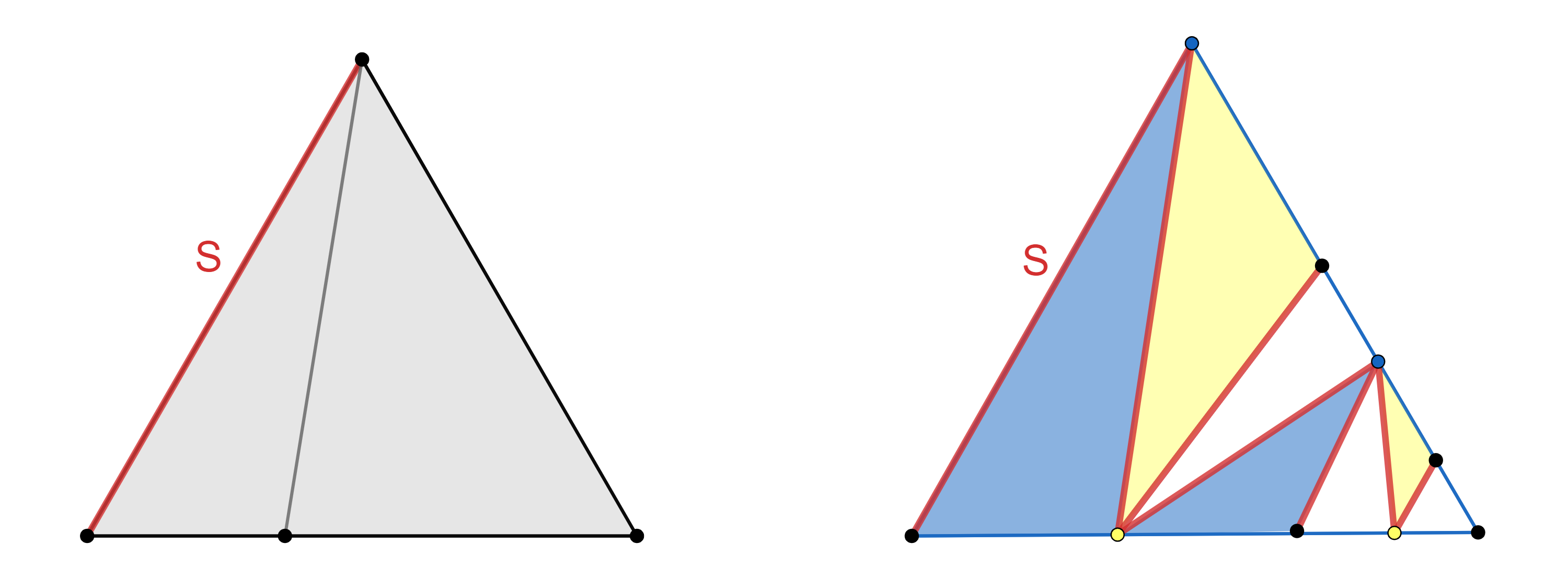}
		\end{center}
    \caption{\label{3v_case} Triangle $\D$ with the segment $S$ (left); sequence of conical facet refinements (left), Lemma \ref{B_1-simplex_lem} is applied to blue and yellow triangles}
\end{figure}

    By the first item there cannot be a $4$-simplex with exactly one codimension one $V_C$-face and base on $H_i$ in the triangulation $T$. So the sequence of conical facet refinements touching $H_i$ is induced by a sequence of separations of the corresponding segment in the triangle $\D$ (see right Figure \ref{3v_case}). Due to $B_1$-simplex Lemma \ref{B_1-simplex_lem} each colored point in Figure \ref{3v_case} has height $1$ with respect to the corresponding codimension $1$ face of the cone $C$. Thus, the conical facet refinement can be considered as a sequence of conical facet refinements with apexes in the colored points.

    \item The only case left is when $i = 1$ and the cone $C = \R^4_{\ge 0}$.

    Suppose that $\rho(W_1, H_1) > 1$. Denote by $\D$ the triangular face of $H_1$ which does not contain the origin.

    The conical facet refinement induces a triangulation of the triangle $\D$. Note every non-trivial triangulation of a triangle satisfies one of the following items (see Figure \ref{4v_case}):

    \begin{enumerate}
        \item There is a vertex in the interior of $\D$.
        \item There is a segment connecting a vertex of the triangle and the interior of the opposite side.
        \item There is a triangle $\D^\prime$ in the triangulation all whose sides do not lie in the boundary of $\D$.
    \end{enumerate}

  \begin{figure}[H]
		\begin{center}
		\includegraphics[scale=4]{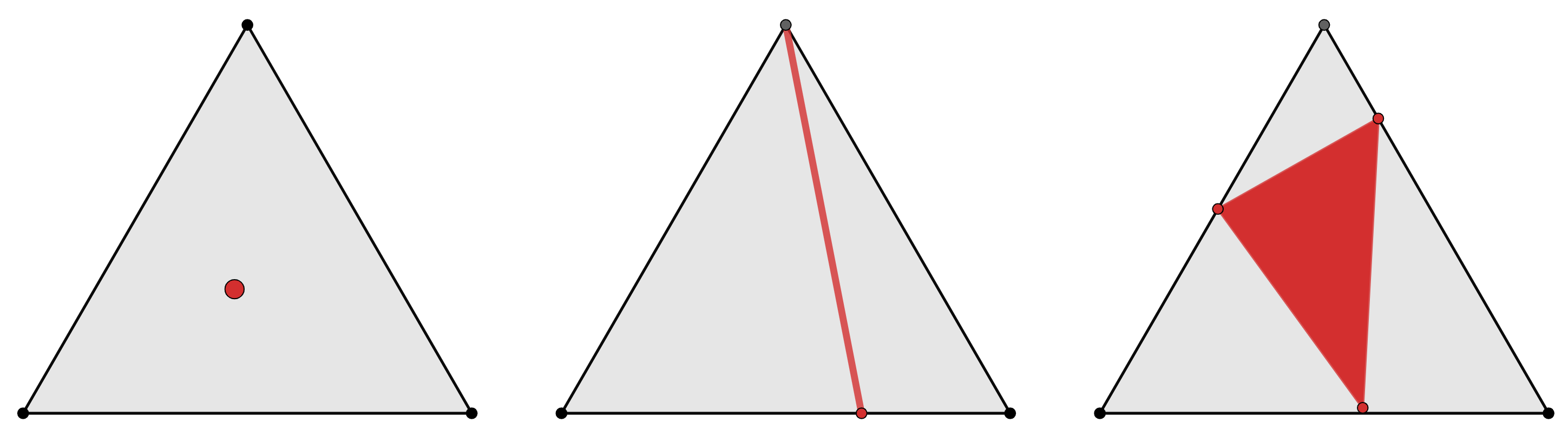}
		\end{center}
    \caption{\label{4v_case} Three possibilities of non-trivial triangulations of a triangle.}
\end{figure}

    Let us consider these cases separately.

    \begin{enumerate}
        \item If there is a point in the interior of $\D$, then by Corollary \ref{internal_ray_cor} we obtain that the triangulation $T/OW_1$ is not thin;  $OW_1$ is not $0$-thin by the assumption, thus $P$ is not negligible.

        \item If there is a segment $XR \subset \D$ connecting a vertex $X$ with a point $R$ on the opposite side, then by Corollary \ref{internal_ray_cor} the triangulation $T/OXW_1$ is not thin. Then $OXW_1$ is $0$-thin, so $X$ is the point $1$ on the corresponding coordinate ray. Thus we can modify the sequence of conical facet refinements so that the first one has apex $X$ and splits the simplex into two separated by the tetrahedron $OW_1XR$ and instead of the rest of the first conical facet refinement we consider new two conical facet refinements (along the cones of the two separated simplices). These new two conical facet refinements are considered in the item 3.

        \item In the last case we have the simplex $Conv(\D^\prime, W_1,O)$ with only one codimension one $V_C$-face, so by $B_1$-simplex Lemma \ref{B_1-simplex_lem} the height $\rho(W_1, H_1)$ is equal to one.
        
    \end{enumerate}
    
\end{enumerate}

\end{proof}

\vspace{1ex}
\noindent
\textit{Krichever Center for Advanced Studies}, Skolkovo Institute of Science and Technology, Moscow\\
\textit{Department of Mathematics}, National Research University ``Higher School of Economics'', Moscow \\
\textit{Email}: Fedor.Selyanin@skoltech.ru

\begin{thebibliography}{xxxxxxxxx}


\bibitem[Ar]{Ar} V.I. Arnold a.o., \textit{Arnold’s problems}. Springer and Phasis, 2004.

\bibitem[AGLV]{AGLV}  V. I. Arnold, V. V. Goryunov, O. V. Lyashko, V. A. Vasil'ev, \textit{Singularity theory I. Local and Global Theory}, Encyclopaedia of Mathematical Sciences, vol 6. Springer, Berlin, Heidelberg. 1993.

\bibitem[Ath12]{Ath12} C Athanasiadis.\textit{ Flag subdivisions and $\gamma$-vectors}, Pac. J. Math.,
259(2):257–278, 2012. \href{https://arxiv.org/abs/1106.4520}{arXiv:1106.4520}

        
\bibitem[Ath]{Ath} C. Athanasiadis, \textit{A survey of subdivisions and local h-vectors}, The mathematical
legacy of Richard P. Stanley, Amer. Math. Soc., Providence, RI, 2016, pp. 39–51. \href{https://arxiv.org/abs/1403.7144}{arXiv:1403.7144}

        



\bibitem[BR]{BR}  M. Beck, S. Robins, \textit{Computing the continuous discretely: Integer-point enumeration in
polyhedra}. Undergraduate Texts in Mathematics. Springer, New York, second edition, (2015). 
\href{https://matthbeck.github.io/papers/ccd.pdf}{https://matthbeck.github.io/papers/ccd.pdf}

\bibitem[BBD81]{BBD81} A. Beilinson, J. Bernstein, P. Deligne, \textit {Faisceaux pervers}, in Analysis and topology on
singular spaces, I (Luminy, 1981), 5–171, Ast\'erisque, 100, Soc. Math. France, Paris, 1982.

\bibitem[BM85] {BM85} U. Betke, P. McMullen, \textit{Lattice points in lattice polytopes}, Monatsh. Math. 99 (1985), no. 4, 253–265. MR 799674


\bibitem[BKN23]{BKN23}  C. Borger, A. Kretschmer, B. Nill, \textit{Thin polytopes: Lattice polytopes with vanishing local $h^*$-polynomial}, Int. Math. Res. Not. IMRN (2023). \href{https://arxiv.org/abs/2207.09323}{arXiv:2207.09323}

\bibitem[BKW19]{BKW19} S. Brzostowski, T. Krasiński, J. Walewska, \textit{Arnold's problem on monotonicity of the Newton number for surface singularities}. J. Math. Soc. Japan 71:4 (2019), 1257–1268 \href{https://arxiv.org/abs/1705.00323}{arXiv:1705.00323}

\bibitem[dCMM18]{dCMM18} M. de Cataldo, L. Migliorini, M. Mustat\u{a}, \textit{Combinatorics and topology of proper toric maps}, J. Reine Angew. Math. 744 (2018), 133–163, \href{https://arxiv.org/abs/1407.3497}{	arXiv:1407.3497}

\bibitem[DL92]{DL92}  J. Denef, F. Loeser, \textit{Caract\'eristiques d’Euler-Poincar\'e, fonctions z\^eta locales et modifications analytiques}, J.
Amer. Math. Soc. 5 (1992), no. 4, 705–720

		
		
		

\bibitem[Ehr62]{Ehr62} E. Ehrhart, \textit{Sur les poly\'edres rationnels homoth\'etiques\'a n dimensions},
C. R. Acad. Sci. Paris, 254:616–618, 1962.



\bibitem[EKK]{EKK} A. Esterov, B. Kazarnovskii, A. Khovanskii, \textit{Newton polytopes and tropical geometry}, Russian
Mathematical Surveys, 76(1):91, 2021. \href{https://www.math.toronto.edu/askold/2021\%20UMN\%201\%20English.pdf}{https://www.math.toronto.edu/askold/2021\%20UMN\%201\%20English.pdf}



\bibitem[ELT22]{ELT22}  A. Esterov, A Lemahieu, K. Takeuchi, \textit{On the monodromy conjecture for non-degenerate hypersurfaces}, Jour.
Eur. Math. Soc. (2022), \href{https://arxiv.org/abs/1309.0630}{arXiv:1309.0630}
		
\bibitem[Ew]{Ew} G. Ewald, \textit{Combinatorial Convexity and Algebraic Geometry}, New York: Springer-Verlag, (1996)

\bibitem[Ful]{Ful} W. Fulton, \textit{Introduction to toric varieties}, Ann. of Math. Stud. 131, The William H. Rover Lectures
in Geometry, Princeton Univ. Press, Princeton, NJ, 1993.

\bibitem[Fur04]{Fur04} M. Furuya, \textit{Lower Bound of Newton Number}, Tokyo J. Math. 27 (2004), 177-186. \href{https://arxiv.org/abs/math/9901107}{arXiv:math/9901107}

\bibitem[GKZ90]{GKZ90} I. M. Gelfand, A. V. Zelevinskii, M. M. Kapranov, \textit{Discriminants of polynomials in several variables and triangulations of Newton polyhedra}, Algebra i Analiz, 2:3 (1990),  1–62  (in Russian); Leningrad Math. J., 2:3 (1991), 499–505

\bibitem[GKZ94]{GKZ94} I. M. Gelfand, M. M. Kapranov, A. V. Zelevinsky, \textit{Discriminants, resultants, and multidimensional
determinants}, Math. Theory Appl., Boston, MA, 1994, x+523 pp.

\bibitem[KS16]{KS16} E. Katz, A. Stapledon, \textit{Local h-polynomials, invariants of subdivisions, and mixed Ehrhart theory}, Adv. Math.
286 (2016), 181–239. \href{https://arxiv.org/abs/1411.7736}{arXiv:1411.7736}

  

		
		
\bibitem[Kou76]{Kou76} A. Kouchnirenko, \textit{Polyèdres de Newton et nombres de Milnor}, Invent Math, 32:1, (1976), 1-31

    
\bibitem[L-AMS20]{L-AMS20} M. Leyton-\'Alvarez, H. Mourtada, M. Spivakovsky, \textit{Newton non-degenerate $\mu$-constant deformations admit simultaneous embedded resolutions}. Compositio Mathematica, 158, (2020) 1268 - 1297. \href{https://arxiv.org/abs/2001.10316}{arXiv:2001.10316}

\bibitem[LPS22a]{LPS22a} M. Larson, S. Payne, A. Stapledon, \textit{The local motivic monodromy conjecture for simplicial nondegenerate singularities}, (2022) \href{https://arxiv.org/abs/2209.03553}{arXiv:2209.03553}
    
\bibitem[LPS23]{LPS23}  M. Larson, S. Payne, A. Stapledon, \textit{Resolutions of local face modules, functoriality, and vanishing of local h-vectors}, Algebraic Combinatorics, (2023), \href{https://arxiv.org/abs/2209.03543}{arXiv:2209.03543}

\bibitem[Mac71]{Mac71} I. G. Macdonald, \textit{Polynomials associated with finite cell-complexes}, J. Lond. Math. Soc. (2) 4 (1971), 181–192.

\bibitem[Mil68]{Mil68} J. Milnor, \textit{Singular points of complex hypersurfaces}, Annals of Mathematics Studies, No. 61, Princeton University Press, Princeton, N.J., 1968.
  
\bibitem[dMGP+20]{dMGP+20} A. de Moura, E. Gunther, S. Payne, J. Schuchardt, A. Stapledon, \textit{Triangulations of simplices with vanishing local h-polynomial}, Algebr. Comb. 3 (2020), no. 6, 1417–1430. \href{https://arxiv.org/abs/1909.10843}{arXiv:1909.10843}
    
        

\bibitem[N12]{N12} B. Nill, \textit{Combinatorial questions related to stringy E-polynomials of Gorenstein polytopes}, in Toric Geometry (K.
Altmann et. al., eds.), 2012, pp. 62–64.

\bibitem[O89]{O89} M. Oka, \textit{On the weak simultaneous resolution of a negligible truncation of the Newton boundary}. In Singularities (Iowa City,
IA, 1986), volume 90 of Contemp. Math., pages 199–210. Amer. Math. Soc., Providence, RI, 1989.


\bibitem[Stan87]{Stan87} R. Stanley, \emph{Generalized H-vectors, intersection cohomology of toric varieties, and related results}, In Commutative algebra and combinatorics (Kyoto,
1985), volume 11 of Adv. Stud. Pure Math., pages 187–213. North-Holland, Amsterdam, 1987

\bibitem[Stan92]{Stan92} R. Stanley, \textit{Subdivisions and local h-vectors}, J. Amer. Math. Soc. 5 (1992), no. 4,
805–851.


\bibitem[Stan]{Stan} R. Stanley, \textit{Enumerative combinatorics 1}, Cambridge Univ. Press., Cambridge, UK, 1997.


\bibitem[Stap08]{Stap08} A. Stapledon, \textit{Weighted Ehrhart theory and orbifold cohomology}, Adv. Math. 219 (2008), no. 1, 63–88. \href{https://arxiv.org/abs/0711.4382}{arXiv:0711.4382}

\bibitem[Stap17]{Stap17} A. Stapledon, \textit{Formulas for monodromy}, Res. Math. Sci. 4 (2017), Paper No. 8, 42 pp. \href{https://arxiv.org/abs/1405.5355}{arXiv:1405.5355}


\bibitem[T]{T} K. Takeuchi, \textit{Geometric monodromies, mixed Hodge numbers of motivic Milnor fibers and Newton polyhedra}, (2023), \href{https://arxiv.org/abs/2308.09418}{arXiv:2308.09418}

\bibitem[Var76]{Var76} A. Varchenko, \textit{Zeta-function of monodromy and Newtons diagram}, Invent. Math. 37 (1976), 253–262.

    
\bibitem[Vas17]{Vas17} V. Vassiliev, \textit{Multiplicities of bifurcation sets of Pham singularities}, Mosc. Math. J., 17:4 (2017), 825–836, \href{https://arxiv.org/abs/1701.03909}{arXiv:1701.03909}

\bibitem[Ve]{Ve} W. Veys, \textit{Introduction to the monodromy conjecture}, (2024), \href{https://arxiv.org/abs/2403.03343}{arXiv:2403.03343}

\bibitem[Y19]{Y19} A. Yuran, \textit{Newton Polytopes of Nondegenerate Quadratic Forms}, Functional Analysis and Its Applications 56 (2019): 152-158. \href{https://arxiv.org/abs/1910.06135}{arXiv:1910.06135}




        
\end{thebibliography}
\end{document}